%% file: paper.tex
\documentclass[a4paper,3p,10pt,times,preprint]{elsarticle}

\usepackage[utf8]{inputenc}
\usepackage{amsmath}
\usepackage{amsfonts}
\usepackage{amssymb}
\usepackage{amsthm}
\usepackage[dvipsnames,table]{xcolor}
\usepackage{graphicx}
\usepackage{adjustbox}
\usepackage{tikz,tikz-3dplot}
\usetikzlibrary{spy}
\usetikzlibrary{matrix}
\usetikzlibrary{arrows}
\usetikzlibrary{positioning}

\usepackage{booktabs}
\usepackage{multirow}
\usepackage{multicol}

\usepackage{pgfplots}
\usepackage{pgfplotstable}
\usetikzlibrary{pgfplots.groupplots}
\usetikzlibrary{arrows.meta,decorations.markings} 
\usepgfplotslibrary{fillbetween}
\usepackage{subcaption}

\usepackage{bm}
\usepackage{mathtools}
\usepackage{enumerate}
\usepackage{hyperref}

\hypersetup{
    colorlinks = true,
}
\usepackage[noabbrev]{cleveref}

\usepackage{physics}
\usepackage{bm}

\usepackage{algorithm}
\usepackage{algpseudocode}

\algrenewcommand\algorithmicrequire{\textbf{Input:}}
\algrenewcommand\algorithmicensure{\textbf{Output:}}

\newcommand{\R}{\mathbb{R}}

\usepackage{lineno}

\biboptions{sort&compress}

\input{TUDColors}

\pgfplotstableset{
    highlightrow/.append style={
        postproc cell content/.append code={
            \count0=\pgfplotstablerow
            \advance\count0 by1
            \ifnum\count0=#1
            \pgfkeysalso{@cell content/.add={$\it}{$}}
            \fi
        },
    },
    highlightcol/.append style={
        postproc cell content/.append code={
            \count0=\pgfplotstablecol
            \advance\count0 by1
            \ifnum\count0=#1
            \pgfkeysalso{@cell content/.add={$\it}{$}}
            \fi
        },
    },
}

\theoremstyle{plain}
\newtheorem{remark}{Remark}
\crefname{rem}{Remark}{Remarks}

\crefname{ex}{Example}{Examples}

\usepackage[compact]{titlesec}

\newcommand{\VEC}[1]{\vb*{#1}}

\newcommand{\RevTWO}[1]{{#1}}
\newcommand{\RevTHREE}[1]{{#1}}

\begin{document}

\begin{frontmatter}

\title{An Adaptive Parallel Arc-Length Method}

\author[label1,label2]{\corref{cor1}H.M. Verhelst}
\ead{h.m.verhelst@tudelft.nl}
\cortext[cor1]{Corresponding Author. }
\author[label1]{J.H. Den Besten}
\ead{henk.denbesten@tudelft.nl}
\author[label2]{M. Möller}
\ead{m.moller@tudelft.nl}
\address[label1]{Delft University of Technology, Department of Maritime and Transport Technology, Mekelweg 2, Delft 2628 CD, The Netherlands}
\address[label2]{Delft University of Technology, Department of Applied Mathematics, Van Mourik Broekmanweg 6, Delft 2628 XE, The Netherlands}

\begin{abstract}
Parallel computing is omnipresent in today’s scientific computer landscape, starting at multicore processors in desktop computers up to massively parallel clusters. While domain decomposition methods have a long tradition in computational mechanics to decompose spatial problems into multiple subproblems that can be solved in parallel, advancing solution schemes for dynamics or quasi-statics are inherently serial processes. For quasi-static simulations, however, there is no accumulating ‘time’ discretization error, hence an alternative approach is required. In this paper, we present an Adaptive Parallel Arc-Length Method (APALM). By using a domain parametrization of the arc-length instead of time, the multi-level error for the arc-length parametrization is formed by the load parameter and the solution norm. Given coarse approximations of arc-length intervals, finer corrections enable the parallelization of the presented method. This results in an arc-length method that is parallel within a branch and inherently adaptive. This concept is easily extended for bifurcation problems. The performance of the method is demonstrated using isogeometric Kirchhoff-Love shells on problems with snap-through and pitch-fork instabilities and applied to the problem of a snapping meta-material. These results show that parallel corrections are performed in a fraction of the time of the serial initialization, achievable on desktop scale.
\end{abstract}


\end{frontmatter}


\section{Introduction}
\input{Sections/Introduction.tex}
\clearpage

\section{Arc-Length methods}\label{sec:APALM_ALM}
\input{Sections/ArcLengthMethods.tex}

\section{Adaptive Parallel Arc-Length Method}\label{sec:APALM_APALM}
\input{Sections/APALM.tex}


\section{Implementation}\label{sec:APALM_implementation}

\input{Sections/Implementation.tex}

\section{Numerical Experiments}\label{sec:APALM_results}
\input{Sections/Results.tex}
\clearpage

\section{Conclusions and outlook}\label{sec:APALM_conclusions}
\input{Sections/Conclusions.tex}

\section*{Acknowledgments}
The authors are greatful for the financial support from Delft University of Technology.

\bibliography{references}

\end{document}

%% file: TUDColors.tex
\definecolor{tudelft-cyan}{cmyk}{1,0,0,0}
\definecolor{tudelft-black}{cmyk}{0,0,0,1}
\definecolor{tudelft-white}{cmyk}{0,0,0,0}
\definecolor{tudelft-dark-blue}{cmyk}{1,0.8,0.08,0.7}
\definecolor{tudelft-turquoise}{cmyk}{0.72,0,0.24,0}
\definecolor{tudelft-blue}{cmyk}{0.98,0.4,0.0,0}
\definecolor{tudelft-purple}{cmyk}{0.65,1,0,0.12}
\definecolor{tudelft-pink}{cmyk}{0,0.7,0,0}
\definecolor{tudelft-raspberry}{cmyk}{0.05,1,0.48,0.3}
\definecolor{tudelft-red}{cmyk}{0,0.85,0.75,0}
\definecolor{tudelft-orange}{cmyk}{0,0.7,0.75,0}
\definecolor{tudelft-yellow}{cmyk}{0,0.31,0.98,0}
\definecolor{tudelft-light-green}{cmyk}{0.63,0,0.84,0}
\definecolor{tudelft-dark-green}{cmyk}{1,0,0.68,0.04}

\newcommand{\normalplotheight}{0.25\textheight}

\colorlet{TUcol1}{tudelft-cyan}
\colorlet{TUcol2}{tudelft-black}
\colorlet{TUcol3}{tudelft-white}
\colorlet{TUcol4}{tudelft-dark-blue}
\colorlet{TUcol5}{tudelft-turquoise}
\colorlet{TUcol6}{tudelft-blue}
\colorlet{TUcol7}{tudelft-purple}
\colorlet{TUcol8}{tudelft-pink}
\colorlet{TUcol9}{tudelft-raspberry}
\colorlet{TUcol10}{tudelft-red}
\colorlet{TUcol11}{tudelft-orange}
\colorlet{TUcol12}{tudelft-yellow}
\colorlet{TUcol13}{tudelft-light-green}
\colorlet{TUcol14}{tudelft-dark-green}

\colorlet{gray1}{black!20}
\colorlet{gray2}{black!40}
\colorlet{gray3}{black!60}
\colorlet{gray4}{black!80}

\colorlet{col0}{TUcol4}
\colorlet{col1}{TUcol9}
\colorlet{col2}{TUcol12}
\colorlet{col3}{TUcol13}
\colorlet{col4}{TUcol6}
\colorlet{col5}{TUcol7}

\colorlet{introcol}{col0}
\colorlet{part1col}{col1}
\colorlet{part2col}{col2}
\colorlet{part3col}{col3}
\colorlet{part4col}{col4}
\colorlet{conclcol}{col0}

\colorlet{plotcol1}{TUcol10}
\colorlet{plotcol2}{TUcol12}
\colorlet{plotcol3}{TUcol13}
\colorlet{plotcol4}{TUcol14}
\colorlet{plotcol5}{TUcol6}
\colorlet{plotcol6}{TUcol7}

\pgfplotscreateplotcyclelist{mark list mod}{
every mark/.append style={solid,fill=\pgfplotsmarklistfill,line width = 0.1},mark size=1.5,mark=*,mark options=solid,fill opacity=0.5\\
every mark/.append style={solid,fill=\pgfplotsmarklistfill,line width = 0.1},mark size=1.5,mark=square*,mark options=solid,fill opacity=0.5\\
every mark/.append style={solid,fill=\pgfplotsmarklistfill,line width = 0.1},mark size=1.5,mark=triangle*,mark options=solid,fill opacity=0.5\\
every mark/.append style={solid,fill=\pgfplotsmarklistfill,line width = 0.1},mark size=1.5,mark=halfsquare*,mark options=solid,fill opacity=0.5\\
every mark/.append style={solid,fill=\pgfplotsmarklistfill,line width = 0.1},mark size=1.5,mark=pentagon*,mark options=solid,fill opacity=0.5\\
every mark/.append style={solid,fill=\pgfplotsmarklistfill,line width = 0.1},mark size=1.5,mark=halfcircle*,mark options=solid,fill opacity=0.5\\
every mark/.append style={solid,fill=\pgfplotsmarklistfill,rotate=180,line width = 0.1},mark size=1.5,mark=halfdiamond*,mark options=solid,fill opacity=0.5\\
every mark/.append style={solid,fill=\pgfplotsmarklistfill!40,line width = 0.1},mark size=1.5,mark=otimes*,mark options=solid,fill opacity=0.5\\
every mark/.append style={solid,fill=\pgfplotsmarklistfill,line width = 0.1},mark size=1.5,mark=diamond*,mark options=solid,fill opacity=0.5\\
every mark/.append style={solid,fill=\pgfplotsmarklistfill,line width = 0.1},mark size=1.5,mark=halfsquare right*,mark options=solid,fill opacity=0.5\\
every mark/.append style={solid,fill=\pgfplotsmarklistfill,line width = 0.1},mark size=1.5,mark=halfsquare left*,mark options=solid,fill opacity=0.5\\
}

\pgfplotscreateplotcyclelist{linestyles mod}{
thin,solid\\
thin,densely dashed\\
thin,densely dotted\\
thin,densely dashdotted\\
thin,densely dashdotdotted\\
thin,dashed\\
thin,dotted\\
thin,dashdotted\\
thin,dashdotdotted\\
}

\pgfplotscreateplotcyclelist{colorlist mod}{
plotcol1\\
plotcol2\\
plotcol3\\
plotcol4\\
plotcol5\\
plotcol6\\
}

\pgfplotsset{
	cycle multiindex* list={
							mark list mod\nextlist
							colorlist mod\nextlist
							linestyles mod\nextlist
							},
    width =\linewidth,
    height=\normalplotheight,
    grid=major,
    ylabel near ticks,
    xlabel near ticks,
    enlargelimits,
    grid=major,
	legend style={fill=white, fill opacity=0.6, draw opacity=1,draw=black,text opacity=1, font=\small},
    legend cell align={left},
    legend image with text/.style={
    legend image code/.code={%
        \node[anchor=center] at (0.3cm,0cm) {#1};
    }
	},
}

\pgfplotsset{select coords between index/.style 2 args={
		x filter/.code={
			\ifnum\coordindex<#1\fi
			\ifnum\coordindex>#2\fi
		}
}}

\tikzset{%
style0/.style = {plotcol1,thin,solid,every mark/.append style={solid,fill=\pgfplotsmarklistfill,line width = 0.1},mark size=1.5,mark=*,mark options=solid,fill opacity=0.75,},
style1/.style = {plotcol2,thin,densely dashed,every mark/.append style={solid,fill=\pgfplotsmarklistfill,line width = 0.1},mark size=1.5,mark=square*,mark options=solid,fill opacity=0.75},
style2/.style = {plotcol3,thin,densely dotted,every mark/.append style={solid,fill=\pgfplotsmarklistfill,line width = 0.1},mark size=1.5,mark=triangle*,mark options=solid,fill opacity=0.75},
style3/.style = {plotcol4,thin,densely dashdotted,every mark/.append style={solid,fill=\pgfplotsmarklistfill,line width = 0.1},mark size=1.5,mark=halfsquare*,mark options=solid,fill opacity=0.75},
style4/.style = {plotcol5,thin,densely dashdotdotted,every mark/.append style={solid,fill=\pgfplotsmarklistfill,line width = 0.1},mark size=1.5,mark=pentagon*,mark options=solid,fill opacity=0.75},
}

\tikzset{%
p2*/.style 	= {style=style0,fill opacity=0.75},
p2/.style 	= {style=style0,fill opacity=0},
p3*/.style 	= {style=style1,fill opacity=0.75},
p3/.style 	= {style=style1,fill opacity=0},
p4*/.style 	= {style=style2,fill opacity=0.75},
p4/.style 	= {style=style2,fill opacity=0},
p5*/.style 	= {style=style3,fill opacity=0.75},
p5/.style 	= {style=style3,fill opacity=0}
}

\tikzset{%
singlepatch/.style 	= {style=style0,thick,black,mark=0},
dpatch/.style 		= {style=style0},
dpatch2/.style 		= {style=style0,plotcol1!50!black},
approxC1/.style 	= {style=style1},
exactC1/.style 		= {style=style2},
almostC1/.style 	= {style=style3},
weak1/.style 		= {style=style4,plotcol5!80!black},
weak2/.style 		= {style=style4,plotcol5!60!black},
weak3/.style 		= {style=style4,plotcol5!40!black},
weak4/.style 		= {style=style4,plotcol5!20!black}
}

\pgfplotstableset{
	highlightrow/.append style={
		postproc cell content/.append code={
			\count0=\pgfplotstablerow
			\advance\count0 by1
			\ifnum\count0=#1
			\pgfkeysalso{@cell content/.add={$\it}{$}}
			\fi
		},
	},
	highlightcol/.append style={
		postproc cell content/.append code={
			\count0=\pgfplotstablecol
			\advance\count0 by1
			\ifnum\count0=#1
			\pgfkeysalso{@cell content/.add={$\it}{$}}
			\fi
		},
	},
}

\pgfdeclarelayer{bg}    
\pgfsetlayers{bg,main}  

\pgfkeys{
    /pgfplots/table/print last/.style={
        row predicate/.code={
            \pgfmathtruncatemacro\firstprintedrownumber{\pgfplotstablerows-#1}
            \ifnum##1<\firstprintedrownumber\relax
            \pgfplotstableuserowfalse
            \fi
        }
    },
    /pgfplots/table/print last/.default=1
}


%% file: Sections/Introduction.tex
Over the last decades, computational power has increased exponentially. In the last year, most improvements were due to an increasing number of threads per processing unit rather than an increase in single-thread performance \cite{Rupp2022}. The trend of increasing logical cores with stagnating single-threaded performance calls for parallelization of existing codes to improve computational efficiency, amongst which numerical algorithms in computational mechanics. In the field of computational mechanics, parallelization in the spatial domain is common practice by using shared-memory assembly routines or distributed-memory parallelization using domain decomposition of meshes. Parallelization can also be achieved in linear solvers or in the temporal domain using parallel-in-time solvers \cite{Gander2015} in the case of dynamic analyses or using parallel continuation for quasi-static or continuation problems - the latter two being sequential by nature.\\


For quasi-static problems, continuation methods can be used when the solution of an equation or a system of equations is desired, given the varying parameters of the system. Such methods, typically referred to as Arc-Length Methods (ALMs), are widely used for (but not limited to) the analysis of the stability of structures. The Riks and Crisfield methods \cite{Riks1972,Crisfield1981} are commonly used and in combination with bifurcation algorithms \cite{Wriggers1988}, whereas ALMs provide a valuable tool in the analysis of the collapse and post-buckling behaviour of structures. Recent developments for ALMs include a new displacement-controlled formulation \cite{Pretti2022}, an improved predictor scheme \cite{Kadapa2021}, and automatic exploration techniques \cite{Thies2021,Wouters2019}. Like time-stepping methods, ALMs are sequential by nature, meaning that the solution at a point is obtained from the solution at a previous point obtained previously.\\

Amongst many parallel time-integration schemes, Parareal is a parallel time-integration method proposed by \cite{Lions2001} and works with a two-level parallel correction scheme of time intervals. The method starts with a series of solutions obtained in serial with a large time step, after which each sub-interval is computed with a finer time step such that a new solution is found at the end-point of the time interval. A multi-level extension of Parareal is proposed in \cite{Falgout2014} and is referred to as Multi-Grid Reduced in Time (MGRIT). This method is similar to Parareal but applies the two-level approach recursively. As a consequence, multi-grid-like cycles can be used to correct previously computed sub-intervals. This method has not only been applied to dynamic problems but also to the training of neural networks \cite{cyr2019multilevel} and \cite{hessenthaler2021time}. Alternative methods for parallel time integration are reviewed in the work of \cite{Gander2015}.\\

Compared to temporal parallelization methods, parallelization of ALMs has received less attention in the academic community. As ALMs are typically used for explorations of solutions across branches, parallel evaluation of branches can be performed as soon as the starting point (and tangent) of each branch is known. The number of branches related to a problem, however, depends typically on the problem that is solved; hence, the parallel scalability of ALMs over branches is not guaranteed. Parallelization within a branch is enabled by the Parallel Adaptive Method for Pseudo-Arclength Continuation (PAMPAC) \cite{Aruliah2016}. This method works with multiple predictors (with different step sizes) and consequently correctors to select an optimal step size, which can be performed in parallel. The PAMPAC method focuses on selecting a maximal step size for the ALM for which the method does not converge.\\


In this paper, a parallelization of the arc-length method is presented that is independent of the physical nature of the underlying problem. That is, the method is developed such that the parallelization can be performed within the branches. In addition to parallelization, the presented arc-length scheme also provides inherent adaptivity; therefore, the method is referred to as the Adaptive Parallel Arc Length Method (APALM). The working principle of the APALM is based on a multi-level approach -- inspired by MGRIT methods -- where a coarse serial approximation of the solution space is refined in parallel until a measure of convergence is achieved. Contrary to PAMPAC, the present method does not maximise the step size for convergence of the ALM iterations, but instead the parallelization is based on convergence of the solution sub-intervals. Without loss of generality, the method is developed given a constraint equation for the arc-length method; thus, it is generalised for the Riks and Crisfield methods, amongst other methods available.\\

The outline of this paper is as follows: \Cref{sec:APALM_ALM} provides a background on arc-length methods. In \cref{sec:APALM_APALM}, the parallelization of arc-length methods is presented, referred to as the APALM. Thereafter, \cref{sec:APALM_implementation} provides algorithms for non-intrusive implementation of the APALM, given an implementation of an existing ALM. \Cref{sec:APALM_results} provides numerical benchmark problems and an application to the analysis of a snapping meta-material, inspired by \cite{Rafsanjani2015}. Finally, \cref{sec:APALM_conclusions} provides conclusions on the presented method.

%% file: Sections/ArcLengthMethods.tex
In this section, the concept of arc-length methods is presented for the sake of completeness. For a detailed overview, one can consult references \cite{Ritto-Correa2008,Ragon2002,Crisfield1981,Riks1972}. Let $\vb*{G}(\vb*{u},\lambda)=\vb*{0}$ be a non-linear system of equations to be solved, with $\vb*{u}$ the solution to the system of equations given a parameter $\lambda$. For structural analyses, $\vb*{u}$ is typically a vector containing discrete displacements of the degrees of freedom, and $\lambda$ is a factor scaling the magnitude of an applied load $P$, i.e.
\begin{equation}\label{eq:APALM_nl}
  \vb*{G}(\vb*{u},\lambda) = \vb*{N}(\vb*{u}) - \lambda\vb*{P},
\end{equation}
where $\vb*{N{(\vb*{u})}}$ is a vector of internal forces, depending on the deformation $\vb*{u}$. For incremental analyses, i.e., quasi-static analyses, a series of solutions $\vb*{w}_i=(\vb*{u}_i,\lambda_i)$ is obtained by computing increments $\Delta\vb*{w}_i=(\Delta\vb*{u}_i,\Delta\lambda_i)$ such that $\vb*{w}_{i+1}=\vb*{w}_i+\Delta\vb*{w}_i$ and \cref{eq:APALM_nl} is satisfied for $\vb*{w}_{i+1}$. These solutions can be obtained by Newton iterations: i) fixing $\lambda$ and finding $\vb*{u}$ (\emph{load control}); ii) fixing some degrees of freedom in $\vb*{u}$ and finding all $\vb*{u}$ and $\lambda$ (\emph{displacement control}); or iii) constraining $\lambda$ and $\vb*{u}$ and solving for both (\emph{arc-length control}); see \cref{fig:control}. In the case of arc-length control, the increment $\Delta\vb*{w}$ is measured by an increment length $d(\Delta\vb*{w})$
\begin{equation}\label{eq:APALM_distance}
  d(\Delta\vb*{w}) = \Delta\vb*{u}^\top\Delta\vb*{u} + \Psi^2\Delta\lambda^2\vb*{P}^\top\vb*{P},
\end{equation}
where $\Psi$ is a scaling parameter given in \cite{Schweizerhof1986, Bellini1987}. The increment $\Delta\vb*{w}$ is constrained by the arc-length $\Delta \ell$ in the constraint equation
\begin{equation}
  f(\Delta\vb*{w}) = d(\Delta\vb*{w})-\Delta\ell=0.
\end{equation}
Since $\VEC{G}(\vb*{u},\lambda)$ is non-linear, the increment $\Delta\vb*{w}_i$ is obtained iteratively, i.e., $\Delta\vb*{w}_{i,k+1}=\Delta\vb*{w}_{i,k} + \delta\vb*{w}_i$ with iteration count $k$. The constraint equation is solved together with \cref{eq:APALM_nl} in every iteration, yielding the Riks and Crisfield methods \cite{Crisfield1981,Riks1972}
\begin{align}
  f(\Delta\vb*{w}_{i,k},\Delta l) &= \Delta\vb*{u}_{i,0}^\top\Delta\vb*{u}_{i,k} + \Psi^2\Delta\lambda_{i,0}\Delta\lambda_{i,0}\vb*{P}^\top\vb*{P}-\Delta \ell^2=0, &&\text{Riks},\\
  f(\Delta\vb*{w}_{i,k},\Delta l) &= \Delta\vb*{u}_{i,k}^\top\Delta\vb*{u}_{i,k} + \Psi^2\Delta\lambda_{i,k}\Delta\lambda_{i,k}\vb*{P}^\top\vb*{P}-\Delta \ell^2=0, &&\text{Crisfield},
\end{align}
where $\Delta\vb*{w}_0$ is the increment in the first iteration. The Crisfield method generally performs well with sharp snap-backs but has the disadvantage that the constraint equation has two intersections with the path formed by \cref{eq:APALM_nl}. Hence, a root has to be selected, which is elaborated in the works \cite{Crisfield1981,Ritto-Correa2008}. When multiple intersections are found, complex roots are found \cite{Carrera1994}, which can be resolved using one of the methods proposed in \cite{Lam1992,Zhou1995}. \RevTWO{It should be noted that any other arc-length method can be used within the scheme proposed in this paper, as long as the constraint equation is satisfied when the arc-length step is converged.}

\begin{figure}
  \centering
  \begin{subfigure}{0.3\linewidth}
    \centering
    \input{Figures/ALM_loadcontrol.pgf}
    \caption{Load control}
  \end{subfigure}
  \hfill
  \begin{subfigure}{0.3\linewidth}
    \centering
    \input{Figures/ALM_displacementcontrol.pgf}
    \caption{Displacement control}
  \end{subfigure}
  \hfill
  \begin{subfigure}{0.3\linewidth}
    \centering
    \input{Figures/ALM_arclengthcontrol.pgf}
    \caption{Arc-length control}
  \end{subfigure}
  \caption{Load (left), displacement (middle), and arc-length control (right) for structural analysis problems. The question mark (?) indicates the iteration where load and displacement control encounter a limit point. In these situations, the next point obtained is typically difficult to find. }
  \label{fig:control}
\end{figure}

%% file: Figures/ALM_loadcontrol.pgf
\begin{tikzpicture}[scale=0.8]
\coordinate (O) at (0,0);
\coordinate (X) at (10,0);
\coordinate (Y) at (0,5);
\draw[-latex] (O) -- (5,0) node [midway,below] {Load $\lambda$};
\draw[-latex] (0,0) -- (0,4.0) node [midway,above ,rotate=90]{$\Vert \vb*{u}(t)\Vert$};

\draw plot[smooth,tension=1] coordinates {(0,1.5) (2,2.5) (2.5,1.5) (1.5,0.5) (3,0.5) (4.5,1.5)};

\draw[-latex](4.0,2.0)node[above,scale=1]{$G(\vb*{u},\lambda)=0$}
        to[out=-60,in=120] (4.2,1.2);

\foreach \x in {0,0.4,...,2.4}{
\draw[dotted](\x,0) -- (\x,4);
}
\draw[dashed](2.8,0) -- node[right]{?}(2.8,4);

\coordinate (p0) at (0.0,1.5);
\coordinate (p1) at (0.4,1.77);
\coordinate (p2) at (0.8,2.03);
\coordinate (p3) at (1.2,2.24);
\coordinate (p4) at (1.6,2.4);
\coordinate (p5) at (2.0,2.5);
\coordinate (p6) at (2.4,2.4);

\foreach \k in {0,1,...,6}{\node at (p\k){$\circ$};}
\end{tikzpicture}

%% file: Figures/ALM_displacementcontrol.pgf
\begin{tikzpicture}[scale=0.8]
\coordinate (O) at (0,0);
\coordinate (X) at (10,0);
\coordinate (Y) at (0,5);
\draw[-latex] (O) -- (5,0) node [midway,below] {Load $\lambda$};
\draw[-latex] (0,0) -- (0,4.0) node [midway,above ,rotate=90]{$\Vert \vb*{u}(t)\Vert$};

\draw plot[smooth,tension=1] coordinates {(0,1.5) (2,2.5) (2.5,1.5) (1.5,0.5) (3,0.5) (4.5,1.5)};

\draw[-latex](4.0,2.0)node[above,scale=1]{$G(\vb*{u},\lambda)=0$}
        to[out=-60,in=120] (4.2,1.2);

\foreach \y in {1.5,1.8,...,2.41}{
\draw[dotted](0,\y) -- (4,\y);
}
\draw[dashed](0,2.7) -- node[above]{?}(4,2.7);

\coordinate (p0) at (0.0,1.5);
\coordinate (p1) at (0.46,1.8);
\coordinate (p2) at (0.94,2.1);
\coordinate (p3) at (1.572,2.4);

\foreach \k in {0,1,...,3}{\node at (p\k){$\circ$};}
\end{tikzpicture}

%% file: Figures/ALM_arclengthcontrol.pgf
\begin{tikzpicture}[scale=0.8]
\coordinate (O) at (0,0);
\coordinate (X) at (10,0);
\coordinate (Y) at (0,5);
\draw[-latex] (O) -- (5,0) node [midway,below] {Load $\lambda$};
\draw[-latex] (0,0) -- (0,4.0) node [midway,above ,rotate=90]{$\Vert \vb*{u}(t)\Vert$};

\draw plot[smooth,tension=1] coordinates {(0,1.5) (2,2.5) (2.5,1.5) (1.5,0.5) (3,0.5) (4.5,1.5)};

\draw[-latex](4.0,2.0)node[above,scale=1]{$G(\vb*{u},\lambda)=0$}
        to[out=-60,in=120] (4.2,1.2);
\draw[-latex](0.75,2.5)node[above,scale=1]{$f(\Delta\vb*{u},\Delta\lambda)=0$} to[out=-90,in=60] (0.5,1.87);

\coordinate (p0) at (0,1.5);
\coordinate (p1) at (0.83,2.04);
\coordinate (p2) at (1.75,2.45);
\coordinate (p3) at (2.57,1.865);
\coordinate (p4) at (2.07,0.99);

\node at (0,1.5) {$\circ$};
\draw[color=black!20,thin,densely dashed] (0,1.5) circle [radius=1];
\node at (p1) {$\circ$};
\draw[color=black!40,thin,densely dashed] (p1) circle [radius=1];
\node at (p2) {$\circ$};
\draw[color=black!60,thin,densely dashed] (p2) circle [radius=1];
\node at (p3) {$\circ$};
\draw[color=black!80,thin,densely dashed] (p3) circle [radius=1];
\node at (p4) {$\circ$};
\draw[color=black!100,thin] (p4) circle [radius=1];
\end{tikzpicture}

%% file: Sections/APALM.tex
In this section, our new approach, the APALM, is presented. Firstly, the method is conceptualised along with some illustrative figures (\cref{subsec:APALM_concept}). Secondly, details are provided on the curve parameterization and the measurement of errors (\cref{subsec:APALM_errors}). Lastly, \cref{subsec:APALM_parameterization} presents (re-)parameterization methods for the solution curve. These parameterizations will be essential to the data structure of the APALM. It should be noted that the method described in this section is presented only for one continuation parameter, $\lambda$.

\subsection{Concept}\label{subsec:APALM_concept}
Learning from parallel-in-time methods like Parareal or MGRIT, parallelization in the APALM is achieved from a subdivision of the \emph{curve length domain}. Contrary to MGRIT and Parareal, where the temporal domain $t\in[T_0,T_1]$ is fixed, the APALM will work with a changing \emph{curve length domain} $s\in[S_0,S_1]$ depending on the length of the traversed path, with an underlying fixed parametric domain with parametric coordinate $\xi\in[0,1]$. The APALM is initialised with an initial coarse grid approximation, in which the parametric and the curve length domains are subdivided into sub-domains $\xi\in[\xi_i,\xi_{i+1}]$ and $s\in[s_i,s_{i+1}]$, respectively, as illustrated in \cref{fig:APALM_APALM_concept}.\\

In the initialization phase of the APALM, the first subdivision into sub-intervals is made (see \cref{fig:APALM_APALM_initialization}). Here, the sizes of the sub-intervals $s\in[s_i^\ell,s^\ell_{i+1}]$ are determined based on the distance measure that is used by the corresponding ALM; see \cref{eq:APALM_distance}. Note that the superscript $\ell$ denotes the $\ell^{\text{th}}$ level. Based on the initial curve-length domain $s\in[0,S]$, where $S$ is the total length of the initial curve, and the corresponding sub-intervals, the curve-length domain can be mapped accordingly onto a parametric domain; see \cref{sec:APALM_implementation} for more details.\\

With an initialised computational domain, the number of sub-intervals determines the degree of parallelization. On any sub-interval, $[s^\ell_i,s^\ell_{i+1}]$ data at the start-point and end-point is known, which can be used to initiate an arc-length method to re-compute the sub-interval with $N$ increments, i.e., with an arc-length of $\Delta L^{\ell+1}_i = \Delta L^\ell_i/N$ (see \cref{fig:APALM_APALM_computation}).\\

After sub-interval $[s^{\ell+1}_i,s^{\ell+1}_{i+1}]$ has been finished, the distance of the end-point of the sub-interval can be compared to the previously known solution at $s^\ell_{i+1}$, which is called \emph{parallel verification of intervals} in \cref{fig:APALM_APALM_verification}. Since the sub-interval is traversed in $N$ increments with length $\Delta s^0_i/N$, the triangle inequality with the arc-length measure implies that there must be a distance greater than or equal to zero between the newly found end-point and the reference end-point. The more `curved' the domain in-between, the larger this distance. Based on an error measure (see \cref{subsec:APALM_errors}), intervals with a relatively large deviation between the coarse-level arc length and the fine-level arc length are to be marked for `refinement'.\\

Lastly, the intervals with a too large deviation in the newly computed curve length need to be reparameterized (see \cref{subsec:APALM_parameterization}). This is because the total curve-length parameterization is elongated exactly by the distance between the newly computed endpoint and the previously known point. By this means, the reference interval is subdivided into $N+1$ sub-intervals, and the data corresponding to the $N$ newly computed point is stored. For sub-intervals that have an error below the tolerance, only $N-1$ points are stored as references, and no reparameterization takes place. The process is sketched in \cref{fig:APALM_APALM_reparameterization}. After reparameterization, the marked interval can be re-computed, and the process can be repeated from \cref{fig:APALM_APALM_computation} onward.\\

\begin{remark}[Difference with parallel-in-time methods]
As mentioned, the multi-level approach that is employed in this method is derived from the idea of parallel-in-time methods. However, the fundamental difference between time integration and continuation comes from the fact that time integration methods typically compute the solution on the next time step with a certain time integration error $\mathcal{O}(\Delta t^p)$. Parallel-in-time methods rely on this time integration error to mark solution intervals as converged or not, and additionally, updated solutions contain smaller time integration errors, so sub-intervals need to be recomputed as soon as solutions previously in time have been updated.\\

For arc-length methods, Newton's method is applied to a system of equations that solves the arc-length constraint equation together with the discretized system $\VEC{G}(\VEC{u},\lambda)=\VEC{0}$. Therefore, the error of the solution that is found after an arc-length increment is independent of the arc-length increment size but solely depends on the convergence tolerance of Newton's method. Therefore, the end-point of an interval does not have to be updated, nor do intervals after the update be recomputed. \RevTHREE{This implies, in principle, that parallel corrections of the arc-length steps are not needed, since the intervals already capture the structural response at the equilibrium path. However, the parallel corrections are still meaningful to capture the equilibrium path in desired detail, in the case where the initial step size is chosen very coarse. As will follow from the results section, cf. \cref{sec:APALM_results}, the parallel performance increases for fewer (hence coarser) initial intervals.}\\
\end{remark}

\RevTHREE{
\begin{remark}[Path-dependency]
The concept presented in this paper assumes path-independence of the equation to be solved, in order to assume that from a given starting point on the equilibrium path, the same end-point could be reached irrespective of the computed intervals on the path in-between the points. Path-dependent problems are out of the scope of this paper.
\end{remark}
}

\begin{figure}
    \centering
    \begin{subfigure}[t]{0.45\linewidth}
        \begin{adjustbox}{width=\linewidth}
            \input{Figures/APALM_it0.pgf}
        \end{adjustbox}
        \caption{\textbf{Initialization based on computed reference solutions.} Without losing generality, the solutions are separated by a fixed distance, $\Delta s$. Given their distances, an initial estimation of the \emph{curve length $s$} can be produced, which can be mapped on the \emph{parametric domain $\xi$}.}
        \label{fig:APALM_APALM_initialization}
    \end{subfigure}
    \hfill
    \begin{subfigure}[t]{0.45\linewidth}
        \begin{adjustbox}{width=\linewidth}
            \input{Figures/APALM_it1.pgf}
        \end{adjustbox}
        \caption{\textbf{Parallel computation of intervals.} On each interval, a finer estimate can be performed by splitting the interval into $n$ sub-intervals ($n=2$ here).}
        \label{fig:APALM_APALM_computation}
    \end{subfigure}

    \begin{subfigure}[t]{0.45\linewidth}
        \begin{adjustbox}{width=\linewidth}
            \input{Figures/APALM_it2.pgf}
        \end{adjustbox}
        \caption{\textbf{Parallel verification of intervals.} When the last sub-interval is computed, the solution is verified with the next known reference solution (here, the solutions following from the initial simulations). When the distance is sufficiently small, the segment can be marked as convergent, and the in-between solutions on the interval can be written.}
        \label{fig:APALM_APALM_verification}
    \end{subfigure}
    \hfill
    \begin{subfigure}[t]{0.45\linewidth}
        \begin{adjustbox}{width=\linewidth}
          \input{Figures/APALM_it3.pgf}
        \end{adjustbox}
        \caption{\textbf{Curve-length reparameterization.} For sub-intervals where the length deviates too much from the previous length, i.e., where the distance between the last computed point and the reference is too large, all solutions of the sub-interval are added to the parameterization, and the known data points ahead of the newly computed points are shifted in the curve parametric coordinates.}
        \label{fig:APALM_APALM_reparameterization}
    \end{subfigure}
    \caption{Concept of the APALM. The large open circles represent \emph{reference solutions} from a previously computed level. The small solid circles represent \emph{new data} on the interval between two reference solutions, computed by the arc-length method (here the large dashed circle). The black dashed line indicates the curve estimation for which the sum is equal to the total curve length.}
    \label{fig:APALM_APALM_concept}
\end{figure}

\subsection{Error Measures}\label{subsec:APALM_errors}
The refinement of computed sub-intervals depends on the distances between the points in the original (coarse) interval and the newly obtained solutions in this sub-interval. Here, error measures are presented, that can be used to mark an interval $[s^\ell_i,s^\ell_{i+1}]$ based on the obtained solutions $\{s^{\ell+1}_k\}_{k=0,...,N}$ at the finer level. \Cref{fig:APALM_distances} presents two possible situations: a nearly straight interval that would not be marked for refinement, and a curved interval that would be marked for refinement. Here, the interval is considered `curved' in the discrete solution space if the hyperdimensional path between two solutions differs from the hyperplane between these solutions. The errors that determine the marking of an interval for refinement are illustrated in \cref{fig:APALM_distance_irregular} and can be interpreted as follows: $\Delta L$ is the original arc length between two coarse solutions, $\Delta L'$ is the newly obtained length between two coarse solutions, the \emph{lower distance} $\overline{\Delta L}$ is the distance between the start of the interval and the last solution at the fine level, and $\delta L$ is the distance between the last obtained solution on the fine level and the final point on the coarse level. Using these distances, the \emph{total error} ($\varepsilon$), the \emph{lower error} ($\varepsilon_l$), and the \emph{upper error} ($\varepsilon_u$) can be defined:
\begin{align}
\varepsilon &= (\Delta L' - \Delta L) / \Delta L, && \quad \text{total error},\label{eq:APALM_total_error}\\
\varepsilon_l &= (\Delta L - \overline{\Delta L}) / \Delta L, && \quad \text{lower error},\label{eq:APALM_lower_error}\\
\varepsilon_u &= (\varepsilon - \varepsilon_l) / \Delta L, && \quad \text{upper error}.\label{eq:APALM_upper_error}
\end{align}
Here, the total error is the total difference between the coarse and fine intervals; the lower error is the contribution of the first $N$ sub-intervals; and the upper error is the contribution of $\delta L$ to the total error. Depending on these errors and specified tolerances, refinement rules can be set up, in particular:
\begin{align}
\text{Refine the first $N$ intervals } \longleftrightarrow &\quad \varepsilon_l > \text{TOL}_l,\\
\text{Refine the last interval } \longleftrightarrow &\quad \varepsilon_u > \text{TOL}_u.
\end{align}

\begin{figure}
    \centering
    \begin{subfigure}[t]{0.45\linewidth}
        \begin{adjustbox}{width=\linewidth}
            \input{Figures/Lengths_regular.pgf}
        \end{adjustbox}
        \caption{Nearly straight interval where the small white dot is sufficiently close to the coarse reference. In this case, the solution of the second increment is not added to the parameterization.}
        \label{fig:APALM_distance_regular}
    \end{subfigure}
    \hfill
    \begin{subfigure}[t]{0.45\linewidth}
        \begin{adjustbox}{width=\linewidth}
            \input{Figures/Lengths_irregular.pgf}
        \end{adjustbox}
        \caption{Curved interval where the solution of the second increment is $\Delta L_u$ away from the coarse reference. In this case, the distance between the coarse reference points $\delta L$ is smaller than the actually traversed distance $\Delta L'$.}
        \label{fig:APALM_distance_irregular}
    \end{subfigure}
    \caption{Error measures on a nearly straight interval (a) and a curved interval (b). For the nearly straight interval, the distance $\delta L$ (see b) is sufficiently small, whereas for the curved interval, it is too big. The measures $\Delta L$ and $\Delta L'$, $\overline{\Delta L}$, and $\delta L$ are, respectively, the coarse arc length, the fine arc length, the lower distance, and the absolute error.}
    \label{fig:APALM_distances}
\end{figure}

\subsection{Curve (Re-)Parameterization}\label{subsec:APALM_parameterization}
As indicated in \cref{fig:APALM_APALM_concept}, the concept of the APALM is supported by the parameterization of the solutions of $G(\VEC{u},\lambda)=\VEC{0}$ by parameterizing the curve length using the increment length $d(\Delta\VEC{w})$ embedded in the arc-length method. As illustrated in \cref{fig:APALM_APALM_concept}, the APALM maps solutions $\VEC{w}$ to a point on the curve-length domain $[0,S]$, and points on the curve-length domain are mapped on a parametric domain $[0,1]$.

Provided a series of solutions from the initialization phase $\{\VEC{w}^0_i\}_{i=0,...,I}$, with $I$ denoting the total number of initial points, and defining solution intervals by $\Delta\VEC{w}^\ell_i=\VEC{w}^\ell_{i+1}-\VEC{w}^\ell_i$, each solution $\VEC{w}^0_i$ can recursively be assigned to the curve-length and parametric domains by
\begin{align}
    s^0_{i+1} &= s^0_{i} + d(\Delta\VEC{w}^0_i),&&\quad i=1,...,I-1, \quad s_0=0,\label{eq:APALM_curve_parameterization}\\
    \xi^0_{i} &= \frac{s_{i}}{s_I},&&\quad i=0,...,I,\label{eq:APALM_parameterization}
\end{align}
where the superscript $0$ represents the $0^\text{th}$ level. In addition, \cref{eq:APALM_curve_parameterization} guarantees that $S=s_I$ marks the total length of the curve that has been traversed, measured by the distance between each solution. Given the curve-length coordinates of each point as an increasing sequence, the parametric domain can simply be obtained by scaling the domain back to $[0,1]$; see \cref{eq:APALM_parameterization}. In the following, two ways of adding solutions to the parameterization are defined: i) interior insertion, and ii) full insertion and stretching. The operations are defined given a parent interval $[s^\ell_i,s^\ell_{i+1})$ in which a set of new solutions $\{s^{\ell+1}_k\}_{k=0,...,N}$, where $s^{\ell+1}_0=s^{\ell}_i$, is computed, with $N$ the total number of points in the interval; see \cref{fig:APALM_domain_ori}. \\

Firstly, the interior insertion operation inserts solutions \emph{within} the sub-interval, see \cref{fig:APALM_domain_interior}, and is later used for intervals where the error is small. The idea behind this operation is that the solutions $\{s^{\ell+1}_k\}_{k=1,...,N-1}$ between $s^\ell_i$ and $s^\ell_{i+1}$ are inserted and that the solution $s^{\ell+1}_N$ is not added to the map. In the case of the interior insertion, the points $s^{\ell+1}_k$ and their parametric coordinates $\xi^{\ell+1}_k$ are added by:
\begin{align}
    s^{\ell+1}_{k+1} &= s^{\ell+1}_{k} + d(\Delta\VEC{w}^{\ell+1}_{k}),&&\quad k=0,...,N-2, \quad s^{\ell+1}_0=s^{\ell}_i, \quad s^{\ell+1}_N=s^{\ell}_{i+1},\\
    \xi^{\ell+1}_{k+1} &= \xi^{\ell+1}_{k} + (\xi^\ell_{i+1}-\xi^\ell_{i})\frac{s^{\ell+1}_{k+1}-s^{\ell+1}_{k}}{s^{\ell}_{i+1}-s^{\ell}_{i}},&&\quad j=0,...,N-2.
\end{align}
Note that $\VEC{w}^{\ell+1}_{k}$ denotes the $k^{\text{th}}$ solution on level $\ell+1$ on the computed sub-interval, here $[s^\ell_i,s^\ell_{i+1}]$.\\

The full insertion and stretching operation inserts the solutions of the sub-interval, including its end point, and also stretches the curve parameterization (see \cref{fig:APALM_domain_stretch}), which is later used for intervals where the error is large, hence intervals that need refinement. The idea behind this operation is that the solutions $\{s^{\ell+1}_k\}_{k=1,...,N-1}$ between $s^\ell_i$ and $s^\ell_{i+1}$ are inserted and that the point $s^\ell_{i+1}$ is shifted such that $s^\ell_{i+1} = s^{\ell+1}_N$ and such that all points further than $s^\ell_{i+1}$ are updated to $\tilde{s^{\ell}_{j}}$ by the a shift using the distance between the last computed solution and the reference solution, i.e. $d(\Delta\VEC{w}^{\ell+1}_N)$, $\Delta\VEC{w}^{\ell+1}_N = \VEC{w}^{\ell}_{i+1} - \VEC{w}^{\ell+1}_N$:
\begin{align}
    s^{\ell+1}_{k+1} &= s^{\ell+1}_{k} + d(\Delta\VEC{w}^{\ell+1}_{k}),&&\quad k=0,...,N-1, \quad s^{\ell+1}_0=s^{\ell}_i, \quad s^{\ell}_{i+1} = s^{\ell+1}_N,\\
    \xi^{\ell+1}_{k+1} &= \xi^{\ell+1}_{k} + (\xi^\ell_{i+1}-\xi^\ell_{i})\frac{s^{\ell+1}_{k+1}-s^{\ell+1}_{k}}{s^{\ell}_{i+1}-s^{\ell}_{i}},&&\quad k=0,...,N-1,\label{eq:APALM_chord-length}\\
    \tilde{s}^{\ell}_{j} &= s^{\ell}_{j} + d(\VEC{w}^{\ell+1}_N,\VEC{w}^{\ell}_{i+1}), &&\quad j=i+1,...,I.
\end{align}

As can be noticed in \cref{eq:APALM_chord-length}, the re-scaling of $\xi$ is done using the parametric length of the original interval at level $\ell$, $(\xi^{\ell}_{i+1}-\xi^{\ell}_{i})$ and the curve coordinate $s_{i+1}$ relative to the beging point of the interval $s_{i}$ with respect to the (updated) total curve length of the interval $s_{i+1}-s_{i}$, which is similar to the well-known \emph{chord-length parameterization} in splines \cite{Piegl1995}.

\begin{figure}
\centering
\begin{subfigure}{\linewidth}
\centering
\input{Figures/domain_ori.pgf}
\caption{Original domain}
\label{fig:APALM_domain_ori}
\end{subfigure}

\begin{subfigure}{\linewidth}
\centering
\input{Figures/domain_interior.pgf}
\caption{Interior insertion}
\label{fig:APALM_domain_interior}
\end{subfigure}

\begin{subfigure}{\linewidth}
\centering
\input{Figures/domain_stretch.pgf}
\caption{Full insertion and stretching}
\label{fig:APALM_domain_stretch}
\end{subfigure}
\caption{Domain parameterizations on the curve-length domain $s$ and the parameter domain $\xi$ with levels $\ell$ and $\ell+1$. \Cref{fig:APALM_domain_ori} illustrates the original domain, \cref{fig:APALM_domain_interior} illustrates the insertion of interior points in the case of a sufficiently close approximation of the end-point of the domain and \cref{fig:APALM_domain_stretch} illustrates the full insertion of all sub-domain solutions combined with the stretching of the curve length domain.}
\label{fig:APALM_domain}
\end{figure}

%% file: Figures/APALM_it0.pgf
\begin{tikzpicture}[scale=1.5,every node/.style={inner sep=10}]
\coordinate (O) at (0,0);
\coordinate (X) at (10,0);
\coordinate (Y) at (0,5);
\draw[-latex] (0,0) -- (5,0) node [right] {$\lambda$};
\draw[-latex] (0,0) -- (0,3.0) node [midway,above ,rotate=90]{$\Vert \mathbf{u}(t)\Vert$};

\draw[-latex] (0,-0.8) -- (5,-0.8)node [right] {$s$};
\draw[-latex] (0,-2) -- (5,-2)node [right] {$\xi$};
\draw (0,-2.1) -- (0,-1.9) node(xi0)[midway,inner sep=0]{};
\draw (4.5,-2.1) -- (4.5,-1.9) node(xi1)[midway,inner sep=0]{};
\node[above] at (xi0) {0};
\node[above] at (xi1) {1};

\draw[black!20] plot[smooth,tension=1] coordinates {(0,1) (2,2) (3,1)  (4.5,2)};

\draw[-latex](4.0,2.5)node[above,scale=1]{$G(\mathbf{u},\lambda)=0$} to[out=-90,in=120] (4.2,1.69);

\coordinate (p0) at (0,1);
\coordinate (p1) at (0.82,1.58);
\coordinate (p2) at (1.74,1.97);
\coordinate (p3) at (2.58,1.39);
\coordinate (p4) at (3.55,1.17);

\draw (p0) circle [radius=0.1];
\draw (p1) circle [radius=0.1];
\draw (p2) circle [radius=0.1];
\draw (p3) circle [radius=0.1];

\draw[fill=white] let \p1 = (p0) in (\x1,-0.8) circle [radius=0.1];
\draw[fill=white] let \p1 = (p1) in (\x1+5,-0.8) circle [radius=0.1];
\draw[fill=white] let \p1 = (p2) in (\x1+7,-0.8) circle [radius=0.1];
\draw[fill=white] let \p1 = (p3) in (\x1+14,-0.8) circle [radius=0.1];
\draw[fill=white] let \p1 = (p4) in (\x1+16,-0.8) circle [radius=0.1];

\draw (1.125,-2.1) -- (1.125,-1.9) node(xi025)[midway,inner sep=0]{};
\draw (2.25,-2.1) -- (2.25,-1.9) node(xi050)[midway,inner sep=0]{};
\draw (3.375,-2.1) -- (3.375,-1.9) node(xi075)[midway,inner sep=0]{};
\node[above] at (xi025) {0.25};
\node[above] at (xi050) {0.50};
\node[above] at (xi075) {0.75};

\draw[dotted] let \p1 = (p0) in (p0) -- (\x1,-0.8) -- (xi0);
\draw[dotted] let \p1 = (p1) in (p1) -- (\x1+5,-0.8) -- (xi025);
\draw[dotted] let \p1 = (p2) in (p2) -- (\x1+7,-0.8) -- (xi050);
\draw[dotted] let \p1 = (p3) in (p3) -- (\x1+14,-0.8) -- (xi075);
\draw[dotted] let \p1 = (p4) in (p4) -- (\x1+16,-0.8) -- (xi1);

{
  \draw (p0) circle [radius=0.1];
}
{
  \draw[color=black!60,thin,densely dotted] let \p1 = (p0) in ([shift=(110:1)]\x1,\y1) arc (110:-110:1);
  \draw (p1) circle [radius=0.1];
}
{
  \draw[color=black!70,thin,densely dotted] (p1) circle [radius=1];
  \draw (p2) circle [radius=0.1];
}
{
  \draw[color=black!80,thin,densely dotted] (p2) circle [radius=1];
  \draw (p3) circle [radius=0.1];
}
{
  \draw[color=black!90,thin,densely dotted] (p3) circle [radius=1];
  \draw (p4) circle [radius=0.1];
}

\draw[densely dashed] (p0)--(p1)--(p2)--(p3)--(p4);

\end{tikzpicture}

%% file: Figures/APALM_it1.pgf
\begin{tikzpicture}[scale=1.5,every node/.style={inner sep=10}]
\coordinate (O) at (0,0);
\coordinate (X) at (10,0);
\coordinate (Y) at (0,5);
\draw[-latex] (0,0) -- (5,0) node [right] {$\lambda$};
\draw[-latex] (0,0) -- (0,3.0) node [midway,above ,rotate=90]{$\Vert \mathbf{u}(t)\Vert$};

\draw[-latex] (0,-0.8) -- (5,-0.8)node [right] {$s$};
\draw[-latex] (0,-2) -- (5,-2)node [right] {$\xi$};
\draw (0,-2.1) -- (0,-1.9) node(xi0)[midway,inner sep=0]{};
\draw (4.5,-2.1) -- (4.5,-1.9) node(xi1)[midway,inner sep=0]{};
\node[above] at (xi0) {0};
\node[above] at (xi1) {1};

\draw[black!20] plot[smooth,tension=1] coordinates {(0,1) (2,2) (3,1)  (4.5,2)};

\draw[-latex](4.0,2.5)node[above,scale=1]{$G(\mathbf{u},\lambda)=0$} to[out=-90,in=120] (4.2,1.69);

\coordinate (p0) at (0,1);
\coordinate (p1) at (0.82,1.58);
\coordinate (p2) at (1.74,1.97);
\coordinate (p3) at (2.58,1.39);
\coordinate (p4) at (3.55,1.17);
\coordinate (p5) at (4.32,1.8);

\draw (p0) circle [radius=0.1];
\draw (p1) circle [radius=0.1];
\draw (p2) circle [radius=0.1];
\draw (p3) circle [radius=0.1];
\draw (p4) circle [radius=0.1];

\draw[fill=white] let \p1 = (p0) in (\x1,-0.8) circle [radius=0.1];
\draw[dotted] let \p1 = (p0) in (\x1,-0.8) -- (p0);
\draw[fill=white] let \p1 = (p1) in (\x1+5,-0.8) circle [radius=0.1];
\draw[dotted] let \p1 = (p1) in (\x1+5,-0.8) -- (p1);
\draw[fill=white] let \p1 = (p2) in (\x1+7,-0.8) circle [radius=0.1];
\draw[dotted] let \p1 = (p2) in (\x1+7,-0.8) -- (p2);
\draw[fill=white] let \p1 = (p3) in (\x1+14,-0.8) circle [radius=0.1];
\draw[dotted] let \p1 = (p3) in (\x1+14,-0.8) -- (p3);
\draw[fill=white] let \p1 = (p4) in (\x1+16,-0.8) circle [radius=0.1];
\draw[dotted] let \p1 = (p4) in (\x1+16,-0.8) -- (p4);

\draw (1.125,-2.1) -- (1.125,-1.9) node(xi025)[midway,inner sep=0]{};
\draw (2.25,-2.1) -- (2.25,-1.9) node(xi050)[midway,inner sep=0]{};
\draw (3.375,-2.1) -- (3.375,-1.9) node(xi075)[midway,inner sep=0]{};
\node[above] at (xi025) {0.25};
\node[above] at (xi050) {0.50};
\node[above] at (xi075) {0.75};

\draw[dotted] let \p1 = (p0) in (\x1,-0.8) -- (xi0);
\draw[dotted] let \p1 = (p1) in (\x1+5,-0.8) -- (xi025);
\draw[dotted] let \p1 = (p2) in (\x1+7,-0.8) -- (xi050);
\draw[dotted] let \p1 = (p3) in (\x1+14,-0.8) -- (xi075);
\draw[dotted] let \p1 = (p4) in (\x1+16,-0.8) -- (xi1);

\coordinate (p01) at (0.4,1.31);
\coordinate (p02) at (0.82,1.58);
\coordinate (p11) at (1.26,1.81);
\coordinate (p12) at (1.73,1.97);
\coordinate (p21) at (2.25,1.96);
\coordinate (p22) at (2.54,1.55);

\draw[color=gray,thin,densely dotted] let \p1 = (p0) in ([shift=(130:0.5)]\x1,\y1) arc (130:-130:0.5);
\draw[fill=gray,fill opacity=0.5]  (p01) circle [radius=0.05];
\draw[densely dashed] (p0)--(p01);
\draw[fill=gray,fill opacity=0.5] let \p1 = (p01) in (\x1+2.5,-0.8) circle [radius=0.05];
\draw[dotted] let \p1 = (p01) in (\x1+2.5,-0.8) -- (p01);
\draw (0.5625,-2.1) -- (0.5625,-1.9) node(xi0125)[midway,inner sep=0]{};
\node[above] at (xi0125) {0.125};
\draw[dotted] let \p1 = (p01) in (\x1+2.5,-0.8) -- (xi0125);


\draw[color=gray,thin,densely dotted] (p1) circle [radius=0.5];
\draw[fill=gray,fill opacity=0.5]  (p11) circle [radius=0.05];
\draw[densely dashed] (p1)--(p11);
\draw[fill=gray,fill opacity=0.5] let \p1 = (p11) in (\x1+6,-0.8) circle [radius=0.05];
\draw[dotted] let \p1 = (p11) in (\x1+6,-0.8) -- (p11);
\draw (1.6875,-2.1) -- (1.6875,-1.9) node(xi0375)[midway,inner sep=0]{};
\node[above] at (xi0375) {0.375};
\draw[dotted] let \p1 = (p11) in (\x1+6,-0.8) -- (xi0375);


\draw[color=gray,thin,densely dotted] (p2) circle [radius=0.5];
\draw[fill=gray,fill opacity=0.5]  (p21) circle [radius=0.05];
\draw[densely dashed] (p2)--(p21);






\end{tikzpicture}

%% file: Figures/APALM_it2.pgf
\begin{tikzpicture}[scale=1.5,every node/.style={inner sep=10}]
\coordinate (O) at (0,0);
\coordinate (X) at (10,0);
\coordinate (Y) at (0,5);
\draw[-latex] (0,0) -- (5,0) node [right] {$\lambda$};
\draw[-latex] (0,0) -- (0,3.0) node [midway,above ,rotate=90]{$\Vert \mathbf{u}(t)\Vert$};

\draw[-latex] (0,-0.8) -- (5,-0.8)node [right] {$s$};
\draw[-latex] (0,-2) -- (5,-2)node [right] {$\xi$};
\draw (0,-2.1) -- (0,-1.9) node(xi0)[midway,inner sep=0]{};
\draw (4.5,-2.1) -- (4.5,-1.9) node(xi1)[midway,inner sep=0]{};
\node[above] at (xi0) {0};
\node[above] at (xi1) {1};

\draw[black!20] plot[smooth,tension=1] coordinates {(0,1) (2,2) (3,1)  (4.5,2)};

\draw[-latex](4.0,2.5)node[above,scale=1]{$G(\mathbf{u},\lambda)=0$} to[out=-90,in=120] (4.2,1.69);

\coordinate (p0) at (0,1);
\coordinate (p1) at (0.82,1.58);
\coordinate (p2) at (1.74,1.97);
\coordinate (p3) at (2.58,1.39);
\coordinate (p4) at (3.55,1.17);
\coordinate (p5) at (4.32,1.8);

\draw (p0) circle [radius=0.1];
\draw (p1) circle [radius=0.1];
\draw (p2) circle [radius=0.1];
\draw (p3) circle [radius=0.1];
\draw (p4) circle [radius=0.1];

\draw[fill=white] let \p1 = (p0) in (\x1,-0.8) circle [radius=0.1];
\draw[dotted] let \p1 = (p0) in (\x1,-0.8) -- (p0);
\draw[fill=white] let \p1 = (p1) in (\x1+5,-0.8) circle [radius=0.1];
\draw[dotted] let \p1 = (p1) in (\x1+5,-0.8) -- (p1);
\draw[fill=white] let \p1 = (p2) in (\x1+7,-0.8) circle [radius=0.1];
\draw[dotted] let \p1 = (p2) in (\x1+7,-0.8) -- (p2);
\draw[fill=white] let \p1 = (p3) in (\x1+14,-0.8) circle [radius=0.1];
\draw[dotted] let \p1 = (p3) in (\x1+14,-0.8) -- (p3);
\draw[fill=white] let \p1 = (p4) in (\x1+16,-0.8) circle [radius=0.1];
\draw[dotted] let \p1 = (p4) in (\x1+16,-0.8) -- (p4);

\draw (1.125,-2.1) -- (1.125,-1.9) node(xi025)[midway,inner sep=0]{};
\draw (2.25,-2.1) -- (2.25,-1.9) node(xi050)[midway,inner sep=0]{};
\draw (3.375,-2.1) -- (3.375,-1.9) node(xi075)[midway,inner sep=0]{};
\node[above] at (xi025) {0.25};
\node[above] at (xi050) {0.50};
\node[above] at (xi075) {0.75};

\draw[dotted] let \p1 = (p0) in (\x1,-0.8) -- (xi0);
\draw[dotted] let \p1 = (p1) in (\x1+5,-0.8) -- (xi025);
\draw[dotted] let \p1 = (p2) in (\x1+7,-0.8) -- (xi050);
\draw[dotted] let \p1 = (p3) in (\x1+14,-0.8) -- (xi075);
\draw[dotted] let \p1 = (p4) in (\x1+16,-0.8) -- (xi1);

\coordinate (p01) at (0.4,1.31);
\coordinate (p02) at (0.82,1.58);
\coordinate (p11) at (1.26,1.81);
\coordinate (p12) at (1.73,1.97);
\coordinate (p21) at (2.25,1.96);
\coordinate (p22) at (2.54,1.55);

\draw[fill=gray,fill opacity=0.5]  (p01) circle [radius=0.05];
\draw[densely dashed] (p0)--(p01);
\draw[fill=gray,fill opacity=0.5] let \p1 = (p01) in (\x1+2.5,-0.8) circle [radius=0.05];
\draw[dotted] let \p1 = (p01) in (\x1+2.5,-0.8) -- (p01);
\draw (0.5625,-2.1) -- (0.5625,-1.9) node(xi0125)[midway,inner sep=0]{};
\node[above] at (xi0125) {0.125};
\draw[dotted] let \p1 = (p01) in (\x1+2.5,-0.8) -- (xi0125);

\draw[color=gray,thin,densely dotted] (p01) circle [radius=0.5];
\draw[fill=gray,fill opacity=0.5]  (p02) circle [radius=0.05];
\draw[densely dashed] (p01)--(p02);
\draw[fill=gray,fill opacity=0.5] let \p1 = (p02) in (\x1+5,-0.8) circle [radius=0.05];

\node[left,inner sep=2] at (p02) {\textcolor{green}{\textbf{$\checkmark$}}};

\draw[fill=gray,fill opacity=0.5]  (p11) circle [radius=0.05];
\draw[densely dashed] (p1)--(p11);
\draw[fill=gray,fill opacity=0.5] let \p1 = (p11) in (\x1+6,-0.8) circle [radius=0.05];
\draw[dotted] let \p1 = (p11) in (\x1+6,-0.8) -- (p11);
\draw (1.6875,-2.1) -- (1.6875,-1.9) node(xi0375)[midway,inner sep=0]{};
\node[above] at (xi0375) {0.375};
\draw[dotted] let \p1 = (p11) in (\x1+6,-0.8) -- (xi0375);

\draw[color=gray,thin,densely dotted] (p11) circle [radius=0.5];
\draw[fill=gray,fill opacity=0.5]  (p12) circle [radius=0.05];
\draw[densely dashed] (p11)--(p12);
\draw[fill=gray,fill opacity=0.5] let \p1 = (p12) in (\x1+7,-0.8) circle [radius=0.05];

\node[left,inner sep=2] at (p12) {\textcolor{green}{\textbf{$\checkmark$}}};

\draw[fill=gray,fill opacity=0.5]  (p21) circle [radius=0.05];
\draw[densely dashed] (p2)--(p21);

\draw[color=gray,thin,densely dotted] (p21) circle [radius=0.5];
\draw[fill=gray,fill opacity=0.5]  (p22) circle [radius=0.05];
\draw[densely dashed] (p21)--(p22);

\node[above right,inner sep=2] at (p22) {\textcolor{red}{\textbf{!}}};


\end{tikzpicture}

%% file: Figures/APALM_it3.pgf
\begin{tikzpicture}[scale=1.5,every node/.style={inner sep=10}]
\coordinate (O) at (0,0);
\coordinate (X) at (10,0);
\coordinate (Y) at (0,5);
\draw[-latex] (0,0) -- (5,0) node [right] {$\lambda$};
\draw[-latex] (0,0) -- (0,3.0) node [midway,above ,rotate=90]{$\Vert \mathbf{u}(t)\Vert$};

\draw[-latex] (0,-0.8) -- (5,-0.8)node [right] {$s$};
\draw[-latex] (0,-2) -- (5,-2)node [right] {$\xi$};
\draw (0,-2.1) -- (0,-1.9) node(xi0)[midway,inner sep=0]{};
\draw (4.5,-2.1) -- (4.5,-1.9) node(xi1)[midway,inner sep=0]{};
\node[above] at (xi0) {0};
\node[above] at (xi1) {1};

\draw[black!20] plot[smooth,tension=1] coordinates {(0,1) (2,2) (3,1)  (4.5,2)};

\draw[-latex](4.0,2.5)node[above,scale=1]{$G(\mathbf{u},\lambda)=0$} to[out=-90,in=120] (4.2,1.69);

\coordinate (p0) at (0,1);
\coordinate (p1) at (0.82,1.58);
\coordinate (p2) at (1.74,1.97);
\coordinate (p3) at (2.58,1.39);
\coordinate (p4) at (3.55,1.17);
\coordinate (p5) at (4.32,1.8);

\draw (p0) circle [radius=0.1];
\draw (p1) circle [radius=0.1];
\draw (p2) circle [radius=0.1];
\draw (p3) circle [radius=0.1];
\draw (p4) circle [radius=0.1];

\draw[fill=white] let \p1 = (p0) in (\x1,-0.8) circle [radius=0.1];
\draw[dotted] let \p1 = (p0) in (\x1,-0.8) -- (p0);
\draw[fill=white] let \p1 = (p1) in (\x1+5,-0.8) circle [radius=0.1];
\draw[dotted] let \p1 = (p1) in (\x1+5,-0.8) -- (p1);
\draw[fill=white] let \p1 = (p2) in (\x1+7,-0.8) circle [radius=0.1];
\draw[dotted] let \p1 = (p2) in (\x1+7,-0.8) -- (p2);
\draw[fill=red!20] let \p1 = (p3) in (\x1+20,-0.8) circle [radius=0.1];
\draw[dotted] let \p1 = (p3) in (\x1+20,-0.8) -- (p3);
\draw[fill=red!20] let \p1 = (p4) in (\x1+25,-0.8) circle [radius=0.1];
\draw[dotted] let \p1 = (p4) in (\x1+25,-0.8) -- (p4);

\draw (1.125,-2.1) -- (1.125,-1.9) node(xi025)[midway,inner sep=0]{};
\draw (2.25,-2.1) -- (2.25,-1.9) node(xi050)[midway,inner sep=0]{};
\draw (3.375,-2.1) -- (3.375,-1.9) node(xi075)[midway,inner sep=0]{};
\node[above] at (xi025) {0.25};
\node[above] at (xi050) {0.50};
\node[above] at (xi075) {0.75};

\draw[dotted] let \p1 = (p0) in (\x1,-0.8) -- (xi0);
\draw[dotted] let \p1 = (p1) in (\x1+5,-0.8) -- (xi025);
\draw[dotted] let \p1 = (p2) in (\x1+7,-0.8) -- (xi050);
\draw[dotted] let \p1 = (p3) in (\x1+20,-0.8) -- (xi075);
\draw[dotted] let \p1 = (p4) in (\x1+25,-0.8) -- (xi1);

\coordinate (p01) at (0.4,1.31);
\coordinate (p02) at (0.82,1.58);
\coordinate (p11) at (1.26,1.81);
\coordinate (p12) at (1.73,1.97);
\coordinate (p21) at (2.25,1.96);
\coordinate (p22) at (2.54,1.55);

\draw[fill=gray,fill opacity=0.5]  (p01) circle [radius=0.05];
\draw[color=gray,thin,densely dotted] (p01) circle [radius=0.5];
\draw[densely dashed] (p0)--(p01);
\draw[fill=gray,fill opacity=0.5] let \p1 = (p01) in (\x1+2.5,-0.8) circle [radius=0.05];
\draw[dotted] let \p1 = (p01) in (\x1+2.5,-0.8) -- (p01);
\draw (0.5625,-2.1) -- (0.5625,-1.9) node(xi0125)[midway,inner sep=0]{};
\draw[dotted] let \p1 = (p01) in (\x1+2.5,-0.8) -- (xi0125);

\draw[fill=gray,fill opacity=0.5]  (p02) circle [radius=0.05];
\draw[densely dashed] (p01)--(p02);
\draw[fill=gray,fill opacity=0.5] let \p1 = (p02) in (\x1+5,-0.8) circle [radius=0.05];

\draw[fill=gray,fill opacity=0.5]  (p11) circle [radius=0.05];
\draw[color=gray,thin,densely dotted] (p11) circle [radius=0.5];
\draw[densely dashed] (p1)--(p11);
\draw[fill=gray,fill opacity=0.5] let \p1 = (p11) in (\x1+6,-0.8) circle [radius=0.05];
\draw[dotted] let \p1 = (p11) in (\x1+6,-0.8) -- (p11);
\draw (1.6875,-2.1) -- (1.6875,-1.9) node(xi0375)[midway,inner sep=0]{};
\draw[dotted] let \p1 = (p11) in (\x1+6,-0.8) -- (xi0375);

\draw[fill=gray,fill opacity=0.5]  (p12) circle [radius=0.05];
\draw[densely dashed] (p11)--(p12);
\draw[fill=gray,fill opacity=0.5] let \p1 = (p12) in (\x1+7,-0.8) circle [radius=0.05];


\draw[fill=gray,fill opacity=0.5]  (p21) circle [radius=0.05];
\draw[color=gray,thin,densely dotted] (p21) circle [radius=0.5];
\draw[densely dashed] (p2)--(p21);

\draw[fill=red,fill opacity=0.5] let \p1 = (p21) in (\x1+9,-0.8) circle [radius=0.05];
\draw[dotted] let \p1 = (p21) in (\x1+9,-0.8) -- (p21);

\draw[fill=gray,fill opacity=0.5]  (p22) circle [radius=0.05];
\draw[densely dashed] (p21)--(p22);

\draw[fill=red,fill opacity=0.5] let \p1 = (p3) in (\x1+14,-0.8) circle [radius=0.05];
\draw[dotted] let \p1 = (p3) in (\x1+14,-0.8) -- (p22);

\draw (2.7125,-2.1) -- (2.7125,-1.9) node(xi0375)[midway,inner sep=0]{};
\draw[dotted] let \p1 = (p21) in (\x1+9,-0.8) -- (xi0375);
\draw (3.200,-2.1) -- (3.200,-1.9) node(xi0375)[midway,inner sep=0]{};
\draw[dotted] let \p1 = (p3) in (\x1+14,-0.8) -- (xi0375);

\node[above right,inner sep=2] at (p22) {\textcolor{red}{\textbf{!}}};

%

%


\end{tikzpicture}

%% file: Figures/Lengths_regular.pgf
\begin{tikzpicture}[scale=1.5]
\coordinate (O) at (0,0);
\coordinate (X) at (10,0);
\coordinate (Y) at (0,5);
\draw[-latex] (0,0) -- (5,0) node [right] {$\lambda$};
\draw[-latex] (0,0) -- (0,3.0) node [midway,above ,rotate=90]{$\Vert \mathbf{u}(\lambda)
\Vert$};

\coordinate (A) at (0.5,1);
\coordinate (B) at (2,1.55);
\coordinate (C) at (3,1.5);
\coordinate (D) at (1.67,1.47);
\coordinate (E) at (2.95,1.52);
\draw[black!20] plot[smooth,tension=1] coordinates {(A) (B) (C)};


\draw[-latex] (0,-0.8) -- (5,-0.8)node [right] {$s$};
\draw[-latex] (0,-2.2) -- (5,-2.2)node [right] {$\xi$};
\draw (0,-2.3) -- (0,-2.1) node(xi0)[midway]{};
\draw (4.0,-2.3) -- (4.0,-2.1) node(xi1)[midway]{};
\draw (1.975,-2.3) -- (1.975,-2.1) node(xi045)[midway]{};
\draw (3.95,-2.3) -- (3.95,-2.1) node(xi090)[midway]{};
\node[above,inner sep=8] at (xi0) {$\xi^\ell_i = \xi^{\ell+1}_0$};
\node[above,inner sep=8] at (xi1) {$\xi^\ell_{i+1}\approx \xi^{\ell+1}_N$};
\node[above,inner sep=8] at (xi045) {$\xi^{\ell+1}_1$};

\draw[black] (A) circle [radius=0.1];
\draw[black] (C) circle [radius=0.1];

\draw[fill=white] let \p1 = (xi0.center) in (\x1,-0.8) node (p0){} circle [radius=0.1];
\draw[fill=white] let \p1 = (xi1.center) in (\x1,-0.8) node (p1){} circle [radius=0.1];
\node[above,inner sep=10] at (p0) {$s^\ell_i=s^{\ell+1}_0$};
\node[above,inner sep=10] at (p1) {$s^\ell_{i+1}\approx s^{\ell+1}_N$};

\draw[latex-latex] let \p0 = (xi0.center),\p1=(xi1.center) in (\x0,-1.35) -- node[midway,below]{$\Delta L$} (\x1,-1.35);
\draw[latex-latex] let \p0 = (xi0.center),\p1=(xi045.center) in (\x0,-0.95) -- node[midway,below]{$\frac{\Delta L}{2}$} (\x1,-0.95);
\draw[latex-latex] let \p0 = (xi045.center),\p1=(xi090.center) in (\x0,-0.95) -- node[midway,below]{$\frac{\Delta L}{2}$} (\x1,-0.95);

\filldraw[draw=black,fill=gray] let \p1 = (xi045.center) in (\x1,-0.8) node (p045){}  circle [radius=0.05];
\filldraw[draw=black,fill=white] let \p1 = (xi090.center) in (\x1,-0.8) node (p090){}  circle [radius=0.05];

\node[above,inner sep=10] at (p045) {$s^{\ell+1}_1$};

\filldraw[draw=black,fill=gray] (D) circle [radius=0.05];
\filldraw[draw=black,fill=white] (E) circle [radius=0.05];

\draw[dotted] (A)--(p0.center);
\draw[dotted] (C)--(p1.center);
\draw[dotted] (D)--(p045.center);
\draw[dotted] (p0.center)--(xi0.center);
\draw[dotted] (p1.center)--(xi1.center);
\draw[dotted] (p045.center)--(xi045.center);

\draw[dashed] (A)-- node[midway,inner sep=0](ACm){} (C);
\draw[dashed] (A)-- node[midway,above,inner sep=0]{$\frac{\Delta L}{2}$}(D);
\draw[dashed] (D)-- node[midway,above,inner sep=0]{$\frac{\Delta L}{2}$}(E);
\node [](ACl) at (2.3,0.5){$\Delta L$};
\draw[-latex](ACl) to[out=100,in=-60] (ACm);
\end{tikzpicture}

%% file: Figures/Lengths_irregular.pgf
\begin{tikzpicture}[scale=1.5]
\coordinate (O) at (0,0);
\coordinate (X) at (10,0);
\coordinate (Y) at (0,5);
\draw[-latex] (0,0) -- (5,0) node [right] {$\lambda$};
\draw[-latex] (0,0) -- (0,3.0) node [midway,above ,rotate=90]{$\Vert \mathbf{u}(\lambda)
\Vert$};

\coordinate (A) at (0.5,1);
\coordinate (B) at (2,2);
\coordinate (C) at (3,1.5);
\coordinate (D) at (1.5,1.8);
\coordinate (E) at (2.78,1.78);
\draw[black!20] plot[smooth,tension=1] coordinates {(A) (B) (C)};
%

\draw[-latex] (0,-0.8) -- (5,-0.8)node [right] {$s$};
\draw[-latex] (0,-2.2) -- (5,-2.2)node [right] {$\xi$};
\draw (0,-2.3) -- (0,-2.1) node(xi0)[midway]{};
\draw (4.0,-2.3) -- (4.0,-2.1) node(xi1)[midway]{};
\draw (1.8,-2.3) -- (1.8,-2.1) node(xi045)[midway]{};
\draw (3.6,-2.3) -- (3.6,-2.1) node(xi090)[midway]{};
\node[above,inner sep=8] at (xi0) {$\xi^\ell_i = \xi^{\ell+1}_0$};
\node[above,inner sep=8] at (xi1) {$\xi^\ell_{i+1}$};
\node[above,inner sep=8] at (xi045) {$\xi^{\ell+1}_1$};
\node[above,inner sep=8] at (xi090) {$\xi^{\ell+1}_N$};

\draw[black] (A) circle [radius=0.1];
\draw[black] (C) circle [radius=0.1];

\draw[fill=white] let \p1 = (xi0) in (\x1,-0.8) node (p0){} circle [radius=0.1];
\draw[fill=white] let \p1 = (xi1) in (\x1+10,-0.8) node (p1){} circle [radius=0.1];
\node[above,inner sep=10] at (p0) {$s^\ell_i$};
\node[above,inner sep=10] at (p1) {$s^\ell_{i+1}$};

\filldraw[draw=black,fill=gray] let \p1 = (xi045) in (\x1,-0.8) node (p045){}  circle [radius=0.05];
\filldraw[draw=black,fill=gray] let \p1 = (xi090) in (\x1,-0.8) node (p090){}  circle [radius=0.05];

\draw[latex-latex] let \p0 = (p0),\p1=(p1) in (\x0,-1.35) -- node[midway,below]{$\Delta L'$} (\x1,-1.35);
\draw[latex-latex] let \p0 = (p0),\p1=(p045) in (\x0,-0.95) -- node[midway,below]{$\frac{\Delta L}{2}$} (\x1,-0.95);
\draw[latex-latex] let \p0 = (p045),\p1=(p090) in (\x0,-0.95) -- node[midway,below]{$\frac{\Delta L}{2}$} (\x1,-0.95);
\draw[latex-latex] let \p0 = (p090),\p1=(p1) in (\x0,-0.95) -- node[midway,below]{$\delta L$} (\x1,-0.95);

\node[above,inner sep=10] at (p045) {$s^{\ell+1}_{1}$};
\node[above,inner sep=10] at (p090) {$s^{\ell+1}_{N}$};

\filldraw[draw=black,fill=gray] (D) circle [radius=0.05];
\filldraw[draw=black,fill=gray] (E) circle [radius=0.05];

\draw[dotted] (A)--(p0);
\draw[dotted] (C)--(p1);
\draw[dotted] (D)--(p045);
\draw[dotted] (E)--(p090);
\draw[dotted] (p0)--(xi0);
\draw[dotted] (p1)--(xi1);
\draw[dotted] (p045)--(xi045);
\draw[dotted] (p090)--(xi090);

\draw[dashed] (A)-- node[midway,inner sep=0](ACm){} (C);
\draw[dashed] (A)-- node[midway,above,inner sep=0]{$\frac{\Delta L}{2}$}(D);
\draw[dashed] (D)-- node[midway,above,inner sep=0]{$\frac{\Delta L}{2}$}(E);
\draw[dashed] (A)-- node[midway,inner sep=0](AEm){}(E);
\draw[dashed] (E)-- node[midway,inner sep=0](ECm){}(C);
\node [](ACl) at (2.3,0.5){$\Delta L$};
\node [](AEl) at (1.4,0.6){$\overline{\Delta L}$};
\node [](ECl) at (3.9,2.0){$\delta L$};
\draw[-latex](ACl) to[out=100,in=-60] (ACm);
\draw[-latex](AEl) to[out=90,in=-100] (AEm);
\draw[-latex](ECl) to[out=180,in=30] (ECm);
\end{tikzpicture}

%% file: Figures/domain_ori.pgf
\begin{tikzpicture}
\def\dy{1.5};
\def\ddy{0.5};
\def\R{0.1};
\def\dxi{0.8};
\def\NN{5};
\def\dx{4};
\pgfmathsetmacro\ddx{\dx / (\NN)}

\draw[white] ($(0,5)+(-2*\dxi,0)$) -- ($(0,5)+(4*\dxi,0)$);

\coordinate (s0begin) at ($(0,5)+(\dxi,0)$);
\coordinate (s0end) at ($(s0begin)+(\dx,0)$);
\draw (s0begin) --(s0end);
\draw[dashed] ($(s0begin)+(-\dxi,0)$) node[left]{$s^{\ell}$} -- (s0begin);
\draw[dashed] (s0end) -- ($(s0end)+(3*\dxi,0)$);
\draw[fill=white] ($(s0end)+(2*\dxi,0)$) circle [radius=\R];

\coordinate (xi0begin) at ($(s0begin)+(0,-\ddy)$);
\coordinate (xi0end) at ($(s0end)+(0,-\ddy)$);
\draw (xi0begin) --(xi0end);
\draw[dashed] ($(xi0begin)+(-\dxi,0)$) node[left]{$\xi^{\ell}$} -- (xi0begin);
\draw[dashed] (xi0end) -- ($(xi0end)+(3*\dxi,0)$);
\draw[fill=white] ($(xi0end)+(2*\dxi,0)$) circle [radius=\R];

\foreach \k in {0,...,1}
{
\draw[fill=white] ($(s0begin)+(\k*\dx,0)$) circle [radius=\R];
\draw[fill=white] ($(xi0begin)+(\k*\dx,0)$) circle [radius=\R];
}
\end{tikzpicture}

%% file: Figures/domain_interior.pgf
\begin{tikzpicture}
\def\dy{1.5};
\def\ddy{0.5};
\def\R{0.1};
\def\dxi{0.8};
\def\NN{5};
\def\dx{4};
\pgfmathsetmacro\ddx{\dx / (\NN)}

\draw[white] ($(0,5)+(-2*\dxi,0)$) -- ($(0,5)+(4*\dxi,0)$);

\coordinate (s1begin) at ($(0,5)+(\dxi,0)$);
\coordinate (s1end) at ($(s1begin)+(\dx,0)$);
\draw (s1begin) --(s1end);
\draw[dashed] ($(s1begin)+(-\dxi,0)$) node[left]{$s^{\ell+1}$} -- (s1begin);
\draw[dashed] (s1end) -- ($(s1end)+(3*\dxi,0)$);
\draw[fill=white] ($(s1end)+(2*\dxi,0)$) circle [radius=\R];

\coordinate (xi1begin) at ($(s1begin)+(0,-\ddy)$);
\coordinate (xi1end) at ($(s1end)+(0,-\ddy)$);
\draw (xi1begin) --(xi1end);
\draw[dashed] ($(xi1begin)+(-\dxi,0)$) node[left]{$\xi^{\ell+1}$} -- (xi1begin);
\draw[dashed] (xi1end) -- ($(xi1end)+(3*\dxi,0)$);
\draw[fill=white] ($(xi1end)+(2*\dxi,0)$) circle [radius=\R];

\foreach \k in {0,...,1}
{
\draw[fill=white] ($(s1begin)+(\k*\dx,0)$) circle [radius=\R];
\draw[fill=white] ($(xi1begin)+(\k*\dx,0)$) circle [radius=\R];
}

\pgfmathsetmacro\NNm{\NN-1}
\pgfmathsetmacro\dxxi{0*\dx}
\foreach \k in {1,...,\NNm}
{
\draw[fill=gray] ($(s1begin)+(\k*\ddx+\dxxi,0)$) circle [radius=\R];
}
\draw[fill=lightgray,dotted,thick] ($(s1begin)+(\NN*\ddx-0.1*\ddx,0)$) circle [radius=\R];

\foreach \k in {1,...,\NNm}
{
\pgfmathsetmacro\dxxi{0*\dx}
\draw[fill=gray] ($(xi1begin)+(\k*\ddx+\dxxi,0)$) circle [radius=\R];
}

\end{tikzpicture}

%% file: Figures/domain_stretch.pgf
\begin{tikzpicture}
\def\dy{1.5};
\def\ddy{0.5};
\def\R{0.1};
\def\dxi{0.8};
\def\NN{5};
\def\dx{4};
\pgfmathsetmacro\ddx{\dx / (\NN)}

\draw[white] ($(0,5)+(-2*\dxi,0)$) -- ($(0,5)+(4*\dxi,0)$);

\coordinate (s2begin) at ($(0,5)+(\dxi,0)$);
\coordinate (s2end) at ($(s2begin)+(\dx,0)$);
\draw (s2begin) --(s2end);
\draw[dashed] ($(s2begin)+(-\dxi,0)$) node[left]{$s^{\ell+1}$} -- (s2begin);
\draw[dashed] (s2end) -- ($(s2end)+(3*\dxi,0)$);
\draw[fill=white] ($(s2end)+(2.5*\dxi,0)$) circle [radius=\R];

\coordinate (xi2begin) at ($(s2begin)+(0,-\ddy)$);
\coordinate (xi2end) at ($(s2end)+(0,-\ddy)$);
\draw (xi2begin) --(xi2end);
\draw[dashed] ($(xi2begin)+(-\dxi,0)$) node[left]{$\xi^{\ell+1}$} -- (xi2begin);
\draw[dashed] (xi2end) -- ($(xi2end)+(3*\dxi,0)$);
\draw[fill=white] ($(xi2end)+(2*\dxi,0)$) circle [radius=\R];

\draw[fill=white] ($(s2begin)+(0*\dx,0)$) circle [radius=\R];
\draw[fill=white] ($(xi2begin)+(0*\dx,0)$) circle [radius=\R];

\draw[fill=white] ($(s2begin)+(1*\dx+0.5*\dxi,0)$) circle [radius=\R];
\draw[fill=white] ($(xi2begin)+(1*\dx,0)$) circle [radius=\R];

\pgfmathsetmacro\NNm{\NN-1}
\pgfmathsetmacro\dxxi{0*\dx}
\foreach \k in {1,...,\NNm}
{
\draw[fill=gray] ($(s2begin)+(\k*\ddx+\dxxi,0)$) circle [radius=\R];
\draw[fill=gray] ($(xi2begin)+(\k*\ddx+\dxxi-0.15*\ddx,0)$) circle [radius=\R];
}
\draw[fill=lightgray] ($(s2begin)+(\NN*\ddx-0.5*\ddx,0)$) circle [radius=\R];
\draw[fill=lightgray] ($(xi2begin)+(\NN*\ddx-0.75*\ddx,0)$) circle [radius=\R];

\end{tikzpicture}

%% file: Sections/Implementation.tex
In this section, data structures and algorithms for the implementation of the APALM are presented. In \cref{subsec:APALM_datastructure}, a data structure is provided for the implementation of the APALM. Thereafter, \cref{subsec:APALM_algorithms} provides algorithms for the implementation of the APALM, and \cref{subsec:APALM_exploration} elaborates on the extension of these methods to multiple branches, hence enabling arc-length exploration.

\subsection{Data Structure}\label{subsec:APALM_datastructure}
Since the APALM is based on a sub-interval approach where the start and end points of each sub-interval are known, a data structure referencing the sub-intervals is essential. Since the curve-length coordinate is subject to change after reparametrisation of the curve and since the curve parameter is fixed, the logical choice is to connect the data to parametric coordinates. That is, a series of discrete maps is constructed such that solutions, levels, and curve-length coordinates can be obtained via a parametric point $\xi^\ell_k$.\\

\Cref{fig:APALM_datastructure} shows the data structure behind the APALM. Firstly, the data structure contains the map $\mathcal{S}(\xi):[0,1]\to[0,S]$, which is the map that maps the parametric coordinate to the curve-length domain. Secondly, the maps $\mathcal{U}(\xi):[0,1]\to\mathbb{R}^{n+1}$ and $\mathcal{U}'(\xi):[0,1]\to\mathbb{R}^{n+1}$ map the solution and the previous solution from the parametric domain to the solution and previous solution domain, respectively. The mapper $\mathcal{U}'(\xi):[0,1]\to\mathbb{R}^{n+1}$ is constructed in order to construct the predictor of the ALM. Lastly, the map $\mathcal{L}(\xi):[0,1]\to\mathbb{N}$ is a map that can be used to obtain the level of a parametric point, which is optional but can be useful to limit the method to a certain depth.\\

\begin{figure}
  \centering
    \input{Figures/DataStructure.pgf}
  \caption{The data structure behind the APALM. The axes represent data sets, which are monotonically increasing when the axis has an arrow. Solid arrows represent mappers from one axis to another, and dashed arrows represent data references. The mappers $\Xi(s_i)$ and $\mathcal{S}(\xi_i)$ map between the curve parametrisation and the curve length axes. The former takes a curve length $s_i$ and returns the curve parameter $\xi_i$, and the latter maps the inverse. The mappers $\mathcal{U}(\xi_k)$ and $\mathcal{U}'(\xi_k)$ return the solution $\VEC{w}_j$ and the previous solution $(\VEC{w}')_j$, respectively, given a parametric coordinate $\xi_j$, and the mapper $\mathcal{L}(\xi_j)$ returns the level on which the coordinate $\xi_j$ was computed. The guess is a data reference to the previous solution. The thick solid intervals represent running jobs assigned with an \texttt{ID}, and the thick dashed intervals represent queued intervals. Each interval is represented by a start-point and an end-point tuple $(\xi_l,\xi_{l+1})$. The \textcolor{red}{red} lines, squares, and arros represent the \texttt{submit} operation when solutions are added to the data structure.}
  \label{fig:APALM_datastructure}
\end{figure}

\subsection{Algorithms}\label{subsec:APALM_algorithms}
Given the underlying data structure of the APALM (see \cref{subsec:APALM_datastructure}), algorithms are defined for its implementation. Firstly, it is assumed that the APALM is based on an ALM with possibly a black-box implementation, striving for the non-intrusiveness of the method. The required routines for the underlying ALM are:
\begin{itemize}
  \item $\VEC{w}^\ell_{i+1}\gets$ \texttt{step}($\VEC{w}^\ell_i,\Delta\VEC{w}^\ell_i\Delta L$): Performs a step with length $\Delta L$ starting at point $\VEC{w}^\ell_i$ and returns the new solution $\VEC{w}^\ell_{i+1}$. Given the current solution and the previous solution, a predictor for the initial iteration that employs $\Delta\VEC{w}^\ell_i=\VEC{w}^\ell_i-\VEC{w}^\ell_{i-1}$ could be available, as well as one for a cold start, i.e., $\Delta\VEC{w}^\ell_i=\VEC{0}$.
  \item $\Delta s\gets$\texttt{distance}$(\VEC{w}^\ell_i,\VEC{w}^\ell_j)$: gets the distance between two points $\VEC{w}^\ell_i$ and $\VEC{w}^\ell_j$, using \cref{eq:APALM_distance} with $\Delta \VEC{w}^{\ell}_i=\VEC{w}^\ell_i-\VEC{w}^\ell_j$.
\end{itemize}
In the following, three implementations are presented. The first implementation is a \emph{serial implementation} without communication but with queueing, referred to as the Adaptive Serial Arc-Length Method (ASALM). The serial implementation provides the building blocks for the \emph{hybrid implementation} and the \emph{parallel implementation}. The \emph{hybrid implementation}, referred to as the Adaptive Serial-Parallel Arc-Length Method (ASPALM), is a hybrid version of the APALM where parallel corrections are performed after a serial solve has finished. The \emph{parallel implementation}, on the other hand, starts parallel corrections as soon as the first interval has been initialised; hence, there is no separation between a serial phase and a parallel phase. The parallel implementation is the final APALM.

\subsubsection*{Serial implementation}
The global workflow for a serial APALM, i.e., the ASALM, is illustrated in \cref{alg:ASALM_full}. As seen in this algorithm, the initialization is performed using a \texttt{serialSolve} routine, which defines the initial solution sequence $\{\VEC{w}^0_i\}_{i=0}^{I}$ with $I$ steps on level $\ell=0$. In addition, this routine also provides a sequence of curve-length coordinates, $\{s^0_i\}_{i=0}^{I}$. Based on these sequences, the mappers from \cref{fig:APALM_datastructure} and a queue $Q$ are initialised in the \texttt{initializeMap} routine. Given the queue $Q$, the \texttt{correctQueueSerial} routine provides a sequence of solutions $\{\VEC{w}_i\}$ and of curve parameters $\{s_i\}$ spanning multiple levels, hence the superscript $\ell$ is omitted. \\

\begin{algorithm}[h]
  \caption{Global ASALM routine (\texttt{ASALM}). The ASALM first consists of a serial solve of the whole curve length domain, followed by an evaluation of subintervals.}
  \label{alg:ASALM_full}
    \begin{algorithmic}[1]
      \Require $\Delta L$, $I$
      \State $\{\VEC{w}^0_i\}_{i=0}^{I},\{s^0_i\}_{i=0}^{I}\gets\texttt{serialSolve}(\Delta L,I)$
      \State $Q\gets\texttt{initializeMap}(\{\VEC{w}^0_i\}_{i=0}^{I},\{s^0_i\}_{i=0}^{I})$
      \State $\{\VEC{w}_i\},\{s_i\}\gets\texttt{correctQueueSerial}(Q)$
      \Ensure $\{\VEC{w}_i\}$, $\{s_i\}$, $\{\xi_i\}$
    \end{algorithmic}
\end{algorithm}

Using basic ALM routines, \cref{alg:serialSolve} defines an algorithm to obtain in serial a coarse approximation to initialise the ASALM. Optionally, a stability computation can be performed after the arc-length step, which could lead to a specialised solution towards a bifurcation point. This allows for automatic exploration of bifurcation diagrams \cite{Thies2021,Wouters2019} and is discussed more in detail in \cref{subsec:APALM_exploration}.\\

\begin{algorithm}[h]
  \caption{Serial solve (\texttt{serialSolve}). This routine provides the initial step for the APALM/ASALM, i.e., the solution data $\{\VEC{w}^\ell_i\}_{i=0,...,I}$ and the corresponding curve-length parameters $\{s^\ell_i\}_{i=0,...,I}$ on level $\ell=0$.}
  \label{alg:serialSolve}
  \begin{algorithmic}[1]
    \Require $\Delta L$, $I$
    \State Initialise $\VEC{w}^\ell_0$, $s^\ell_0=0$
    \State $\VEC{w}^\ell_1 \gets $\texttt{step}$(\VEC{w}^\ell_0,\VEC{0},\Delta L)$ \Comment{Compute first solution}
    \State $s_1\gets\Delta L$
    \For{$k=1,...,I-1$}
      \State $\VEC{w}^\ell_{k+1},\Delta s^\ell_{k+1} \gets $\texttt{initiate}$(\VEC{w}^\ell_{k},\Delta\VEC{w}_{k-1}^\ell,\Delta L)$ \Comment{Compute new solution}
      \State $s^\ell_{k+1}=s^\ell_{k}+\Delta s^\ell_{k+1}$ \Comment{Compute curve coordinate}
    \EndFor
    \Ensure $\{\VEC{w}^\ell_i\}_{i=0}^{I}$, $\{s^\ell_i\}_{i=0}^{I}$
  \end{algorithmic}
\end{algorithm}

In \cref{alg:serialSolve}, the \texttt{initiate} routine (see \cref{alg:initiate}) computes an arc-length interval and returns a new point $\VEC{w}^\ell_{k}$ and the traversed distance $\Delta s_k^\ell$, provided the previous point $\VEC{w}^\ell_{k-1}$, the previous solution interval $\Delta\VEC{w}_k^\ell$ and the intended arc-length step size $\Delta L$ and using the $\texttt{step}$ and $\texttt{distance}$ functions. Typically, $\Delta s_k^\ell$ is equal to $\Delta L$ unless the arc-length step does not converge and needs to be bisected.

\begin{algorithm}[h]
  \caption{Initiation of an interval (\texttt{initiate}). Given a previous solution $\VEC{w}^\ell_{k}$, the previous step size $\VEC{w}^\ell_{k-1}$, and the desired arc-length $\Delta L$, this routine returns a new solution $\VEC{w}^\ell_{k+1}$ and a distance with respect to the previous solution $\VEC{w}^\ell_{k}$, denoted by $\Delta s^\ell_{k+1}$.}
  \label{alg:initiate}
  \begin{algorithmic}[1]
    \Require $\VEC{w}^\ell_{k}$ , $\VEC{w}^\ell_{k-1}$, $\Delta L$
    \State $\Delta\VEC{w}^\ell_k = \VEC{w}^\ell_k-\VEC{w}^\ell_{k-1}$ \Comment{Compute previous step}
    \State $\VEC{w}^\ell_{k+1} \gets $\texttt{step}$(\VEC{w}^\ell_{k-1},\Delta\VEC{w}^\ell_k,\Delta L)$ \Comment{Compute new solution}
    \State $\Delta s^\ell_{k+1}\gets $\texttt{distance}$(\VEC{w}^\ell_{k+1},\VEC{w}^\ell_{k})$ \Comment{Compute curve coordinate}
    \Ensure $\VEC{w}^\ell_{k+1}$, $\Delta s^\ell_{k+1}$
  \end{algorithmic}
\end{algorithm}

As soon as a set of solution data and curve parameters, $\{\VEC{w}^\ell_i\}_{i=0,...,I}$ and $\{s^\ell_i\}_{i=0,...,I}$, respectively, are known, the initialization of the parallel computations can take place. In this initialization, the maps from \cref{fig:APALM_datastructure} are constructed and a queue $Q$ of jobs is created, c.f. \cref{alg:initialise}.

\begin{algorithm}[h]
  \caption{Parallel initialization (\texttt{initializeMap}). Within this algorithm, the maps $\mathcal{U}$, $\mathcal{U}'$, $\mathcal{S}$, and $\Xi$ are constructed from a series of solutions and corresponding curve length coordinates, $\{\VEC{w}^\ell_i\}_{i=0}^{I}$ and $\{s^\ell_i\}_{i=0}^{I}$, respectively, both on level $\ell=0$. Furthermore, the queue $Q$ is initialised by adding all subintervals to the queue.}
  \label{alg:initialise}
  \begin{algorithmic}[1]
    \Require $\{\VEC{w}^\ell_i\}_{i=0}^{I}$, $\{s^\ell_i\}_{i=0}^{I}$
    \State Add $\VEC{w}^\ell_0$ to $\mathcal{U}$, $s^\ell_0$ to $\mathcal{S}$, and $\xi^\ell_0$ to $\Xi$ \Comment{Add the start of the level 0 solutions to the map $\mathcal{U}$}
    \For{$k=1,...,I$}
    \State Add $\VEC{w}^\ell_k$ to $\mathcal{U}$ and $\VEC{w}^\ell_{k-1}$ to $\mathcal{U}'$.
    \State Add $s^\ell_k$ to $\mathcal{S}$ and $\xi^0_k$ to $\Xi$.
    \State Add $Q_{k-1}=[\xi^\ell_{k-1},\xi^\ell_{k})$ to Q. \Comment{Construct elements of the queue $Q$}
    \EndFor
    \Ensure $Q=\{Q_k=[\xi^\ell_k,\xi^\ell_{k+1}]\}_{k=0,...,I}$
  \end{algorithmic}
\end{algorithm}

After initialization, the computation of the sub-intervals can take place. This requires the routine \texttt{correctQueueSerial} as defined in \cref{alg:correctQueueSerial} for fully serial computations. That is, no communication between manager and worker takes place since everything will be done on the same node.\\

\begin{algorithm}[h]
  \caption{The routine that solves the queue (\texttt{correctQueueSerial}). Given a queue $Q$, this algorithm takes an entry from the queue using \texttt{pop} and solves the defined interval using \texttt{correct}. The new solution is added to the solution maps, and if required, new jobs are added to the queue $Q$ using \texttt{submit}. The final solutions are collected from the maps using the \texttt{collectSolutions} routine.}
  \label{alg:correctQueueSerial}
  \begin{algorithmic}[1]
  \Require $Q$
  \While{$Q\neq\emptyset$}
  \State $Q$, $\texttt{ID}$, $\Delta L$, $\VEC{w}^\ell_i$, $\VEC{w}^\ell_{i-1},\VEC{w}^\ell_{i+1}$$\gets\texttt{pop}(Q)$
  \State $\{d^{\ell+1}_k\}_{k=0,...,N-1}$, $\{\VEC{w}^{\ell+1}_k\}_{k=0,...,N}$, $\overline{\Delta L}$, $\delta L$$\gets\texttt{correct}(N,\Delta L_0, \VEC{w}^\ell_{i}, \VEC{w}^\ell_{i-1}, \VEC{w}^\ell_{i+1})$
  \State $Q$$\gets\texttt{submit}(\texttt{ID},\{d^{\ell+1}_k\}_{k=0,...,N-1},\{\VEC{w}^{\ell+1}_k\}_{k=0,...,N},\overline{\Delta L},\delta L,Q)$ \Comment{Submit the job; adds new jobs to $Q$ if needed}
  \EndWhile
  \State $\{\VEC{w}_i\}$, $\{s_i\}$$\gets\texttt{collectSolutions}$
  \Ensure $\{\VEC{w}_i\}$, $\{s_i\}$
  \end{algorithmic}
\end{algorithm}

The computation of the sub-interval takes place in the \texttt{correct} routine of \cref{alg:correct} and can be used both in a serial and a parallel implementation. Given a number of sub-intervals $N$, the original distance between the end-points of the interval $\Delta L_0$ and the start point, previous solutions, and reference solutions $\VEC{w}^\ell_i$, $\VEC{w}^\ell_{i-1}$ and $\VEC{w}^\ell_{i+1}$, respectively, this routine computes a series of solutions of the sub-interval $\{\VEC{w}^{\ell+1}_k\}_{k=0,...,N}$, their distances $\{d_k\}^{\ell+1}_{k=0,...,N-1}$ and the distances $\overline{\Delta L}$ and $\delta L$ for error computation. Note that the distance $\Delta L'$ can be computed by taking the sum of the distances.

\begin{algorithm}[h]
  \caption{The routine that solves an interval (\texttt{correct}). This routine takes a number of subintervals $N$, the desired step length for the total interval, the start point $\VEC{w}^\ell_i$, the previous point $\VEC{w}^\ell_{i-1}$, and the next point $\VEC{w}^\ell_{i+1}$. It returns the solutions on the subinterval and the distances between them, respectively $\{d^{\ell+1}_j\}_{j=0,...,N-1}$ and $\{\VEC{w}^{\ell+1}_j\}_{j=0,...,N}$, as well as the distances $\overline{\Delta L}$ and $\delta L$. When the \texttt{step} does not converge, it is assumed that step size modification takes place and that the data points are adjusted accordingly.}
  \label{alg:correct}
  \begin{algorithmic}[1]
    \Require $N$, $\Delta L_0$, $\VEC{w}^\ell_{i}$, $\VEC{w}^\ell_{i-1}$, $\VEC{w}^\ell_{i+1}$
    \State Initialise output vectors $\{d^{\ell+1}_k\}_{k=0,...,N-1}$, $\{\VEC{w}^{\ell+1}_k\}_{k=0,...,N}$
    \State $\Delta L=\Delta L_0/N$ \Comment{Defines the size of the sub-intervals}
    \State $\VEC{w}^{\ell+1}_0=\VEC{w}^{\ell}_i$
    \State $\Delta \VEC{w}^{\ell+1}=\VEC{w}^\ell_i-\VEC{w}^\ell_{i-1}$ \Comment{Determine previous arc-length step, to be used to predict the step on $\ell+1$}
    \For{$k=0,...,N-1$}
    \State $\VEC{w}_{k+1}^{\ell+1}\gets\texttt{step}(\VEC{w}_k^{\ell+1},\Delta\VEC{w},\Delta L)$ \Comment{Perform the ALM iteration}
    \State $\Delta\VEC{w}=\VEC{w}^{\ell+1}_{k+1}-\VEC{w}^{\ell+1}_{k}$\Comment{Update the solution step}
    \State $d_k\gets\texttt{distance}(\VEC{w}^{\ell+1}_{k+1},\VEC{w}^{\ell+1}_{k})$\Comment{Gets the distance}
    \EndFor
    \State $\delta L\gets\texttt{distance}(\VEC{w}^\ell_{i+1},\VEC{w}^{\ell+1}_N)$
    \State $\overline{\Delta L}\gets\texttt{distance}(\VEC{w}_N,\VEC{w}_0)$
    \Ensure $\{d^{\ell+1}_j\}_{j=0,...,N-1}$, $\{\VEC{w}^{\ell+1}_j\}_{j=0,...,N}$, $\overline{\Delta L}$, $\delta L$
  \end{algorithmic}
\end{algorithm}

The \texttt{correctQueueSerial} routine includes the \texttt{pop}, \texttt{submit}, and \texttt{collectSolutions} routines. These routines mainly involve read and write operations for the mappers defined in \cref{fig:APALM_datastructure}; hence, only a brief description is provided:
\begin{itemize}
\item $Q$, $\texttt{ID}$, $\Delta L$, $\VEC{w}^\ell_i$, $\VEC{w}^\ell_{i-1},\VEC{w}^\ell_{i+1}$$\gets\texttt{pop}(Q)$: Takes the first available interval from the queue $Q$ and returns a job \texttt{ID}, an interval length $\Delta L$, the start solution $\VEC{w}^\ell_i$, the previous solution $\VEC{w}^\ell_{i-1}$ and the next available solution $\VEC{w}^\ell_{i+1}$. It also updates the queue $Q$ internally by removing the current entry.
\item $Q$$\gets\texttt{submit}(\texttt{ID},\{d^{\ell+1}_k\}_{k=0,...,N-1},\{\VEC{w}^{\ell+1}_k\}_{k=0,...,N},\overline{\Delta L},\delta L,Q)$: Takes a job \texttt{ID}, the series of solutions $\{\VEC{w}^{\ell+1}_k\}_{k=0,...,N}$ and their distances $\{d^{\ell+1}_k\}_{k=0,...,N-1}$ and the distances $\overline{\Delta L}$ and $\delta L$. Using \cref{eq:APALM_total_error,eq:APALM_lower_error,eq:APALM_upper_error}, the errors are computed and solution intervals are added to the queue $Q$ if needed.
\item $\{\VEC{w}_i\}$, $\{s_i\}$$\gets$ \texttt{collectSolutions}: Based on the underlying mappers, the solutions of all levels are collected into $\{\VEC{w}_i\}$ and $\{s_i\}$.
\end{itemize}

\subsubsection*{Hybrid implementation}
The hybrid implementation of the APALM is referred to as the Adaptive Serial-Parallel Arc-Length Method (ASPALM), since it is a two-stage method with a serial initalization and a parallel correction. This concept is similar to the concept presented for the ASALM, but communication between the manager and worker processes is added so that the correction stage can be performed in parallel. To this end, the \texttt{correctQueueSerial} routine is re-defined into \texttt{correctQueueParallel} and communications between workers and the manager are defined. In the following, the global solution algorithm for the ASPALM is defined in \cref{alg:ASPALM}. As for the ASALM (see \cref{alg:ASALM_full}), the initialization is performed in serial by the \texttt{serialSolve} routine, and the maps are initialised, both by the manager process. As soon as queue $Q$ is established, the queue can be processed in parallel. Indeed, this implies that the worker processes are idle until the queue $Q$ is fully available.\\

\begin{algorithm}[h]
  \caption{Global ASPALM routine (\texttt{ASPALM}). The ASPALM first consists of a serial solve of the whole curve length domain, followed by a parallel evaluation of subintervals. This algorithm is specified for simple \texttt{manager}-\texttt{worker} parallelization; more advanced parallelization schemes, e.g., with multiple managers, are easily achieved.}
  \label{alg:ASPALM}
    \begin{algorithmic}[1]
      \Require $\Delta L$, $I$
      \If{\texttt{manager}}
      \State $\{\VEC{w}^\ell_i\}_{i=0}^{I},\{s^\ell_i\}_{i=0}^{I}\gets\texttt{serialSolve}(\Delta L,I)$ \Comment{See \cref{alg:serialSolve}}
      \State $Q\gets\texttt{initializeMap}(\{\VEC{w}^\ell_i\}_{i=0}^{I},\{s^\ell_i\}_{i=0}^{I})$
      \State $\{\VEC{w}^\ell_i\},\{s^\ell_i\}\gets\texttt{correctQueueParallel}(Q)$
      \Else
      \State Initialise \texttt{stop}$=$\texttt{false}
      \While{$\texttt{stop}=\texttt{false}$}
      \State $\texttt{stop}\gets\texttt{receiveStop}$
      \State $\texttt{workerCorrect}()$
      \EndWhile
      \EndIf
      \Ensure $\{\VEC{w}_i\}$, $\{s_i\}$
    \end{algorithmic}
\end{algorithm}

As seen in \cref{alg:ASPALM}, the manager process uses the \texttt{correctQueueParallel} routine; see \cref{alg:correctQueueParallel}. This routine is called by the manager process and sends and receives data to and from the workers, updates the queue, and heavily relies on communication functions as defined in \cref{tab:comms}. The communication function \texttt{sendMetaData} can be omitted for the ASPALM since it is primarily used in the APALM to distinguish between initiation and correction jobs. Furthermore, as described in \cref{tab:comms}, the routine \texttt{sendJob} can be called with and without the reference solution $\VEC{w}_{i+1}^{\ell}$, depending whether it is available (correction) or not (initiation). In the case of the ASPALM, all jobs that are popped from the queue $Q$ in the \texttt{correctQueueParallel} routine are by definition correction jobs.\\

As shown in \cref{alg:ASPALM}, the worker processes will perform the \texttt{workerCorrect} from \cref{alg:workerCorrect}. This algorithm contains the correction step for any interval that is received from the communications coming from the manager process. The \texttt{workerCorrect} is executed until a stop signal is received from the manager process. The latter is broadcast to all workers as soon as queue $Q$ is empty. \Cref{tab:comms} gives an overview of the communication functions that are used for communications between the manager and worker processes in the hybrid and parallel implementations.

\begin{table}
\caption{Required communications between manager and worker processes for the ASPALM and APALM.}
\label{tab:comms}
\centering
\begin{tabular}{p{0.15\linewidth}p{0.07\linewidth}p{0.07\linewidth}p{0.2\linewidth}p{0.35\linewidth}}
\toprule
\textbf{Send/Receive} & \textbf{From} & \textbf{To} & \textbf{Data} & \textbf{Description}\\
\midrule
\texttt{sendMetaData} \texttt{receiveMetaData} & Manager & Worker & $\texttt{ID}$, $\ell$ & Communicates meta-data between the manager and the worker processes. In this case, only the \texttt{ID} and the level $\ell$ are needed.\\
\midrule
\texttt{sendJob} \texttt{receiveJob} & Manager & Worker & $\texttt{ID}$, $\Delta L_0$, $\VEC{w}^\ell_i$, $\VEC{w}^\ell_{i-1}$, ,($\VEC{w}^\ell_{i+1}$) & Communicates information to perform the computation of an interval between the manager and the worker processes. The reference solution $\VEC{w}^\ell_{i+1}$ is optional since it is not available for initiation jobs.\\
\midrule
\texttt{sendData} \texttt{receiveData} & Worker & Manager & $W_j$, $\texttt{ID}$, $\{d^{\ell+1}_j\}_{j=0,...,N-1}$, $\{\VEC{w}^{\ell+1}_j\}_{j=0,...,N}$, $\overline{\Delta L}$, $\delta L$ & Communicates the data resulting from a sub-interval computation between the manager and the worker processes. The receive communication also provides the worker from whom the data is received.\\
\midrule
\texttt{sendStop} \texttt{receiveStop} & Manager & Worker (all) & Boolean & Communicates a stop signal between the manager and the worker processes.\\
\bottomrule
\end{tabular}
\end{table}

\subsubsection*{Parallel implementation}
Contrary to the serial and hybrid implementations, the fully parallel implementation does not work with a two-staged procedure of serial initialization and parallel correction. Instead, the fully parallel solve consists of one stage with a single queue consisting of initialization and correction jobs. For the Adaptive Parallel Arc-Length Method (APALM), the global routine is provided in \cref{alg:APALM_full}.\\

\begin{algorithm}[h]
  \caption{Global APALM routine (\texttt{APALM}). The APALM first consists of a serial solve of the whole curve length domain, followed by a parallel evaluation of subintervals. This algorithm is specified for simple \texttt{manager}-\texttt{worker} parallelization; more advanced parallelization schemes, e.g., with multiple managers, are easily achieved.}
  \label{alg:APALM_full}
  \begin{algorithmic}[1]
  \Require $\Delta L$, $I$
  \If{\texttt{manager}}
  \State $Q\gets\texttt{initializeQueue}(\Delta L, I)$
  \State $\{\VEC{w}^\ell_i\},\{s^\ell_i\}\gets\texttt{correctQueueParallel}(Q)$
  \State $\{\VEC{w}_i\}$, $\{s_i\}$$\gets\texttt{collectSolutions}$
  \Else
  \While{\texttt{true}}
  \State $\texttt{stop}\gets\texttt{receiveStop}$
  \If{stop}
  \State Break loop
  \EndIf
  \State $\texttt{ID},\ell\gets\texttt{receiveMetaData}$
  \If{$\ell=0$}
  \State $\texttt{workerInitiate}()$
  \Else
  \State $\texttt{workerCorrect}()$
  \EndIf
  \EndWhile
  \EndIf
  \Ensure $\{\VEC{w}_i\}$, $\{s_i\}$, $\{\xi_i\}$
  \end{algorithmic}
\end{algorithm}

As can be seen in \cref{alg:APALM_full}, the manager process in the APALM only initialises the queue using the \texttt{initializeQueue} routine, and it contains the \texttt{correctQueueParallel} routine. The former is not specified explicitly since it only initialises a map with zero points and allocates the maximum number of intervals $I$ as well as the interval length $\Delta L$. The \texttt{correctQueueParallel} routine is given in \cref{alg:correctQueueParallel}. It applies meta-data communication, see \cref{tab:comms}, to the worker processes so that the distinction based on the data level $\ell$ in \cref{alg:APALM_full} can be made by the workers. On the side of the worker processes, the meta-data is received, and if the level is equal to 0, an interval is initiated using \texttt{workerInitiate} (see \cref{alg:workerInitiate}), and if the level is larger than 0, an interval is corrected using \texttt{workerCorrect} (see \cref{alg:workerCorrect}); see \cref{alg:APALM_full}. Inside the \texttt{workerInitiate} and \texttt{workerCorrect} routines, communications from the worker to the manager process (see \cref{tab:comms}) are included.

\begin{algorithm}[h]
  \caption{The \texttt{correctQueueParallel} routine, accompanied by the \texttt{workerCorrect} routine from \cref{alg:workerCorrect} and communication functions defined in \cref{tab:comms}. This routine takes the queue $Q$ and assigns jobs from the queue to the available workers. Then, while the queue $Q$ is non-empty, data is communicated to and from the workers, and solutions are submitted. Note that the \texttt{pop} and \texttt{submit} routines are equivalent to the ones in \cref{alg:correctQueueSerial}.}
  \label{alg:correctQueueParallel}
  \begin{algorithmic}[1]
  \Require $Q$
  \State Initialise a pool of worker processes $W=\{W_j,\:j=1,...,N_{\text{workers}}\}$
  \While{$Q\neq\emptyset$ and $W\neq\emptyset$}
  \State $Q$, $\texttt{ID}$, $\Delta L$, $\VEC{w}^\ell_i$, $\VEC{w}^\ell_{i-1}$, $\VEC{w}^\ell_{i+1}$$\gets\texttt{pop}(Q)$ \Comment{See line 2 of \cref{alg:correctQueueSerial}}
  \State $\texttt{sendStop}(\texttt{false})$
  \State $\texttt{sentMetaData}(\texttt{ID},\ell)$
  \State $\texttt{sendJob}(\texttt{ID},\Delta L_0,\VEC{w}^\ell_i,\VEC{w}^\ell_{i-1},\VEC{w}^\ell_{i+1},W_j)$
  \State Remove $Q_i$ from $Q$ and $W_j$ from $W$
  \EndWhile
  \State \Comment{Send jobs to workers when they are available}
  \While{$\vert W \vert \neq N_\text{workers}$}
  \State $W_j,\texttt{ID},\{d^{\ell+1}_j\}_{j=0,...,N-1}, \{\VEC{w}^{\ell+1}_j\}_{j=0,...,N}, \overline{\Delta L}, \delta L\gets\texttt{receiveDataWorker2Manager}$
  \State $Q$$\gets\texttt{submit}(\texttt{ID},\{d^{\ell+1}_j\}_{j=0,...,N-1},\{\VEC{w}^{\ell+1}_j\}_{j=0,...,N},\overline{\Delta L}, \delta L,Q)$ \Comment{See line 4 of \cref{alg:correctQueueSerial}}
  \State Add $W_j$ to $W$
  \While{$Q\neq\emptyset$ and $W\neq\emptyset$}
  \State $Q$, $\texttt{ID}$, $\Delta L$, $\VEC{w}^\ell_i$, $\VEC{w}^\ell_{i-1}$, $\VEC{w}^\ell_{i+1}$$\gets\texttt{pop}(Q)$
  \State $\texttt{sendStop}(\texttt{false})$ \Comment{The worker always expects a stop signal, now it is false.}
  \State $\texttt{sentMetaData}(\texttt{ID},\ell)$
  \State $\texttt{sendJob}(\texttt{ID},\Delta L_0,\VEC{w}^\ell_i,\VEC{w}^\ell_{i-1},\VEC{w}^\ell_{i+1},W_j)$
  \State Remove $Q_i$ from $Q$ and $W_j$ from $W$
  \EndWhile
  \EndWhile
  \State $\texttt{sendStopManager2All}(\texttt{true})$
  \Ensure $\{\VEC{w}_i\}$, $\{s_i\}$
  \end{algorithmic}
\end{algorithm}

\begin{algorithm}[h]
  \caption{Solve routine for a worker (\texttt{workerInitiate}). This routine performs the initiation steps (see \cref{alg:initiate}) on jobs received from the manager, until a stop signal is received. More information on the communication functions can be found in \cref{tab:comms}}
  \label{alg:workerInitiate}
  \begin{algorithmic}[1]
  \Require
  \State $(\texttt{ID},\Delta L_0,\VEC{w}_i,\VEC{w}_{i-1},\VEC{w}_\text{ref})\gets\texttt{receiveJob}$
  \State $\VEC{w}^\ell_{k+1},\Delta s^\ell_{k+1} \gets $\texttt{initiate}$(\VEC{w}^\ell_{k},\Delta\VEC{w}_{k-1}^\ell,\Delta L)$
  \State $\texttt{sendDataWorker2Manager}(\texttt{ID},\{d_j\}^{\ell+1}_{j=0,...,N-1}, \{\VEC{w}_j\}^{\ell+1}_{j=0,...,N}, \overline{\Delta L})$
  \Ensure
  \end{algorithmic}
\end{algorithm}

\begin{algorithm}[h]
  \caption{Solve routine for a worker (\texttt{workerCorrect}). This routine performs correction steps (see \cref{alg:correct}) on jobs received from the manager, until a stop signal is received. More information on the communication functions can be found in \cref{tab:comms}}
  \label{alg:workerCorrect}
  \begin{algorithmic}[1]
  \Require
  \State $(\texttt{ID},\Delta L_0,\VEC{w}_i,\VEC{w}_{i-1},\VEC{w}_\text{ref})\gets\texttt{receiveJob}$
  \State $\{d^{\ell+1}_k\}_{k=0,...,N-1}$, $\{\VEC{w}^{\ell+1}_k\}_{k=0,...,N}$, $\overline{\Delta L}$, $\delta L$$\gets\texttt{correct}(N,\Delta L_0, \VEC{w}^\ell_{i}, \VEC{w}^\ell_{i-1}, \VEC{w}^\ell_{i+1})$
  \State $\texttt{sendDataWorker2Manager}(\texttt{ID},\{d_j\}^{\ell+1}_{j=0,...,N-1}, \{\VEC{w}_j\}^{\ell+1}_{j=0,...,N}, \overline{\Delta L})$
  \Ensure
  \end{algorithmic}
\end{algorithm}

\subsection{Arc-length Exploration}\label{subsec:APALM_exploration}
To enable multi-branch parallelization of APALM, small modifications are required to the data structure and the algorithms presented in \cref{subsec:APALM_datastructure,subsec:APALM_algorithms}. The easiest multi-branch parallelization is achieved by identifying branch switches only in the serial solve, such that the initialization of the APALM can be done across different branches. In this case, the serial solve is performed on the main branch, and any bifurcation point is stored such that a restart can be performed from this point; see \cite{Wouters2019} for details. As soon as such a bifurcation point is identified, a branch switch can be performed, and a new serial solve can be started from that point. As a result, a series of solutions $\{\VEC{w}_i\}_{i=0}^{I}$ is computed for each branch. Similar to the single-branch case, a data structure and a queue can be initialised using \cref{alg:initialise} per branch. Depending on the parallel configuration, the queues $Q^b$, $b=0,...,n_{\text{branches}}$ of all branches can be treated separately by multiple manager processes, or they can be combined into one large queue $Q$ and handled by one single manager process. In the latter case, each job will also contain a branch identifier to refer to the corresponding data structure.\\

The advantages of the above approach combining multi-branch and within-branch parallelization using the APALM are that the extension from a single-branch APALM to a multi-branch APALM is straightforward. The disadvantage, however, is that the identification of bifurcations is only taken into account in the serial solve step; hence, any bifurcations that are identified in the parallel solve will not be taken into account. A remedy would be to rebuild the map and the data structure on the manager process as soon as a worker process finds a bifurcation point; this requires all active workers to terminate.

%% file: Figures/DataStructure.pgf
\begin{tikzpicture}
\def\N{7}
\def\L{7}
\def\dh{0.1}
\def\hX{0}
\def\hS{0.8}
\def\hL{1.6}
\def\hU{2.4}
\def\hUg{3.2}
\def\dI{0.3}

\def\hI{-1.5}

\pgfmathsetmacro\min{0-0.5)}
\pgfmathsetmacro\max{\L+0.5}

\draw[-latex] (\min,\hS) -- (\max,\hS)node [right] {$s$};
\draw[-latex] (\min,\hX) -- (\max,\hX)node [right] {$\xi$};
\draw[-] (\min,\hL) -- (\max,\hL)node [right] {$l$};
\draw[-] (\min,\hU) -- (\max,\hU)node [right] {$\vb*{w}$};
\draw[-] (\min,\hUg) -- (\max,\hUg)node [right] {$\vb*{w}'$};

\foreach[evaluate={\kmm=int(\k-1)}] \k in {0,...,\N}
{
	\pgfmathsetmacro\x{\k*\L/\N}
	\pgfmathsetmacro\Xi{\k*\L/\N/\L}
	\draw (\x,\hX+\dh) -- node[midway](X\k){} (\x,\hX-\dh);
	\draw (\x,\hS+\dh) -- node[midway](S\k){} (\x,\hS-\dh);
	\draw (\x,\hL+\dh) -- node[midway](L\k){} (\x,\hL-\dh);
	\draw (\x-\dh,\hU-\dh) rectangle (\x+\dh,\hU+\dh);
	\node[] (U\k) at (\x,\hU){};
	\draw[dashed] (\x-\dh,\hUg-\dh) rectangle (\x+\dh,\hUg+\dh);
	\node[] (Ug\k) at (\x,\hUg){};
	\ifthenelse{\k>0} {\draw[dashed,latex-](Ug\kmm) -- (U\k);}{\draw[dashed,latex-](Ug\k) -- (U\k);}
}

\draw[-latex](X1) to[in=-60,out=60] node[inner sep=1,midway,above right]{$\mathcal{S}(\xi_i)$} (S1);
\draw[latex-](X1) to[in=240,out=-240] node[inner sep=1,midway,above left]{$\Xi(s_i)$} (S1);
\draw[-latex](X2) to[in=-60,out=60] node[inner sep=1,midway,below right]{$\mathcal{L}(\xi_j)$} (L2);
\draw[-latex](X4) to[in=-60,out=60] node[inner sep=1,midway,below right]{$\mathcal{U}(\xi_k)$} (U4);
\draw[-latex](X4) to[in=240,out=-240] node[inner sep=1,midway,below left]{$\mathcal{U}'(\xi_k)$} (Ug4);

\foreach[evaluate={\Ipp=int(\I+1)},count=\c from 0] \I in {2,4}
{
	\draw[ultra thick] (X\I.center)--(X\Ipp.center) node[midway,below]{$A_\c$};
}

\draw[ultra thick,densely dashed] (X0.center)--(X1.center) node[midway,below]{$Q_{m}$};
\draw[ultra thick,densely dashed] (X3.center)--(X4.center) node[midway,below]{$Q_{n}$};

\def\Nsub{3}
\foreach[evaluate={\Ipp=int(\I+1)}] \I in {6}
{
	\pgfmathsetmacro\Nmm{int(\Nsub-1)}
	\pgfmathsetmacro\dhn{\dh*1cm/1pt};
	\foreach[evaluate={\kmm=int(\k-1)}] \k in {1,...,\Nmm}
{
		\draw[red] let \p0 = (X\I), \p1 = (X\Ipp) in (\x0+\k*\x1/\Nsub-\k*\x0/\Nsub,\hX-\dh/2) --  (\x0+\k*\x1/\Nsub-\k*\x0/\Nsub,\hX+\dh/2);
		\draw[red] let \p0 = (S\I), \p1 = (S\Ipp) in (\x0+\k*\x1/\Nsub-\k*\x0/\Nsub,\hS-\dh/2) --  (\x0+\k*\x1/\Nsub-\k*\x0/\Nsub,\hS+\dh/2);
		\draw[red] let \p0 = (L\I), \p1 = (L\Ipp) in (\x0+\k*\x1/\Nsub-\k*\x0/\Nsub,\hL-\dh/2) --  (\x0+\k*\x1/\Nsub-\k*\x0/\Nsub,\hL+\dh/2);
		\draw[red] let \p0 = (U\I), \p1 = (U\Ipp) in (\x0+\k*\x1/\Nsub-\k*\x0/\Nsub-\dhn/2,\hU-\dh/2) rectangle (\x0+\k*\x1/\Nsub-\k*\x0/\Nsub+\dhn/2,\hU+\dh/2);
	\path let \p0 = (U\I), \p1 = (U\Ipp) in node (Un\k) at (\x0+\k*\x1/\Nsub-\k*\x0/\Nsub,\hU){};
		\draw[red,densely dashed] let \p0 = (Ug\I), \p1 = (Ug\Ipp) in (\x0+\k*\x1/\Nsub-\k*\x0/\Nsub-\dhn/2,\hUg-\dh/2) rectangle  (\x0+\k*\x1/\Nsub-\k*\x0/\Nsub+\dhn/2,\hUg+\dh/2);
	\path let \p0 = (Ug\I), \p1 = (Ug\Ipp) in node (Ugn\k) at (\x0+\k*\x1/\Nsub-\k*\x0/\Nsub,\hUg){};

	\ifthenelse{\k>1} 
	{
		\draw[dashed,latex-,red](Ugn\kmm) -- (Un\k);
	}
	{
		\draw[dashed,latex-,red](Ug\I) -- (Un\k);
	}
}

}
\end{tikzpicture}

%% file: Sections/Results.tex
In this section, the APALM scheme is demonstrated on a series of benchmark problems. The first two benchmark problems are structural analysis problems with limit-point instabilities and complex collapsing mechanisms involving strongly curved solution paths. The third benchmark is a buckling problem containing a bifurcation with multiple branches to illustrate the concept of the APALM in a multi-branch setting. All benchmark problems are computed using isogeometric Kirchhoff-Love shell elements based on the works \cite{Kiendl2009,Kiendl2015,Verhelst2021}. Furthermore, a scaling test with respect to the number of worker processes is performed for all benchmark problems. Here, a scaling analysis of the ASPALM is performed to demonstrate the relative computational costs of the serial initialization phase compared to the parallel correction phase. Furthermore, a scaling analysis of the APALM is performed to show the advantage of the fully parallel APALM scheme over the two-stage ASPALM scheme. For the scaling tests, the communications from \cref{tab:comms} are performed using the Message Passing Interface (MPI), and they are performed on the Delft High Performance Computing Centre (DHPC) \cite{DHPC2022} with Intel XEON E5-6248R 24C 3.0GHz nodes with 96GB of memory per CPU. The code will be made available within the Geometry + Simulation modules \cite{Juttler2014} upon publication.

\subsection{Collapse of a Shallow Roof}
The first benchmark problem involves a shallow roof subject to a point load at the midpoint. The roof is discretized with $4\times4$ NURBS elements of degree $3$. The roof is composed of a lay-up of composites with material properties as presented in \cref{fig:APALM_roof_case}, inspired by \cite{Leonetti2019}. It is modelled using isogeometric Kirchhoff-Love shell elements \cite{Kiendl2009} supporting composite laminates \cite{Herrema2019}. A Crisfield ALM is used with an initial arc length of 30 and a scaling parameter of $\Psi=1$. The tolerance of the APALM is set to $\varepsilon_l = \varepsilon_u = 10^{-2}$. \RevTHREE{This tolerance implies that intervals are marked for refinement when the traversed length is deviates than $1\%$ of the original interval length. A smaller tolerance would imply that more elements are marked for refinement and that the corrections are performed up to lower levels.} The material, and load parameters for this benchmark can be found in \cref{fig:APALM_roof_case}. Reference solutions are obtained using a serial arc-length method with a sufficiently fine increment size.\\

\begin{figure}[h]
\centering
\input{Figures/Roof.pgf}
\caption{The problem definition for the benchmark of the collapsing roof with length $L=508\:[\text{mm}]$, with radius $R=2540\:[\text{mm}]$, an angle $\theta=10\:[\text{rad}]$ with a thickness of $t=6.35\:[\text{mm}]$. The boundaries $\Gamma_1$ and $\Gamma_3$ have fixed displacements, and the other sides are free. The material is modelled using a Saint-Venant Kirchhoff laminate with $E_{11}=3300\:[\text{N}/\text{mm}^2]$, $E_{22}=E_{33}=1100\:[\text{N}/\text{mm}^2]$, $G_{12}=G_{13}=660\:[\text{N}/\text{mm}^2]$, $E_{23}=440\:[\text{N}/\text{mm}^2]$ and $\nu_{12}=\nu_{13}=\nu_{23}=0.25\:[-]$ and with lay-up angles $[0/90/0]^\circ$. The load is variable with a magnitude $P=10\:[\text{N}]$ and magnification factor $\lambda$.}
\label{fig:APALM_roof_case}
\end{figure}

\Cref{fig:APALM_roof} provides the results of the APALM applied to the collapse of the shallow roof. As can be seen from this figure, the serial computation provides a coarse estimate of the reference curve. Especially on the first limit point (between $\lambda\times w_A\in[-15,-8]\times[20,25]$), the data is sparse, similar to the collapse itself (see inset in \cref{fig:APALM_roof}). The results of the APALM show that a lot of refinements are needed to represent the collapsing behaviour correctly in the region of the inset in \cref{fig:APALM_roof}. These regions do not necessarily involve extremely curved paths in the axes of \cref{fig:APALM_roof}, but the solutions $\VEC{w}$ are most likely curved in the higher-dimensional solution space.\\

\begin{figure}[h]
\def\minXi{7}
\def\maxXi{16}
\def\minYi{22}
\def\maxYi{28}
\def\minXii{10}
\def\maxXii{15}
\def\minYii{8}
\def\maxYii{14}
\def\minXiii{12}
\def\maxXiii{17}
\def\minYiii{-25}
\def\maxYiii{-19}
\centering
\begin{subfigure}[t]{\linewidth}
    \centering
    \begin{tikzpicture}
        \begin{axis}
        [
            xlabel=Load factor $\lambda$,
            ylabel=Mid-point displacement $u_A$,
            legend pos = south west,
            width=0.7\linewidth,
            height=0.3\textheight,
            grid=major,
            legend columns=2,
            restrict x to domain = 0:27
        ]
        \addplot+[solid,no markers,black] table[header=true,x expr = -\thisrowno{2},y index = {4}, col sep = comma]{Results/roof_ref.csv};
        \addlegendentry{Reference}
        \addplot+[only marks,mark=o,mark size=2,black] table[header=true,x expr = -\thisrowno{2},y index = {4}, col sep = comma]{Results/roof_serial.csv};
        \addlegendentry{Serial}
        \pgfplotsset{cycle list shift=-2}
        \foreach \level in {1,...,6}
        {
            \addplot+[only marks] table[header=true,x expr = -\thisrowno{2},y index = {4}, col sep = comma,restrict expr to domain={\thisrowno{6}}{\level:\level}]{Results/roof_parallel.csv};
            \addlegendentryexpanded{Level \level}
        }

        \draw[gray,dashed] (axis cs: \minXi,\minYi) -- (axis cs: \maxXi,\minYi) -- (axis cs: \maxXi,\maxYi) node[above right, inner sep=1pt] {A} -- (axis cs: \minXi,\maxYi) --cycle;

        \draw[gray,dashed] (axis cs: \minXii,\minYii) -- (axis cs: \maxXii,\minYii) -- (axis cs: \maxXii,\maxYii) node[above right,inner sep=1pt] {B} -- (axis cs: \minXii,\maxYii) --cycle;

        \draw[gray,dashed] (axis cs: \minXiii,\minYiii) -- (axis cs: \maxXiii,\minYiii) -- (axis cs: \maxXiii,\maxYiii) node[above right,inner sep=1pt] {C} -- (axis cs: \minXiii,\maxYiii) --cycle;
     \end{axis}
     \end{tikzpicture}
     \begin{tikzpicture}
         \begin{groupplot}
        [
            group style={
                group size=1 by 3,
                xlabels at=edge bottom,
                yticklabels at=edge right,
                vertical sep=0.025\textheight,
                horizontal sep=0.025\linewidth,
            },
            xlabel=Load factor $\lambda$,
            xlabel style={color=gray},
            x tick label style={color = gray},
            ylabel style={color=gray},
            y tick label style={color = gray},
            legend pos = north east,
            width=0.3\linewidth,
            height=0.13\textheight,
            grid=major,
            legend columns=2,
            restrict x to domain = 0:27,
            ylabel={\phantom{.}},
        ]
        \nextgroupplot[
            xmin = \minXi,
            xmax = \maxXi,
            ymin = \minYi,
            ymax = \maxYi,
            draw = gray,
        ]
        \addplot+[solid,no markers,black] table[header=true,x expr = -\thisrowno{2},y index = {4}, col sep = comma]{Results/roof_ref.csv};
        \addplot+[only marks,mark=o,mark size=2,black] table[header=true,x expr = -\thisrowno{2},y index = {4}, col sep = comma]{Results/roof_serial.csv};
        \pgfplotsset{cycle list shift=-2}
        \foreach \level in {1,...,6}
        {
            \addplot+[only marks] table[header=true,x expr = -\thisrowno{2},y index = {4}, col sep = comma,restrict expr to domain={\thisrowno{6}}{\level:\level}]{Results/roof_parallel.csv};
        }
        \node[gray] at (rel axis cs: 0.9,0.85){A};

        \nextgroupplot[
            xmin = \minXii,
            xmax = \maxXii,
            ymin = \minYii,
            ymax = \maxYii,
            draw = gray,
        ]
        \addplot+[solid,no markers,black] table[header=true,x expr = -\thisrowno{2},y index = {4}, col sep = comma]{Results/roof_ref.csv};
        \addplot+[only marks,mark=o,mark size=2,black] table[header=true,x expr = -\thisrowno{2},y index = {4}, col sep = comma]{Results/roof_serial.csv};
        \pgfplotsset{cycle list shift=-2}
        \foreach \level in {1,...,6}
        {
            \addplot+[only marks] table[header=true,x expr = -\thisrowno{2},y index = {4}, col sep = comma,restrict expr to domain={\thisrowno{6}}{\level:\level}]{Results/roof_parallel.csv};
        }
        \node[gray] at (rel axis cs: 0.9,0.85){B};

        \nextgroupplot[
            xmin = \minXiii,
            xmax = \maxXiii,
            ymin = \minYiii,
            ymax = \maxYiii,
            draw = gray,
        ]
        \addplot+[solid,no markers,black] table[header=true,x expr = -\thisrowno{2},y index = {4}, col sep = comma]{Results/roof_ref.csv};
        \addplot+[only marks,mark=o,mark size=2,black] table[header=true,x expr = -\thisrowno{2},y index = {4}, col sep = comma]{Results/roof_serial.csv};        \pgfplotsset{cycle list shift=-2}
        \foreach \level in {1,...,6}
        {
            \addplot+[only marks] table[header=true,x expr = -\thisrowno{2},y index = {4}, col sep = comma,restrict expr to domain={\thisrowno{6}}{\level:\level}]{Results/roof_parallel.csv};
        }
        \node[gray] at (rel axis cs: 0.9,0.85){C};
      \end{groupplot}
 \end{tikzpicture}
\end{subfigure}
\caption{Results of the collapsing roof. The figure on the left indicates the full solution path, and the figures on the right depict the insets indicated in the left figure. The reference and serial solutions are represented by the solid line and the black markers, respectively. The solutions computed by the APALM are indicated per level. The simulation is performed with a tolerance of $\varepsilon_l=\varepsilon_u=10^{-2}$ and an increment length of $\Delta L=30$.}
\label{fig:APALM_roof}
\end{figure}

In order to assess the parallelization of the APALM in this example, the test from \cref{fig:APALM_roof} is performed with an increasing number of workers using the ASPALM and the APALM schemes. The results in \cref{tab:roof_dL30} show that for the computation with arc-length parameter $\Delta L=30$, the parallel correction step of the ASPALM scales optimally (i.e., with a factor 2) up to around 8 workers, after which the scalability decreases and the parallel correction phase takes around 15\% of the total computational time. Using the APALM scheme, the total computational time is decreased compared to the ASPALM scheme, and parallel corrections can be started as soon as the first interval has been initialized. The computational time of the APALM stagnates around 4 workers, at a computational time similar to the serial initialization time for the ASPALM method, showing that adaptive parallel corrections can be performed without significantly more computational costs compared to a serial arc-length method. When the number of intervals is increased, e.g., by decreasing the arc-length parameter to $\Delta L=2.5$ (\cref{tab:roof_dL2.5}), it can be seen that the parallel stage of the ASPALM scales up to a higher number of workers, in this case 64, up to the point that it takes around 5\% of the total computational time for 256 workers. The APALM again provides computational times similar to a serial ALM without corrections. The improved scalability is explained by the fact that the queue is in general longer; therefore, the time that workers are idle waiting for the last job to be finished is smaller relative to the total computational time.

\begin{table}[h]
    \centering
    \caption{Computational time in $[\text{s}]$ for the benchmark of the collapsing roof for the ASPALM and APALM for different numbers of worker processes. The times for the ASPALM are presented for the serial initialization and the parallel correction phases, and the sum of the two is given as the total computational time. The numbers in the \textit{Serial} column should theoretically be the same, but they provide a representation of the variation in the time measurements. The results are presented for simulations with increment lengths $\Delta L=30$ (\subref{tab:roof_dL30}) and $\Delta L=2.5$ (\subref{tab:roof_dL2.5}), and the italic row with $0$ workers denotes the ASALM method.}
    \label{tab:roof}
    \begin{subtable}[t]{0.45\linewidth}
    \centering
    \footnotesize
    \caption{$\Delta L=30$}
    \label{tab:roof_dL30}
    \pgfkeys{/pgf/number format/.cd,fixed,precision=1,fixed zerofill}
    \begin{filecontents*}{Results/Roof_times_50_steps_table.csv}
0,115.747,195.303,311.05,287.077
1,119.155,209.011,328.166,318.762
2,113.995,100.781,214.776,162.405
4,109.491,46.1302,155.6212,115.772
8,115.019,27.0378,142.0568,115.881
16,115.133,17.8133,132.9463,116.268
32,114.918,15.8918,130.8098,113.036
64,114.537,13.2838,127.8208,115.98
\end{filecontents*}
\pgfplotstabletypeset[
        header=false,
        every head row/.style={
            output empty row,
            before row={%
                \toprule
                \textbf{Workers} & \multicolumn{3}{c}{\textbf{ASPALM}} & \textbf{APALM}\\
                    \#  & \multicolumn{1}{c@{\hspace*{\tabcolsep}\makebox[0pt]{+}}}{Serial} & \multicolumn{1}{c@{\hspace*{\tabcolsep}\makebox[0pt]{=}}}{Parallel} & Total       & Parallel \\
            },
                        after row=\cmidrule(lr){0-0}\cmidrule(lr){2-4}\cmidrule(lr){5-5},
        },
        every last row/.style={after row=\bottomrule},
        columns/0/.style ={column name={},fixed,precision=0,fixed zerofill},
        col sep=comma,
        highlightrow={1},
    ]
    {Results/Roof_times_50_steps_table.csv}
    \end{subtable}
    \hfill
    \begin{subtable}[t]{0.45\linewidth}
    \centering
    \footnotesize
    \caption{$\Delta L=2.5$}
    \label{tab:roof_dL2.5}
    \pgfkeys{/pgf/number format/.cd,fixed,precision=1,fixed zerofill}
    \begin{filecontents*}{Results/Roof_times_600_steps_table.csv}
0,507.165,1778.09,2285.255,2187.09
1,500.457,1757.7,2258.157,2310.15
2,447.541,835.329,1282.87,1114.02
4,493.377,449.409,942.786,558.066
8,496.782,223.248,720.03,453.913
16,503.266,112.967,616.233,483.588
32,493.209,58.1247,551.3337,510.879
64,504.244,29.1852,533.4292,498.33
128,501.032,20.2423,521.2743,494.655
256,505.501,18.7973,524.2983,509.607
\end{filecontents*}
\pgfplotstabletypeset[
        header=false,
        every head row/.style={
            output empty row,
            before row={%
                \toprule
                \textbf{Workers} & \multicolumn{3}{c}{\textbf{ASPALM}} & \textbf{APALM}\\
                    \#  & \multicolumn{1}{c@{\hspace*{\tabcolsep}\makebox[0pt]{+}}}{Serial} & \multicolumn{1}{c@{\hspace*{\tabcolsep}\makebox[0pt]{=}}}{Parallel} & Total       & Parallel \\
            },
                        after row=\cmidrule(lr){0-0}\cmidrule(lr){2-4}\cmidrule(lr){5-5},
        },
        every last row/.style={after row=\bottomrule},
        columns/0/.style ={column name={},fixed,precision=0,fixed zerofill},
        col sep=comma,
        highlightrow={1},
    ]
    {Results/Roof_times_600_steps_table.csv}
    \end{subtable}
\end{table}

\subsection{Collapse of a Truncated Cone}
The second benchmark example is based on \cite{Verhelst2021} and involves the collapsing behaviour of a truncated cone with a hyperelastic Mooney-Rivlin material model. This benchmark is based on \cite{Basar1998}, but in the work of \cite{Verhelst2021}, the full collapsing behaviour was revealed using arc-length methods. The geometry, material, and load specifications can be found in \cref{fig:APALM_frustrum_case}.\\

The truncated cone is modelled using a quarter geometry using symmetry conditions to represent the axisymmetry as used in the original case of \cite{Basar1998}. The geometry is modelled with $32$ NURBS elements of degree $2$ over the height. Further, an initial arc length of $0.5$ is used, and the scaling factor is $\Psi=0$. The top boundary $\Gamma_2$ is free, and on the bottom boundary $\Gamma_4$, all displacements are fixed. The other boundaries have symmetric boundary conditions. The governing material model is an incompressible Mooney-Rivlin material model with a strain energy density function (with a slight abuse of notation)
\begin{equation}
 \Psi(\vb{C})=\frac{c_1}{2}\left(I_1-3\right) + \frac{c_2}{2}\left(I_2-3\right),
\end{equation}
with $I_1$ and $I_2$ the first and second invariants of the deformation tensor $\vb{C}=\VEC{F}^\top\VEC{F}$. More information on the problem set-up and the material models can be found in \cite{Verhelst2021}. The reference results are obtained from a serial ALM computation with a sufficiently small increment size.\\

The results of the collapsing truncated cone problem are presented in \cref{fig:APALM_frustrum}. As seen in this picture, the serial initialization provides a coarse approximation of the path but leaves out details, e.g., the rotated ``S''-shaped curve in the inset in \cref{fig:APALM_frustrum}. From the results, it is clear that the APALM focuses its refinement on the curved parts of the path and reveals the ``S''-shaped curve among other features of the path.\\

\begin{figure}[h]
\centering
\input{Figures/Frustrum.pgf}
\caption{The problem definition for the benchmark of the collapsing truncated cone with inner radii $R_1=1\:[\text{m}]$ and $R_2=2\:[\text{m}]$ and height $H=1\:[\text{m}]$. The thickness of the cone is $t=0.1\:[\text{m}]$. The cone is modelled by using a quarter of the geometry, using symmetry conditions on $\Gamma_1$ and $\Gamma_3$. The displacements at the bottom boundary ($\Gamma_4$) are fixed, and on the top boundary, a variable line load is applied and is variable with magnitude $p=1\:[\text{N}/\text{mm}]$ and magnification factor $\lambda$. The material of the cone is modelled using an incompressible Mooney-Rivlin model with parameters $\mu=c_1+c_2=4.225\:[\text{N}/\text{mm}^2]$, $c_1/c_2=7$.}
\label{fig:APALM_frustrum_case}
\end{figure}

\begin{figure}[h]
\def\minX{0}
\def\maxX{1.1}
\def\minY{-0.02}
\def\maxY{0.055}
\centering
\begin{subfigure}[t]{\linewidth}
\centering
\begin{tikzpicture}
        \begin{axis}
        [
            xlabel=Load factor $\lambda$,
            ylabel=Mid-point displacement $u_A$,
            legend pos = north east,
            ymin = -0.08,
            ymax = 0.08,
            width=0.6\linewidth,
            height=0.3\textheight,
            grid=major,
            legend columns=2,
        ]
        \addplot+[solid,no markers,black] table[header=true,x index = {0},y index = {1}, col sep = comma]{Results/frustrum_ref.csv};
        \addlegendentry{Reference}
        \addplot+[only marks,mark=o,mark size=2,black] table[header=true,x expr = -\thisrowno{2},y index = {4}, col sep = comma]{Results/frustrum_serial.csv};
        \addlegendentry{Serial}

        \pgfplotsset{cycle list shift=-2}
        \foreach \level in {1,...,6}
        {
            \addplot+[only marks] table[header=true,x expr = -\thisrowno{2},y index = {4}, col sep = comma,restrict expr to domain={\thisrowno{6}}{\level:\level}]{Results/frustrum_parallel.csv};
            \addlegendentryexpanded{Level \level}
        }

        \draw[gray,dashed] (axis cs: \minX,\minY) -- (axis cs: \maxX,\minY) -- (axis cs: \maxX,\maxY) -- (axis cs: \minX,\maxY) --cycle;
        \end{axis}
    \end{tikzpicture}
    \begin{tikzpicture}
        \begin{axis}
        [
            xlabel={Load factor $\lambda$},
            xlabel style={color=gray},
            x tick label style={color = gray},
            ylabel style={color=gray},
            y tick label style={color = gray},
            legend pos = north east,
            xmin = \minX,
            xmax = \maxX,
            ymin = \minY,
            ymax = \maxY,
            width=0.4\linewidth,
            height=0.3\textheight,
            grid=major,
            draw = gray,
            legend columns=2,
        ]
        \addplot+[solid,no markers,black] table[header=true,x index = {0},y index = {1}, col sep = comma]{Results/frustrum_ref.csv};
        \addplot+[only marks,mark=o,mark size=2,black] table[header=true,x expr = -\thisrowno{2},y index = {4}, col sep = comma]{Results/frustrum_serial.csv};
        \pgfplotsset{cycle list shift=-2}
        \foreach \level in {1,...,6}
        {
            \addplot+[only marks] table[header=true,x expr = -\thisrowno{2},y index = {4}, col sep = comma,restrict expr to domain={\thisrowno{6}}{\level:\level}]{Results/frustrum_parallel.csv};
        }
      \end{axis}
 \end{tikzpicture}
\end{subfigure}
\caption{Results of the collapsing truncated cone. The figure on the left depicts the full solution path, and the figure on the right depicts the inset indicated in the left figure. The reference and serial solutions are represented by the solid line and the black markers, respectively. The solutions computed by the APALM are indicated per level. The simulation is performed with a tolerance of $\varepsilon_l=\varepsilon_u=10^{-2}$ and an increment length of $\Delta L=0.5$.}
\label{fig:APALM_frustrum}
\end{figure}

Similar to the collapse of the roof, a scaling analysis of the parallel evaluations is performed. The results in \cref{tab:frustrum_dL0.5} verify that, as for the benchmark example with the collapsing roof, the scalability of the parallel correction phase scales optimally up to 8 workers, where the parallel correction phase takes around 15 \% of the total computational time when using 64 workers. When using the APALM scheme, the collapsing cone also shows that the computation times of the APALM are similar to the times needed for serial initialization, in other words, a classical ALM without adaptive corrections. When the number of intervals is increased, i.e., when the arc-length parameter is decreased to $\Delta L=0.0625$ (\cref{tab:frustrum_dL0.0625}), the scalability of the parallel phase of the ASPALM and of the full APALM reaches further, up to 64 workers.

\begin{table}[h]
    \centering
    \caption{Computational time in $[\text{s}]$ for the benchmark of the collapsing truncated cone for the ASPALM and APALM for different numbers of worker processes. The times for the ASPALM are presented for the serial initialization and the parallel correction phases, and the sum of the two is given as the total computational time. The numbers in the \textit{Serial} column should theoretically be the same, but they provide a representation of the variation in the time measurements. The results are presented for simulations with increment lengths $\Delta L=0.5$ (\subref{tab:frustrum_dL0.5}) and $\Delta L=0.0625$ (\subref{tab:frustrum_dL0.0625}), and the italic row with $0$ workers denotes the ASALM method.}
    \label{tab:frustrum}
    \begin{subtable}[t]{0.45\linewidth}
    \centering
    \footnotesize
    \caption{$\Delta L=0.5$}
    \label{tab:frustrum_dL0.5}
    \pgfkeys{/pgf/number format/.cd,fixed,precision=1,fixed zerofill}
    \begin{filecontents*}{Results/Frustrum_times_80_steps_table.csv}
0,160.233,244.015,404.248,436.56
1,162.537,247.17,409.707,424.762
2,169.495,130.064,299.559,207.121
4,170.582,68.0801,238.6621,172.857
8,162.612,43.017,205.629,160.54
16,175.291,32.015,207.306,173.276
32,175.476,27.3013,202.7773,170.784
64,170.063,23.3084,193.3714,169.506
\end{filecontents*}
\pgfplotstabletypeset[
        header=false,
        every head row/.style={
            output empty row,
            before row={%
                \toprule
                \textbf{Workers} & \multicolumn{3}{c}{\textbf{ASPALM}} & \textbf{APALM}\\
                    \#  & \multicolumn{1}{c@{\hspace*{\tabcolsep}\makebox[0pt]{+}}}{Serial} & \multicolumn{1}{c@{\hspace*{\tabcolsep}\makebox[0pt]{=}}}{Parallel} & Total       & Parallel \\
            },
                        after row=\cmidrule(lr){0-0}\cmidrule(lr){2-4}\cmidrule(lr){5-5},
        },
        every last row/.style={after row=\bottomrule},
        columns/0/.style ={column name={},fixed,precision=0,fixed zerofill},
        col sep=comma,
        highlightrow={1},
    ]
    {Results/Frustrum_times_80_steps_table.csv}
    \end{subtable}
    \hfill
    \begin{subtable}[t]{0.45\linewidth}
    \centering
    \footnotesize
    \caption{$\Delta L=0.0625$}
    \label{tab:frustrum_dL0.0625}
    \pgfkeys{/pgf/number format/.cd,fixed,precision=1,fixed zerofill}
    \begin{filecontents*}{Results/Frustrum_times_640_steps_table.csv}
0,499.66,2575.94,3075.6,3055.28
1,467.496,2232.51,2700.006,2783.82
2,496.28,1336.95,1833.23,1573.41
4,467.822,654.402,1122.224,789.503
8,490.113,322.555,812.668,489.392
16,467.553,167.592,635.145,495.96
32,494.06,97.0809,591.1409,483.911
64,491.42,55.6696,547.0896,493.593
128,485.046,41.4755,526.5215,494.498
256,493.78,32.8887,526.6687,491.432
512,491.835,25.7793,517.6143,488.525
\end{filecontents*}
\pgfplotstabletypeset[
        header=false,
        every head row/.style={
            output empty row,
            before row={%
                \toprule
                \textbf{Workers} & \multicolumn{3}{c}{\textbf{ASPALM}} & \textbf{APALM}\\
                    \#  & \multicolumn{1}{c@{\hspace*{\tabcolsep}\makebox[0pt]{+}}}{Serial} & \multicolumn{1}{c@{\hspace*{\tabcolsep}\makebox[0pt]{=}}}{Parallel} & Total       & Parallel \\
            },
            after row=\cmidrule(lr){0-0}\cmidrule(lr){2-4}\cmidrule(lr){5-5},
        },
        every last row/.style={after row=\bottomrule},
        columns/0/.style ={column name={},fixed,precision=0,fixed zerofill},
        col sep=comma,
        highlightrow={1},
    ]
    {Results/Frustrum_times_640_steps_table.csv}
    \end{subtable}
\end{table}

\subsection{Strip buckling}
The third example involves a benchmark problem consisting of a bifurcation instability. The problem consists of a flat strip that is clamped on one edge and free on all the others, with an in-plane compressive load applied on the free end opposite to the clamped edge; see \cref{fig:APALM_beam_case} for the problem set-up and \cite{Pagani2018} for the reference results. The ALM that is used is a Crisfield method with $\Psi=0$, with a pre-buckling arc-length of $5\cdot 10^{-5}$, a post-buckling arc-length of $5$, and a tolerance of the APALM of $\varepsilon_l = \varepsilon_u = 10^{-3}$. The serial ALM is equipped with an extension for the computation of singular points (Wriggers 1988); see \cite{Verhelst2019} for more details on this implementation. Using these methods, an initially flat strip is compressed until the bifurcation point has been computed. As soon as the strip becomes unstable, the bifurcation point is computed, and a branch switch is performed, marking the transition between the pre-buckling and post-buckling branches.\\

The results for the buckled strip are presented in \cref{fig:APALM_beam}. In this figure, the non-dimensional horizontal and vertical displacements of the end point are plotted with respect to the non-dimensional applied load. In the plots, the pre- and post-buckling branches are plotted separately for clarity, but the branches should obviously be connected at the bifurcation point. As can be seen from the results, a rather coarse serial approximation of the post-buckling branch gives a good starting point for a multi-level approximation of the curve, providing additional detail in the sharp corner in $W/L\in[0.7,0.8]$. In addition, it can be seen that the pre-buckling branch requires no more levels than the first, as the behaviour there is just a linear axial compression, hence the solution path is straight.\\

\begin{figure}[h]
\centering
\input{Figures/ShellBeam.pgf}
\caption{The problem definition for the benchmark of the buckling of a strip with length $L=1\:[\text{m}]$, width $W=0.01\:[\text{m}]$ and thickness $t=0.01\:[\text{m}]$ subject to a horizontal load with magnitude $p=0.1\:[\text{N}]$ with magnification factor $\lambda$. The strip has fixed displacements and rotations at $\Gamma_1$ and fixed displacements in $y$-direction on $\Gamma_2$ and $\Gamma_4$. The material is modelled using a Saint-Venant Kirchhoff material model with Young's modulus $E=75\cdot 10^6\:[\text{N}/\text{mm}^2]$ and Poisson ratio $\nu=0\:[-]$.}
\label{fig:APALM_beam_case}
\end{figure}

\begin{figure}[h]
\begin{subfigure}{0.45\linewidth}
 \begin{tikzpicture}
    \begin{groupplot}
        [
            group style={
                group name=my plots,
                group size=1 by 2,
                xlabels at=edge bottom,
                xticklabels at=edge bottom,
                vertical sep=5pt
            },
            xlabel=$w/L$,
            ylabel=$4PL^2 / \pi^2EI$,
            legend pos = north west,
            xmin = 0.0,
            xmax = 0.8,
            enlarge x limits=true,
            enlarge y limits=true,
            width=\linewidth,
            grid=major,
            extra y ticks={1},
        ]
        \nextgroupplot
        [
            ymin=1,
            ymax=15,
            height=0.3\textheight
        ]
        \addplot+[black,no markers] table[header=true,x expr = \thisrowno{14},y index = {13}, col sep = comma]{Results/beam_ref_1.csv};
        \addlegendentry{Reference};

        \addplot+[only marks,mark=o,mark size=2,black] table[header=true,x expr = -\thisrowno{9},y index = {8}, col sep = comma]{Results/beam_serial_1.csv};
        \addlegendentry{Serial}

        \pgfplotsset{cycle list shift=-2}

        \foreach \level in {1,...,6}
        {
            \addplot+[only marks] table[header=true,x expr = -\thisrowno{9},y index = {8}, col sep = comma,restrict expr to domain={\thisrowno{6}}{\level:\level},restrict y to domain={0.8:20}]{Results/beam_parallel_1.csv};
            \addlegendentryexpanded{Level \level}
        }

        \nextgroupplot
        [
            ymin=0.0,
            ymax=1.0,
            height=0.15\textheight
        ]
        \addplot+[black,no markers] table[header=true,x expr = \thisrowno{14},y index = {13}, col sep = comma]{Results/beam_ref_0.csv};

        \addplot+[only marks,mark=o,mark size=2,black] table[header=true,x expr = -\thisrowno{9},y index = {8}, col sep = comma]{Results/beam_serial_0.csv};
        \pgfplotsset{cycle list shift=-2}
        \foreach \level in {1,...,6}
        {
            \addplot+[only marks] table[header=true,x expr = -\thisrowno{9},y index = {8}, col sep = comma,restrict expr to domain={\thisrowno{6}}{\level:\level},restrict y to domain={0:1}]{Results/beam_parallel_0.csv};
        }

      \end{groupplot}
 \end{tikzpicture}
 \caption{Non-dimensional out-of-plane displacement of the beam with respect to the non-dimesional buckling load.}
 \label{fig:APALM_beam_vert}
\end{subfigure}
\hfill
\begin{subfigure}{0.45\linewidth}
 \begin{tikzpicture}
        \begin{groupplot}
        [
            group style={
                group name=my plots,
                group size=1 by 2,
                xlabels at=edge bottom,
                xticklabels at=edge bottom,
                vertical sep=5pt
            },
            xlabel=$-u/L$,
            ylabel=$4PL^2 / \pi^2EI$,
            legend pos = north west,
            xmin = 0.0,
            xmax = 2.0,
            width=\linewidth,
            grid=major,
            extra y ticks={1},
        ]
        \nextgroupplot
        [
            ymin=1.0,
            ymax=15,
            height=0.3\textheight
        ]
        \addplot+[black,no markers] table[header=true,x expr = \thisrowno{15},y index = {13}, col sep = comma]{Results/beam_ref_1.csv};
        \addlegendentry{Reference};

        \addplot+[only marks,mark=o,mark size=2,black] table[header=true,x expr = \thisrowno{10},y index = {8}, col sep = comma]{Results/beam_serial_1.csv};
        \addlegendentry{Serial};

        \pgfplotsset{cycle list shift=-2}
        \foreach \level in {1,...,6}
        {
            \addplot+[only marks] table[header=true,x expr = \thisrowno{10},y index = {8}, col sep = comma,restrict expr to domain={\thisrowno{6}}{\level:\level},restrict y to domain={1:20}]{Results/beam_parallel_1.csv};
            \addlegendentryexpanded{Level \level}
        }

        \nextgroupplot
        [
            ymin=0.0,
            ymax=1.0,
            height=0.15\textheight
        ]
        \addplot+[black,no markers] table[header=true,x expr = \thisrowno{15},y index = {13}, col sep = comma]{Results/beam_ref_0.csv};

        \addplot+[only marks,mark=o,mark size=2,black] table[header=true,x expr = \thisrowno{10},y index = {8}, col sep = comma]{Results/beam_serial_0.csv};
        \pgfplotsset{cycle list shift=-2}
        \foreach \level in {1,...,6}
        {
            \addplot+[only marks] table[header=true,x expr = \thisrowno{10},y index = {8}, col sep = comma,restrict expr to domain={\thisrowno{6}}{\level:\level},restrict y to domain={0:1}]{Results/beam_parallel_0.csv};
        }
      \end{groupplot}
 \end{tikzpicture}
  \caption{Non-dimensional in-of-plane displacement (length direction) of the beam with respect to the non-dimensional buckling load.}
 \label{fig:APALM_beam_hor}
\end{subfigure}
\caption{Results of the buckling of a clamped strip. The left figure provides the out-of-plane displacement of the free end with respect to the non-dimensional load $4PL^2 / \pi^2EI$, and the right figure represents the horizontal displacement of the free end with respect to the same non-dimensional load. In both figures, buckling occurs when $4PL^2 / \pi^2EI=1$ and the axes are split on this point to make the pre- and post-buckling branches both visible. The simulation is performed with a tolerance of $\varepsilon_l=\varepsilon_u=10^{-3}$ and a step length of $\Delta L=5\cdot 10^{-5}$ (pre-buckling) and $\Delta L=5$ (post-buckling).}
 \label{fig:APALM_beam}
\end{figure}

As for the previous two benchmark examples, a scaling analysis of the parallel evaluations is performed. The main difference between the previous two examples is that the present example involves a bifurcation point. However, since the job queue includes the jobs from all branches together, there is no idle time to wait for a branch to finish before starting a new branch; hence, it is expected that the parallel scaling for a bifurcation problem should have the same scaling properties. Indeed, \cref{tab:beam_dL2.5} shows that optimal scaling is achieved in the parallel correction phase of the ASPALM up to 8 nodes, after which the idle time to wait for the last job to finish significantly impacts the scaling, as observed in the other benchmarks. In addition, it is found that the APALM reaches efficient computation of the full adaptive load-displacement curve within the time of a serial ALM computation, using 8 workers. When increasing the number of curve segments by decreasing the arc-length parameter to $\Delta L=0.025$, it can again be seen that the parallel scalability increases. Optimal scaling of the parallel correction of the ASPALM is achieved up to 256 workers, with the correction phase taking only 1\% of the total computational time. The APALM provides an adaptively refined solution curve in the time of a serial ALM computation, using only 4 to 8 workers.\\

\begin{table}[h]
    \caption{Computational time in $[\text{s}]$ for benchmark of the buckled strip using the ASPALM and APALM with different numbers of worker processes. The times for the ASPALM are presented for the serial initialization and the parallel correction phases, and the sum of the two is given as the total computational time. The numbers in the \textit{Serial} column should theoretically be the same, but they provide a representation of the variation in the time measurements. The results are presented for simulations with increment lengths $\Delta L=2.5$ (\subref{tab:beam_dL2.5}) and $\Delta L=0.025$ (\subref{tab:beam_dL0.025}), and the italic row with $0$ workers denotes the ASALM method.}
    \label{tab:beam}
    \begin{subtable}[t]{0.45\linewidth}
    \centering
    \footnotesize
    \caption{$\Delta L=2.5$}
    \label{tab:beam_dL2.5}
    \pgfkeys{/pgf/number format/.cd,fixed,precision=1,fixed zerofill}
    \begin{filecontents*}{Results/Beam_times_5_steps_table.csv}
0,50.7268,171.759,222,244
1,58.6929,198.133,257,244
2,58.462,102.578,161,112
4,58.9623,51.0642,110,77
8,58.9625,27.9202,87,67
16,59.3445,24.9421,84,68
32,60.2796,23.046,83,68
64,58.6245,24.3947,83,66
\end{filecontents*}
\pgfplotstabletypeset[
        header=false,
        every head row/.style={
            output empty row,
            before row={%
                \toprule
                \textbf{Workers} & \multicolumn{3}{c}{\textbf{ASPALM}} & \textbf{APALM}\\
                    \#  & \multicolumn{1}{c@{\hspace*{\tabcolsep}\makebox[0pt]{+}}}{Serial} & \multicolumn{1}{c@{\hspace*{\tabcolsep}\makebox[0pt]{=}}}{Parallel} & Total       & Parallel \\
            },
                        after row=\cmidrule(lr){0-0}\cmidrule(lr){2-4}\cmidrule(lr){5-5},
        },
        every last row/.style={after row=\bottomrule},
        columns/0/.style ={column name={},fixed,precision=0,fixed zerofill},
        col sep=comma,
        highlightrow={1},
    ]
    {Results/Beam_times_5_steps_table.csv}
    \end{subtable}
    \hfill
    \begin{subtable}[t]{0.45\linewidth}
    \centering
    \footnotesize
    \caption{$\Delta L=0.025$}
    \label{tab:beam_dL0.025}
    \pgfkeys{/pgf/number format/.cd,fixed,precision=1,fixed zerofill}
    \begin{filecontents*}{Results/Beam_times_500_steps_table.csv}
0,963.348,1963.26,2926.608,3017.55
1,1022.3,2065.66,3087.96,3020.56
2,1025,1053.73,2078.73,1509.7
4,943.939,463.985,1407.924,1034.1
8,1019.5,256.035,1275.535,1028.19
16,1006.21,129.267,1135.477,1027.72
32,1026.61,68.8346,1095.4446,1028.3
64,1032.72,33.6395,1066.3595,1027.96
128,935.134,18.853,953.987,1023.51
256,1012.59,11.9872,1024.5772,1026.83
\end{filecontents*}
\pgfplotstabletypeset[
        header=false,
        every head row/.style={
            output empty row,
            before row={%
                \toprule
                \textbf{Workers} & \multicolumn{3}{c}{\textbf{ASPALM}} & \textbf{APALM}\\
                    \#  & \multicolumn{1}{c@{\hspace*{\tabcolsep}\makebox[0pt]{+}}}{Serial} & \multicolumn{1}{c@{\hspace*{\tabcolsep}\makebox[0pt]{=}}}{Parallel} & Total       & Parallel \\
            },
                        after row=\cmidrule(lr){0-0}\cmidrule(lr){2-4}\cmidrule(lr){5-5},
        },
        every last row/.style={after row=\bottomrule},
        columns/0/.style ={column name={},fixed,precision=0,fixed zerofill},
        col sep=comma,
        highlightrow={1},
    ]
    {Results/Beam_times_500_steps_table.csv}
    \end{subtable}
\end{table}

\subsection{Snapping Meta-Material}\label{sec:APALM_Results_Snapping}
As a final example, the APALM is applied to a problem of larger scale. In particular, the snapping behaviour of a snapping meta-material is modelled, inspired by \cite{Rafsanjani2015}. The meta-material consists of $N_x\times N_Y=3\times2.5$ building blocks (see \cref{fig:APALM_snapping_meta-material}) with a snapping and a bearing segment (see \cref{fig:APALM_snapping_element}), and the material is modelled with a compressible Neo-Hookean material model. The full problem details are provided in \cref{fig:APALM_snapping_case}. The snapping behaviour of the meta-material is investigated by using arc-length methods on the varying load $\lambda P$, with a step size of $\Delta L=5\cdot 10^{-2}$, until $1.5\%$ strain. The simulation is modelled using 2D elasticity equations using the plane-stress assumption, which are discretised using B-splines with mixed degrees 2 and 3 and maximal regularity. The mesh is uniformly distributed, and the full system of equations has $16563$ degrees of freedom. The simulations are performed using shared-memory parallelization of the assembly routines and distributed memory parallelization of the ASPALM and APALM. For reference, a displacement-controlled (DC) simulation is performed.\\

\begin{figure}
    \centering
    \begin{subfigure}{0.45\linewidth}
    \centering
    \resizebox{\linewidth}{!}{\input{Figures/SnappingGlobal.pgf}}
    \caption{A snapping meta-material with $3\times2.5$ building blocks, of which one is outlined. The total multi-patch consists of 132 patches.}
    \label{fig:APALM_snapping_meta-material}
    \end{subfigure}
    \hfill
    \begin{subfigure}{0.45\linewidth}
    \centering
    \resizebox{\linewidth}{!}{\input{Figures/SnappingLocal.pgf}}
    \caption{The snapping building block, composed of 15 patches outlines in black.}
    \label{fig:APALM_snapping_element}
    \end{subfigure}
    \caption{The problem definition for the snapping meta material using a grid of $3\times2.5$ elements (\subref{fig:APALM_snapping_meta-material}) with the element geometry as defined in (\subref{fig:APALM_snapping_element}). The element dimensions are defined using the thickness of the load-bearing part $t_b=1.5\:[\text{mm}]$ and the thickness of the snapping part $t_s=1.0\:[\text{mm}]$, the thickness of the gap $t_g=1.0\:[\text{mm}]$ and the thickness of the connectors $t_w=1.5\:[\text{mm}]$, such that the height $h=t_b+t_s+2t_g$. The length of the element is $\ell=10\:[\text{mm}]$, and the amplitude of the cosine wave defining the element shape is given by $a=0.3l$. Since the meta-material has $3\times2.5$ elements, the total width is $W=3\ell$. The height of the total metamaterial is given by $H=3h+2t_g+t_s+h_B+h_T$, where $h_B=h_T=5t_g$ are the buffer zones on the top and the bottom. The thickness of the specimen (in out-of-plane direction) is $b=3\:[\text{mm}]$. The material is defined using a compressible Neo-Hookean material model with Young's modulus $E=78\:[\text{N}/\text{mm}^2]$ and Poisson ratio $\nu=0.4\:[-]$. The bottom boundary $\Gamma_1$ is fixed using $u_x=u_y=0$, and the top boundary $\Gamma_2$ is fixed in the horizontal direction ($u_x=0$) and coupled in the vertical direction $u_y$. The load applied on the top boundary is a variable defined by $\lambda P$.}
    \label{fig:APALM_snapping_case}
\end{figure}

\Cref{fig:APALM_snapping} depicts the stress-strain curve for the snapping meta-material depicted in \cref{fig:APALM_snapping_case}. Firstly, it can be seen that the curve computed using a serial ALM coincides with the curves obtained from DC simulations, but the ALM shows additional snap-through behaviour on the points where the DC curve has kinks. These kinks coincide with the instability in the metamaterial. Furthermore, the figure shows that the initial course approximation at level 0 is refined up to level 4 in the adaptive arc-length method scheme proposed in this paper, provided the tolerance of $\varepsilon_l=\varepsilon_u=10^{-3}$. The refinements of the adaptive scheme are mainly present in the highly curved segments of the load displacement curve. The reader is referred to Video 1 from the supplementary material for the deformations corresponding to the stress-strain curve in \cref{fig:APALM_snapping}.\\

As for the other numerical experiments, the parallelization properties of the ASPALM and APALM schemes are investigated. For the snapping meta-material simulation, the computational times are presented in \cref{tab:snapping}. The results in the table show high computational times for the ASALM (i.e., the ASPALM and APALM with 0 workers) compared to the DC simulation. However, the scalability observed in the previous benchmark problems can also be observed in the simulation of the snapping metamaterial. In fact, the APALM with 8 workers requires a factor of 4 less computational time, again equivalent to the computational time required for only the serial initialization phase of the ASPALM. Lastly, the case of the snapping meta-material shows that, compared to a naturally serial displacement-controlled method, the APALM achieves a speed-up of a factor of 2.5 while providing snapping behaviour with greater accuracy.\\

\begin{figure}
\def\minXi{0.5}
\def\maxXi{1.2}
\def\minYi{0.12}
\def\maxYi{0.23}
\centering
\begin{subfigure}[t]{\linewidth}
    \centering
    \begin{tikzpicture}
        \begin{axis}
        [
            xlabel={Strain $\varepsilon\:[-]$},
            ylabel={Equivalent stress $\sigma\:[\text{MPa}]$},
            legend pos = south east,
            xmin = 0,
            xmax = 1.5,
            ymin = 0,
            ymax = 0.25,
            enlarge x limits = true,
            enlarge y limits = true,
            width=0.5\linewidth,
            height=0.3\textheight,
            grid=major,
            legend columns=2,
            restrict x to domain = 0:1.5
        ]
        \addplot+[dashed,no markers,black,thick] table[header=true,x expr = \thisrowno{3},y expr = -\thisrowno{4}/1e6, col sep = comma]{Results/Snapping_DC.csv};
        \addlegendentry{DC}
        \def\level{0}
        \addplot+[only marks,mark=o,mark size=2.0,black] table[header=true,x expr =\thisrowno{3},y expr = \thisrowno{4}/1e6, col sep = comma,restrict expr to domain={\thisrowno{7}}{\level:\level}]{Results/Snapping_parallel.csv};
        \addlegendentry{Serial}
        \addplot+[no markers,gray,thin,solid,] table[header=true,x expr = \thisrowno{3},y expr = \thisrowno{4}/1e6, col sep = comma]{Results/Snapping_parallel.csv};
        \addlegendentry{Parallel}

        \def\level{1}
        \pgfplotsset{cycle list shift=-3}
        \addplot+[only marks,mark size=1.0] table[header=true,x expr =\thisrowno{3},y expr = \thisrowno{4}/1e6, col sep = comma,restrict expr to domain={\thisrowno{7}}{\level:\level}]{Results/Snapping_parallel.csv};
        \addlegendentry{Level 1}
        \draw[gray,dashed] (axis cs: \minXi,\minYi) -- (axis cs: \maxXi,\minYi) -- (axis cs: \maxXi,\maxYi) -- (axis cs: \minXi,\maxYi) node[above right, inner sep=1pt] {A} --cycle;
     \end{axis}
     \end{tikzpicture}
     \begin{tikzpicture}
         \begin{groupplot}
        [
            group style={
                group size=2 by 2,
                xlabels at=edge bottom,
                xticklabels at=edge bottom,
                ylabels at=edge right,
                yticklabels at=edge right,
                vertical sep=0.025\textheight,
                horizontal sep=0.025\linewidth,
            },
            xlabel style={color=gray},
            x tick label style={color = gray},
            ylabel style={color=gray},
            y tick label style={color = gray},
            legend pos = south east,
            width=0.29\linewidth,
            height=0.172\textheight,
            grid=major,
            legend columns=2,
            restrict x to domain = 0:1.5,
            ylabel={\phantom{.}},
            xlabel={Strain $\varepsilon\:[-]$}
        ]
        \nextgroupplot[
            xmin = \minXi,
            xmax = \maxXi,
            ymin = \minYi,
            ymax = \maxYi,
            draw = gray,
        ]
        \addplot+[dashed,no markers,black,thick,forget plot] table[header=true,x expr = \thisrowno{3},y expr = -\thisrowno{4}/1e6, col sep = comma]{Results/Snapping_DC.csv};
        \addplot+[no markers,gray,thin,solid,forget plot] table[header=true,x expr = \thisrowno{3},y expr = \thisrowno{4}/1e6, col sep = comma]{Results/Snapping_parallel.csv};
        \pgfplotsset{cycle list shift=-1}
        \def\level{2}
        \addplot+[only marks,mark size=0.5] table[header=true,x expr =\thisrowno{3},y expr = \thisrowno{4}/1e6, col sep = comma,restrict expr to domain={\thisrowno{7}}{\level:\level}]{Results/Snapping_parallel.csv};
        \addlegendentry{Level \level}
        \nextgroupplot[
            xmin = \minXi,
            xmax = \maxXi,
            ymin = \minYi,
            ymax = \maxYi,
            draw = gray,
        ]
        \addplot+[dashed,no markers,black,thick,forget plot] table[header=true,x expr = \thisrowno{3},y expr = -\thisrowno{4}/1e6, col sep = comma]{Results/Snapping_DC.csv};
        \addplot+[no markers,gray,thin,solid,forget plot] table[header=true,x expr = \thisrowno{3},y expr = \thisrowno{4}/1e6, col sep = comma]{Results/Snapping_parallel.csv};
        \pgfplotsset{cycle list shift=2}
        \def\level{3}
        \addplot+[only marks,mark size=0.5] table[header=true,x expr =\thisrowno{3},y expr = \thisrowno{4}/1e6, col sep = comma,restrict expr to domain={\thisrowno{7}}{\level:\level}]{Results/Snapping_parallel.csv};
        \addlegendentry{Level \level}
        
        \nextgroupplot[
            xmin = \minXi,
            xmax = \maxXi,
            ymin = \minYi,
            ymax = \maxYi,
            draw = gray,
        ]
        \addplot+[dashed,no markers,black,thick,forget plot] table[header=true,x expr = \thisrowno{3},y expr = -\thisrowno{4}/1e6, col sep = comma]{Results/Snapping_DC.csv};
        \addplot+[no markers,gray,thin,solid,forget plot] table[header=true,x expr = \thisrowno{3},y expr = \thisrowno{4}/1e6, col sep = comma]{Results/Snapping_parallel.csv};
        \pgfplotsset{cycle list shift=3}
        \def\level{4}
        \addplot+[only marks,mark size=0.5] table[header=true,x expr =\thisrowno{3},y expr = \thisrowno{4}/1e6, col sep = comma,restrict expr to domain={\thisrowno{7}}{\level:\level}]{Results/Snapping_parallel.csv};
        \addlegendentry{Level \level}

        \nextgroupplot[
            xmin = \minXi,
            xmax = \maxXi,
            ymin = \minYi,
            ymax = \maxYi,
            draw = gray,
        ]
        \addplot+[dashed,no markers,black,thick,forget plot] table[header=true,x expr = \thisrowno{3},y expr = -\thisrowno{4}/1e6, col sep = comma]{Results/Snapping_DC.csv};
        \addplot+[no markers,gray,thin,solid,forget plot,
        postaction={decorate, decoration={markings,
        mark=at position 0.2 with {\arrow{latex};},
        mark=at position 0.3 with {\arrow{latex};},
        mark=at position 0.4 with {\arrow{latex};},
        mark=at position 0.5 with {\arrow{latex};},
        mark=at position 0.6 with {\arrow{latex};},
        mark=at position 0.7 with {\arrow{latex};},
        mark=at position 0.8 with {\arrow{latex};},
        mark=at position 0.9 with {\arrow{latex};},
        }}] table[header=true,x expr = \thisrowno{3},y expr = \thisrowno{4}/1e6, col sep = comma ]{Results/Snapping_parallel.csv};
        \def\level{5}
        \pgfplotsset{cycle list shift=4}
        \addplot+[only marks,mark size=0.5] table[header=true,x expr =\thisrowno{3},y expr = \thisrowno{4}/1e6, col sep = comma,restrict expr to domain={\thisrowno{7}}{\level:\level}]{Results/Snapping_parallel.csv};
        \addlegendentry{Level \level}
      \end{groupplot}
 \end{tikzpicture}
\end{subfigure}
\caption{Stress-strain diagram for the snapping meta-material from \cref{fig:APALM_snapping_case}. The vertical axis depicts the equivalent stress $\sigma=\lambda P / (b W)$, and the horizontal axis represents the strain $\varepsilon=u_y/H$, where $u_y$ is the displacement of the top boundary $\Gamma_2$. The complete curve with the displacement-controlled (DC) results, the points obtained in serial initialization, and the line obtained by parallel corrections are presented on the left. The figures on the right present the points from different hierarchical levels at the inset depicted in the left diagram. The simulation is performed with a tolerance of $\varepsilon_l=\varepsilon_u=10^{-3}$ and an increment length of $\Delta L=0.05$.}
\label{fig:APALM_snapping}
\end{figure}

\begin{table}
	\footnotesize
    \centering
    \caption{Computational time in $[\text{s}]$ for the example of the snapping meta-material for the ASPALM and APALM for different numbers of worker processes. The computational time for a displacement-controlled (DC) simulation with step $\Delta u_y=0.0005\:[mm]$ is provided as a reference. The times for the ASPALM are presented for the serial initialization and the parallel correction phases, and the sum of the two is given as the total computational time. The numbers in the \textit{Serial} column should theoretically be the same, but they provide a representation of the variation in the time measurements. The italic row with $0$ workers denotes the ASALM method.}
    \label{tab:snapping}
    \pgfkeys{/pgf/number format/.cd,fixed,precision=1,fixed zerofill}
    \begin{filecontents*}{Results/snapping_times_table.csv}
0,1571.61,5204.75,6776.36,7022.94,4400.75
1,1686.88,4593.23,6280.11,5319.06,
2,1237.49,3005.89,4243.38,3827.88,
4,1742.71,1548.22,3290.93,2137.34,
8,1445.4,717.406,2162.806,1711.84,
16,1931.11,352.172,2283.282,1632.92,
32,1746.89,219.745,1966.635,1755.55,
\end{filecontents*}
\pgfplotstabletypeset[
        header=false,
        every head row/.style={
            output empty row,
            before row={%
                \toprule
                \textbf{Workers} & \multicolumn{3}{c}{\textbf{ASPALM}} & \textbf{APALM} & \textbf{DC}\\
                    \#  & \multicolumn{1}{c@{\hspace*{\tabcolsep}\makebox[0pt]{+}}}{Serial} & \multicolumn{1}{c@{\hspace*{\tabcolsep}\makebox[0pt]{=}}}{Parallel} & Total       & Parallel & Serial\\
            },
                        after row=\cmidrule(lr){0-0}\cmidrule(lr){2-4}\cmidrule(lr){5-5}\cmidrule(lr){6-6},
        },
        every last row/.style={after row=\bottomrule},
        columns/0/.style ={column name={},fixed,precision=0,fixed zerofill},
        col sep=comma,
        highlightrow={1},
    ]
    {Results/snapping_times_table.csv}
\end{table}

%% file: Figures/Roof.pgf
\tdplotsetmaincoords{-20}{0}

\begin{tikzpicture}[scale=2,tdplot_main_coords]
\tdplotsetrotatedcoords{0}{-20}{0}
\def\angle{30}
\def\radius{2.54}
\def\length{5.08}
\begin{scope}[tdplot_rotated_coords]
\coordinate (start1) at (0,0,0);

\draw[ultra thin] (start1) arc [start angle=-\angle+90,end angle=\angle+90,x radius=\radius,y radius=\radius] coordinate (end1) coordinate[midway](mid1);

\coordinate (start2) at (0,0,\length);
\draw[ultra thin] (start2) arc [start angle=-\angle+90,end angle=\angle+90,x radius=\radius,y radius=\radius] coordinate (end2) coordinate[midway](mid2);

\draw(mid1)--(mid2) coordinate[midway](mid);

\filldraw[right color=black!50,left color=black!10] (start1) arc [start angle=-\angle+90,end angle=\angle+90,x radius=\radius,y radius=\radius] node[above right,midway] {$\Gamma_4$}-- (end2) node[above left,midway]{$\Gamma_3$} arc [start angle=\angle+90,end angle=-\angle+90,x radius=\radius,y radius=\radius] node[above right,midway] {$\Gamma_2$}--(start1) node[above left,midway]{$\Gamma_1$};

\begin{scope}[canvas is xz plane at y=0.5]
\filldraw[fill=white,draw=black,rotate=0] (mid) ellipse (0.04 and 0.04);
\end{scope}
\draw[latex-,thick] (mid) node[below]{$A$} --($ (mid)+(0,0.5,0)$) node[above]{$\lambda P$};

\draw [latex-latex] ($(start2)+(0.2,0,0)$) -- ($(start1)+(0.2,0,0)$) node[below right,midway] {$L$};
\draw [latex-] ($(mid2)+(0,0,0)$) -- ($(mid2)+(0,-2,0)$) node[right,midway] {$R$} coordinate(bot) coordinate[midway] (midtheta2);
\draw [dashed] (bot)--(end2) coordinate[midway] (midtheta1);
\draw (midtheta1) to[in=-200,out=40](midtheta2);
\node [above left,inner sep=12] at (bot) {$\theta$};
\end{scope}
\end{tikzpicture}

%% file: Figures/Frustrum.pgf
\begin{tikzpicture}
\def\R{4}
\def\r{2}
    \begin{axis}[
        hide axis,
        width=\linewidth,
        height=\linewidth,
        axis equal image,
        view = {100}{20},
        xmin = -\R-0.5,
        xmax = \R+0.5,
        ymin = -\R-0.5,
        ymax = \R+0.5,
        zmin = 0.0,
        zmax = 3.0,
    ]
        \addplot3[
            surf,
            samples = 50,
            samples y = 2,
            domain = 0.5*pi:2*pi,
            domain y = 0:1,
z buffer = sort,
no marks,
mesh/interior colormap={blueblack}{color=(black!20) color=(white)},
colormap ={blueblack}{color=(black!20) color=(white)},
shader=interp,
point meta=0.5*y*y+2*x*x-0.4*z,
        ](
            {((\y*(\r-\R)+\R)*cos(deg(\x)},
            {(\y*(\r-\R)+\R)*sin(deg(\x))},
            {2*\y}
        );

        \addplot3[
                    surf,
                    samples = 15,
            samples y = 2,
            domain = 0:0.5*pi,
            domain y = 0:1,
            draw=none,
                    z buffer = sort,
                    no marks,
                    mesh/interior colormap={blueblack}{color=(black!50) color=(white)},
                    colormap ={blueblack}{color=(black!50) color=(white)},
                    shader=interp,
					point meta=0.5*y*y+1*x*x-0.4*z,
        ](
            {((\y*(\r-\R)+\R)*cos(deg(\x)},
            {(\y*(\r-\R)+\R)*sin(deg(\x))},
            {2*\y}
        );

\addplot3[domain=-0.5*pi:0.5*pi, samples=100, samples y=0, no marks, smooth, thick,black](
            {\R*cos(deg(\x)},
            {\R*sin(deg(\x))},
            {0}
    );
\addplot3[domain=-0.5*pi:0.5*pi, samples=100, samples y=0, no marks, smooth, thick,gray](
            {-\R*cos(deg(\x)},
            {-\R*sin(deg(\x))},
            {0}
    );

\addplot3[domain=0:2*pi, samples=100, samples y=0, no marks, smooth, thick,black](
            {\r*cos(deg(\x)},
            {\r*sin(deg(\x))},
            {2}
    );

\foreach \theta in {0.50,0.6,...,2} {

	\edef\x{\r*cos(deg(\theta*3.1415} 
	\edef\y{\r*sin(deg(\theta*3.1415} 
    \edef
	\temp{
				\noexpand
				\draw [gray,latex-] (axis cs:\x,\y,2) -- (axis cs:\x,\y,3);
			}
    \temp
}
\node[above left] at (axis cs: 0,-\r,3) {$\lambda p$};

\addplot3[domain=0.5*pi:2*pi, samples=100, samples y=0, no marks, smooth, thick,gray](
            {\r*cos(deg(\x)},
            {\r*sin(deg(\x))},
            {3}
    );

\foreach \theta in {0,0.1,...,0.50} {

    \edef\x{\r*cos(deg(\theta*3.1415} 
    \edef\y{\r*sin(deg(\theta*3.1415} 
    \edef
    \temp{
                \noexpand
                \draw [latex-] (axis cs:\x,\y,2) -- (axis cs:\x,\y,3);
            }
    \temp
}

\addplot3[domain=0:0.5*pi, samples=100, samples y=0, no marks, smooth, thick,black](
            {\r*cos(deg(\x)},
            {\r*sin(deg(\x))},
            {3}
    );

\draw[thick] (axis cs: 0,\R,0) -- (axis cs: 0,\r,2) node[midway,below,]{$\Gamma_1$}; 
\draw[thick] (axis cs: \R,0,0) -- (axis cs: \r,0,2) node[midway,below right,]{$\Gamma_3$}; 


\node[below] at (axis cs: 0.707*\r,0.707*\r,2) {$\Gamma_2$};
\node[above] at (axis cs: 0.707*\R,0.707*\R,0) {$\Gamma_4$};

\draw[thick] (axis cs: 0,\R,0) -- (axis cs: 0,\r,2);
\draw[thick] (axis cs: \R,0,0) -- (axis cs: \r,0,2);

\draw[latex-latex] (axis cs: 0,0,2) -- (axis cs: -0.707*\r,-0.707*\r,2) node[midway,below]{$R_1$};
\draw[latex-latex] (axis cs: 0,0,0) -- (axis cs: -0.707*\R,-0.707*\R,0) node[midway,below]{$R_2$};

\draw[latex-latex] (axis cs:-0.707*\R-0.5,0.707*\R+0.5,0) -- (axis cs: -0.707*\R-0.5,0.707*\R+0.5,2) node[midway,left]{$H$};
\draw[dashed] (axis cs:-0.707*\R-0.5,0.707*\R+0.5,0) -- (axis cs: 0,0,0);
\draw[dashed] (axis cs:-0.707*\R-0.5,0.707*\R+0.5,2) -- (axis cs: 0,0,2);
\node at (axis cs: 0,0,2) {$\times$};
\node at (axis cs: 0,0,0) {$\times$};

\draw[-latex] (axis cs: -\R+0.5,-\R-0.5,0) -- (axis cs: -\R+0.5,-\R-0.5,1.0) node[left]{$z$};
\draw[-latex] (axis cs: -\R+0.5,-\R-0.5,0) -- (axis cs: -\R+0.5,-\R-0.5+1.0,0) node[above]{$y$};
\draw[-latex] (axis cs: -\R+0.5,-\R-0.5,0) -- (axis cs: -\R+0.5+1.0,-\R-0.5,0) node[below]{$x$};

\end{axis}
\end{tikzpicture}

%% file: Figures/ShellBeam.pgf
\tdplotsetmaincoords{0}{0}
\tdplotsetrotatedcoords{-30}{40}{40}
\begin{tikzpicture}[tdplot_rotated_coords]
\draw[top color=black!60,bottom color=black!20](-0.5,0,0) -- 
node[midway,below left]{$\Gamma_4$} (-0.5,0,5) -- 
node[midway,above left]{$\Gamma_3$} (0.5,0,5)  -- 
node[midway,above right]{$\Gamma_2$}(0.5,0,0)  --
node[midway,above left]{$\Gamma_1$}
cycle;


\draw[latex-] (0,0,5) -- (0,0,6) node[ right]{$\lambda p$};

\foreach \x in {-0.5,-0.4,-0.3,...,0.5}
{
	 \draw[] (\x,0,0) -- (\x-0.1,0,-0.3);
}

\draw[latex-latex] (-2,0,0) -- (-2,0,5) node[midway,below]{L};
\draw[latex-latex] (-0.5,0,7) -- (0.5,0,7) node[midway,below right]{W};

\draw[-latex] (-1.5,0,0) -- (-1,0,0) node[below right] {$y$};
\draw[-latex] (-1.5,0,0) -- (-1.5,0.5,0) node[above] {$z$};
\draw[-latex] (-1.5,0,0) -- (-1.5,0.0,0.5) node[below] {$x$};
\end{tikzpicture}

%% file: Figures/SnappingGlobal.pgf
	\def\Nx{3}
	\def\Ny{2}
	\def\al{0.2}
	\def\tw{0.15}
	\def\tb{0.15}
	\def\tg{0.1}
	\def\ts{0.1}
	\def\l{1}
	\def\samples{50}
	\begin{tikzpicture}[scale=2]
	\pgfmathsetmacro\L{\Nx*\l}
	\pgfmathsetmacro\Nym{\Ny-1}
	\pgfmathsetmacro\Nxx{2*\Nx-1}
	\foreach \y in {0,...,\Ny}
	{
		\pgfmathsetmacro\hm{\y*(\ts+2*\tg+\tb)+\al*\l/2-\tg/2}
		\pgfmathsetmacro\hp{\y*(\ts+2*\tg+\tb)+\al*\l/2+1.5*\tg}
		\foreach \x in {1,...,\Nxx}
		{
			\pgfmathsetmacro\lrm{\x*\l/2-\tw/2}
			\pgfmathsetmacro\lrp{\x*\l/2+\tw/2}
			\filldraw[color=gray2] (\lrm,\hm) rectangle (\lrp,\hp);
		}
		\pgfmathsetmacro\lrm{0}
		\pgfmathsetmacro\lrp{\tw/2}
		\filldraw[color=gray2] (\lrm,\hm) rectangle (\lrp,\hp);
		\pgfmathsetmacro\lrm{\L-\tw/2}
		\pgfmathsetmacro\lrp{\L}
		\filldraw[color=gray2] (\lrm,\hm) rectangle (\lrp,\hp);
	}

	\filldraw[color=gray1,domain=0:\L,samples=\samples] plot (\x,{-\al*\l/2*cos(\x/\l*360)})--(\L,-1) -- (0,-1) --cycle;
	\draw[black,thick] (\L,-1) -- (0,-1) node[below,midway]{$\Gamma_1$};
	\draw[black,thick,latex-latex] ($(\L,-1)+(0,0.1)$) -- ($(0,-1)+(0,0.1)$) node[above,midway]{$W$};
	\draw[black,thick,latex-latex] ($(0,-1)-(0.1,0)$) -- ($(0,0)-(0.1,0)$) node[left,midway]{$h_B$};

	\foreach \k in {0,...,\Nym}
	{
		\pgfmathsetmacro\h{\k*(\ts+\tg+\tb)+(\k+1)*\tg}
		\filldraw[color=col2,domain=0:\L,samples=\samples] plot (\x,{-\al*\l/2*cos(\x/\l*360)+\h}) -- plot[domain=\L:0] (\x,{-\al*\l/2*cos(\x/\l*360)+\h+\ts}) --cycle;
		\pgfmathsetmacro\h{\h+\ts+\tg}
		\filldraw[color=col1,domain=0:\L,samples=\samples] plot (\x,{-\al*\l/2*cos(\x/\l*360)+\h}) -- plot[domain=\L:0] (\x,{-\al*\l/2*cos(\x/\l*360)+\h+\tb}) --cycle;
	}

	\pgfmathsetmacro\h{\Ny*(\ts+\tg+\tb)+(\Ny+1)*\tg}
	\filldraw[color=col2,domain=0:\L,samples=\samples] plot (\x,{-\al*\l/2*cos(\x/\l*360)+\h}) -- plot[domain=\L:0] (\x,{-\al*\l/2*cos(\x/\l*360)+\h+\ts}) --cycle;
	\pgfmathsetmacro\H{(\Ny+1)*(2*\tg+\ts+\tb)-\tb}
	\filldraw[color=gray1,domain=0:\L,samples=\samples] plot (\x,{-\al*\l/2*cos(\x/\l*360)+\H})--(\L,\H+1) -- (0,\H+1) --cycle;
	\draw[black,thick] (\L,\H+1) -- (0,\H+1) node[below,midway]{$\Gamma_2$};
	\draw[black,thick,-latex] (\L/2,\H+1) -- (\L/2,\H+1+0.5) node[above]{$\lambda P$};
	\draw[black,thick,latex-latex]($(0,\H)-(0.1,0)$) -- ($(0,\H+1)-(0.1,0)$) node[left,midway]{$h_T$};
	\draw[black,thick,latex-latex] ($(\L,-1)+(0.1,0)$) --  ($(\L,\H+1)+(0.1,0)$)  node[right,midway]{$H$};

	\pgfmathsetmacro\h{\tg}
	\pgfmathsetmacro\hm{\al*\l/2+\tg/2}
	\pgfmathsetmacro\hmg{\al*\l/2+\ts+2*\tg+\tb+\tg/2}
	\draw[color=black,thick] (\l/2,\hm) -- (\l/2+\tw/2,\hm) -- plot[domain=\l/2+\tw/2:3*\l/2-\tw/2,samples=\samples] (\x,{-\al*\l/2*cos(\x/\l*360)+\h}) -- (3*\l/2-\tw/2,\hm) -- (3*\l/2,\hm) -- plot[domain=3*\l/2:\l+\tw/2] (\x,{-\al*\l/2*cos(\x/\l*360)+\h+\ts}) -- plot[domain=\l+\tw/2:3*\l/2] (\x,{-\al*\l/2*cos(\x/\l*360)+\h+\ts+\tg}) -- (3*\l/2,\hmg) -- (3*\l/2-\tw/2,\hmg) -- plot[domain=3*\l/2-\tw/2:\l/2+\tw/2] (\x,{-\al*\l/2*cos(\x/\l*360)+\h+\ts+\tg+\tb}) -- (\l/2+\tw/2,\hmg) -- (\l/2,\hmg)  -- plot[domain=\l/2:\l-\tw/2] (\x,{-\al*\l/2*cos(\x/\l*360)+\h+\ts+\tg})  -- plot[domain=\l-\tw/2:\l/2] (\x,{-\al*\l/2*cos(\x/\l*360)+\h+\ts}) --cycle;

	\end{tikzpicture}

%% file: Figures/SnappingLocal.pgf
\def\Nx{3}
\def\Ny{3}
\def\al{0.2}
\def\tw{0.15}
\def\tb{0.15}
\def\tg{0.1}
\def\ts{0.1}
\def\l{1}
\def\samples{50}
\begin{tikzpicture}[scale=5]
	\node (A) at (0,\tg/2){};
	\node (B) at (0,\tg){};
	\node (C) at (\tw/2,\tg/2){};
	\pgfmathsetmacro\x{\tw/2}
	\pgfmathsetmacro\y{\al*\l/2*cos(\x/\l*360)}
	\node (D) at (\x,\y){};
	\node (E) at ($(B)+(0,\ts)$){};
	\node (F) at ($(D)+(0,\ts)$){};
	\pgfmathsetmacro\x{\l/2-\tw/2}
	\pgfmathsetmacro\y{\al*\l/2*cos(\x/\l*360)}
	\node (G) at (\x,\y){};
	\node (H) at ($(G)+(0,\ts)$){};
	\pgfmathsetmacro\x{\l/2+\tw/2}
	\pgfmathsetmacro\y{\al*\l/2*cos(\x/\l*360)}
	\node (I) at (\x,\y){};
	\node (J) at ($(I)+(0,\ts)$){};
	\pgfmathsetmacro\x{\l-\tw/2}
	\node (K) at (\x,\tg/2){};
	\pgfmathsetmacro\x{\l-\tw/2}
	\pgfmathsetmacro\y{\al*\l/2*cos(\x/\l*360)}
	\node (L) at (\x,\y){};
	\node (M) at (\l,\tg/2){};
	\node (N) at (\l,\tg){};
	\node (O) at ($(L)+(0,\ts)$){};
	\node (P) at ($(N)+(0,\ts)$){};

	\node (Q) at ($(E)+(0,\tg)$){};
	\node (R) at ($(F)+(0,\tg)$){};
	\node (S) at ($(H)+(0,\tg)$){};
	\node (T) at ($(J)+(0,\tg)$){};
	\node (U) at ($(O)+(0,\tg)$){};
	\node (V) at ($(P)+(0,\tg)$){};

	\node (W) at ($(Q)+(0,\tb)$){};
	\node (X) at ($(R)+(0,\tb)$){};
	\node (Y) at ($(S)+(0,\tb)$){};
	\node (Z) at ($(T)+(0,\tb)$){};
	\node (AA) at ($(U)+(0,\tb)$){};
	\node (AB) at ($(V)+(0,\tb)$){};

	\node (AC) at (0,\tg+\tb+\ts+\al*\l/2+\tg/2){};
	\node (AD) at (\tw/2,\tg+\tb+\ts+\al*\l/2+\tg/2){};
	\node (AE) at (\x,\tg+\tb+\ts+\al*\l/2+\tg/2){};
	\node (AF) at (\l,\tg+\tb+\ts+\al*\l/2+\tg/2){};

	\draw[fill=gray2]  plot[domain=0:\tw/2,samples=\samples] (\x,{\al*\l/2*cos(\x/\l*360)}) -- (C.center) -- (A.center) -- cycle;
	\draw[fill=gray2]  plot[domain=\l-\tw/2:\l,samples=\samples] (\x,{\al*\l/2*cos(\x/\l*360)}) -- (M.center) -- (K.center) -- cycle;
	\draw[fill=col2]  plot[domain=0:\tw/2,samples=\samples] (\x,{\al*\l/2*cos(\x/\l*360)}) -- plot[domain=\tw/2:0,samples=\samples] (\x,{\al*\l/2*cos(\x/\l*360)+\ts}) -- cycle;
	\draw[fill=col2]  plot[domain=\tw/2:\l/2-\tw/2,samples=\samples] (\x,{\al*\l/2*cos(\x/\l*360)}) -- plot[domain=\l/2-\tw/2:\tw/2,samples=\samples] (\x,{\al*\l/2*cos(\x/\l*360)+\ts}) -- cycle;
	\draw[fill=col2]  plot[domain=\l/2-\tw/2:\l/2+\tw/2,samples=\samples] (\x,{\al*\l/2*cos(\x/\l*360)}) -- plot[domain=\l/2+\tw/2:\l/2-\tw/2,samples=\samples] (\x,{\al*\l/2*cos(\x/\l*360)+\ts}) -- cycle;
	\draw[fill=col2]  plot[domain=\l/2+\tw/2:\l-\tw/2,samples=\samples] (\x,{\al*\l/2*cos(\x/\l*360)}) -- plot[domain=\l-\tw/2:\l/2+\tw/2,samples=\samples] (\x,{\al*\l/2*cos(\x/\l*360)+\ts}) -- cycle;
	\draw[fill=col2]  plot[domain=\l-\tw/2:\l,samples=\samples] (\x,{\al*\l/2*cos(\x/\l*360)}) -- plot[domain=\l:\l-\tw/2,samples=\samples] (\x,{\al*\l/2*cos(\x/\l*360)+\ts}) -- cycle;

	\draw[fill=gray2]  plot[domain=\l/2-\tw/2:\l/2+\tw/2,samples=\samples] (\x,{\al*\l/2*cos(\x/\l*360)+\tg}) -- plot[domain=\l/2+\tw/2:\l/2-\tw/2,samples=\samples] (\x,{\al*\l/2*cos(\x/\l*360)+\ts+\tg}) -- cycle;

	\draw[fill=col1]  plot[domain=0:\tw/2,samples=\samples] (\x,{\al*\l/2*cos(\x/\l*360)+\ts+\tg}) -- plot[domain=\tw/2:0,samples=\samples] (\x,{\al*\l/2*cos(\x/\l*360)+\ts+\tg+\tb}) -- cycle;
	\draw[fill=col1]  plot[domain=\tw/2:\l/2-\tw/2,samples=\samples] (\x,{\al*\l/2*cos(\x/\l*360)+\ts+\tg}) -- plot[domain=\l/2-\tw/2:\tw/2,samples=\samples] (\x,{\al*\l/2*cos(\x/\l*360)+\ts+\tg+\tb}) -- cycle;
	\draw[fill=col1]  plot[domain=\l/2-\tw/2:\l/2+\tw/2,samples=\samples] (\x,{\al*\l/2*cos(\x/\l*360)+\ts+\tg}) -- plot[domain=\l/2+\tw/2:\l/2-\tw/2,samples=\samples] (\x,{\al*\l/2*cos(\x/\l*360)+\ts+\tg+\tb}) -- cycle;
	\draw[fill=col1]  plot[domain=\l/2+\tw/2:\l-\tw/2,samples=\samples] (\x,{\al*\l/2*cos(\x/\l*360)+\ts+\tg}) -- plot[domain=\l-\tw/2:\l/2+\tw/2,samples=\samples] (\x,{\al*\l/2*cos(\x/\l*360)+\ts+\tg+\tb}) -- cycle;
	\draw[fill=col1]  plot[domain=\l-\tw/2:\l,samples=\samples] (\x,{\al*\l/2*cos(\x/\l*360)+\ts+\tg}) -- plot[domain=\l:\l-\tw/2,samples=\samples] (\x,{\al*\l/2*cos(\x/\l*360)+\ts+\tg+\tb}) -- cycle;

	\draw[fill=gray2]  plot[domain=0:\tw/2,samples=\samples] (\x,{\al*\l/2*cos(\x/\l*360)+\ts+\tg+\tb}) -- (AD.center) -- (AC.center) -- cycle;
	\draw[fill=gray2]  plot[domain=\l-\tw/2:\l,samples=\samples] (\x,{\al*\l/2*cos(\x/\l*360)+\ts+\tg+\tb}) -- (AF.center) -- (AE.center) -- cycle;

	\path (H)--(S) node[midway] (HS){};
	\path (J)--(T) node[midway] (JT){};
	\draw[latex-latex] (HS.north) -- (JT.north) node[midway,below]{$t_w$};
	\draw[latex-latex] ($(N)+(\tw/4,0)$) -- ($(M)+(\tw/4,0)$) node[midway,right]{$\frac{t_g}{2}$};
	\draw[latex-latex] ($(P)+(\tw/4,0)$) -- ($(N)+(\tw/4,0)$) node[midway,right]{$t_s$};
	\draw[latex-latex] ($(V)+(\tw/4,0)$) -- ($(P)+(\tw/4,0)$) node[midway,right]{$t_g$};
	\draw[latex-latex] ($(AB)+(\tw/4,0)$) -- ($(V)+(\tw/4,0)$) node[midway,right]{$t_b$};
	\draw[latex-latex] ($(AF)+(\tw/4,0)$) -- ($(AB)+(\tw/4,0)$) node[midway,right]{$\frac{t_g}{2}$};

	\pgfmathsetmacro\y{-\al*\l/2}
	\draw[densely dotted] (-\tw/8,\y) -- (\l/2,\y);
	\draw[densely dotted] ($(B)-(\tw/8,0)$)  -- (B.center);
	\draw[latex-latex] (-\tw/8,\y)  -- ($(B)-(\tw/8,0)$) node[midway,left]{$a$};
	\draw[latex-latex] ($(A)-(\tw/4,0)$) -- ($(AC)-(\tw/4,0)$) node[midway,left]{$h$};

	\draw[latex-latex] (0,-\al*\l/2-\tg/2) -- (\l,-\al*\l/2-\tg/2) node[midway,below]{$l$};
\end{tikzpicture}

%% file: Sections/Conclusions.tex
In this paper, an Adaptive Parallel Arc Length Method (APALM) is presented. Contrary to existing parallel implementations of the Arc-Length Method (ALM), the present method provides within-branch parallelization, hence providing scalable parallelization independent of the physics of the problem, i.e., the number of branches. The method employs a multi-level approach, where parallel corrections are performed on solution intervals that have been initialised before. Given the sub-intervals provided by the serial computation, computations with finer arc lengths can be performed and evaluated using suitable error measures, marking intervals for further refinement when needed. Employing the multi-level approach, their implementations are discussed: the Adaptive Serial ALM (ASALM), the Adaptive Serial-Parallel ALM (ASPALM), and the APALM. The ASALM is a serial implementation, employing only the inherent adaptivity of the concept provided in this paper. The ASPALM is a two-stage approach, separating a serial initalization of the full equilibrium path from the parallel corrections. The APALM is a fully parallel implementation, where parallel corrections are performed as soon as the first path segments have been initialized. Conceptually, the APALM has a higher degree of parallelization since the workers are not idle until the full solution curve is obtained. Given a basic step function and distance computation, the present paper provides all algorithms necessary for implementing the APALM with manager-worker parallelization.\\

The implementation of the APALM is evaluated using three benchmark problems and an application example. The first problem involves the collapse of a composite shallow cylindrical shell. The second problem involves the collapse of a truncated conical rubber shell, and the third example involves the bifurcation problem of a strip subject to an in-plane load. Moreover, the method is applied to the modelling of a snapping metamaterial to investigate its performance on a larger-scale problem. In all examples, it can be observed that the APALM provides an accurate description of the reference solution, given a (sufficiently) coarse serial initialization of the curve. Through refinement, the APALM provides refinements (hence details in the solution), typically on sharp corners in the load-displacement diagrams. In addition, the bifurcation example also shows that the APALM is able to work within an exploration scheme for bifurcations. In all benchmark problems, the ASPALM and APALM have been used to evaluate the parallelization of the schemes. The natural separation of the serial and parallel stages of the ASPALM reveals the scalability of the parallel correction with respect to the number of workers, showing that the parallel correction can take only a fraction of the total computational time for a larger number of workers. Furthermore, the comparison between the ASPALM and the APALM shows that the full parallelization of the APALM provides a more efficient scheme than the two-stage approach of the ASPALM, as expected. The benchmarks and example also show that the APALM provides a full solution curve -- including adaptive refinements -- in the same computational time needed to compute only the initialisation of the ASPALM. This reveals the potential of the APALM: it can provide detailed solution paths without significantly increasing the computational time. \RevTHREE{The coarser the initial step size, the more arc-length intervals are computed during the parallel corrections of the method until a sufficient; hence, the higher the computational merit of the method to reach a desired level of detail. Moreover,} the scaling analyses also show that the benefits of the APALM are already achieved with a small number of workers, e.g., 8 workers, making the APALM interesting on a desktop scale. For larger clusters, the APALM can be employed using dynamic load balancing within OpenMP.\\

As the APALM enables parallelization in the arc-length domain, future applications of this method include quasi-static computations for solid and fluid dynamics, among other problems, especially those with a large number of load steps. Therefore, future works with this method include automatic exploration of solution spaces, e.g., following the work of \cite{Thies2021,Wouters2019}, or applications with large numbers of degrees of freedom, for instance with phase-field models for fracture mechanics \cite{Borden2014}. Other future work includes combining the APALM with a spatial refinement scheme to enable space-quasi-time refinements. MPI scalability and distribution of cores per worker are topics to investigate for different applications. \RevTWO{Another topic for further investigation is the convergence of the underlying arc-length method for large steps. Since a fewer number of initial intervals reduces the serial initialization time of the APALM, the parallel performance can be increased significantly when the initial step size is maximized. For example, the Mixed Integration Point (MIP) method increases the convergence of the ALM, allowing for larger step sizes. The performance of the MIP is demonstrated for isogeometric Kirchhoff--Love shells in \cite{Magisano2017,Magisano2017b,Leonetti2023,Leonetti2018a}}. \RevTHREE{Lastly, since the presented APALM is developed for path-independent problems, an extension to path-dependent problems is a natural direction for future research}.

%% file: paper.bbl
\begin{thebibliography}{38}
\expandafter\ifx\csname natexlab\endcsname\relax\def\natexlab#1{#1}\fi
\providecommand{\url}[1]{\texttt{#1}}
\providecommand{\href}[2]{#2}
\providecommand{\path}[1]{#1}
\providecommand{\DOIprefix}{doi:}
\providecommand{\ArXivprefix}{arXiv:}
\providecommand{\URLprefix}{URL: }
\providecommand{\Pubmedprefix}{pmid:}
\providecommand{\doi}[1]{\href{http://dx.doi.org/#1}{\path{#1}}}
\providecommand{\Pubmed}[1]{\href{pmid:#1}{\path{#1}}}
\providecommand{\bibinfo}[2]{#2}
\ifx\xfnm\relax \def\xfnm[#1]{\unskip,\space#1}\fi
\bibitem[{Rupp(2022)}]{Rupp2022}
\bibinfo{author}{K.~Rupp}, \bibinfo{title}{Microprocessor trend data},
  \bibinfo{howpublished}{\url{https://github.com/karlrupp/microprocessor-trend-data}},
  \bibinfo{year}{2022}.
\bibitem[{Gander(2015)}]{Gander2015}
\bibinfo{author}{M.~J. Gander},
\newblock \bibinfo{title}{{50 Years of Time Parallel Time Integration}}
  (\bibinfo{year}{2015}) \bibinfo{pages}{69--113}.
\bibitem[{Riks(1972)}]{Riks1972}
\bibinfo{author}{E.~Riks},
\newblock \bibinfo{title}{{The Application of Newton's Method to the Problem of
  Elastic Stability}},
\newblock \bibinfo{journal}{Journal of Applied Mechanics} \bibinfo{volume}{39}
  (\bibinfo{year}{1972}) \bibinfo{pages}{1060}.
\bibitem[{Crisfield(1981)}]{Crisfield1981}
\bibinfo{author}{M.~M. Crisfield},
\newblock \bibinfo{title}{{A Fast Incremental/Iterative Solution Procedure That
  Handles “Snap-Through”}},
\newblock in: \bibinfo{booktitle}{Computational Methods in Nonlinear Structural
  and Solid Mechanics}, \bibinfo{publisher}{Pergamon}, \bibinfo{year}{1981},
  pp. \bibinfo{pages}{55--62}.
\bibitem[{Wriggers et~al.(1988)Wriggers, Wagner, and Miehe}]{Wriggers1988}
\bibinfo{author}{P.~Wriggers}, \bibinfo{author}{W.~Wagner},
  \bibinfo{author}{C.~Miehe},
\newblock \bibinfo{title}{{A quadratically convergent procedure for the
  calculation of stability points in finite element analysis}},
\newblock \bibinfo{journal}{Computer Methods in Applied Mechanics and
  Engineering} \bibinfo{volume}{70} (\bibinfo{year}{1988})
  \bibinfo{pages}{329--347}.
\bibitem[{Pretti et~al.(2022)Pretti, Coombs, and Augarde}]{Pretti2022}
\bibinfo{author}{G.~Pretti}, \bibinfo{author}{W.~M. Coombs},
  \bibinfo{author}{C.~E. Augarde},
\newblock \bibinfo{title}{{A displacement-controlled arc-length solution
  scheme}},
\newblock \bibinfo{journal}{Computers \& Structures} \bibinfo{volume}{258}
  (\bibinfo{year}{2022}) \bibinfo{pages}{106674}.
\bibitem[{Kadapa(2021)}]{Kadapa2021}
\bibinfo{author}{C.~Kadapa},
\newblock \bibinfo{title}{{A simple extrapolated predictor for overcoming the
  starting and tracking issues in the arc-length method for nonlinear
  structural mechanics}},
\newblock \bibinfo{journal}{Engineering Structures} \bibinfo{volume}{234}
  (\bibinfo{year}{2021}) \bibinfo{pages}{111755}.
\bibitem[{Thies et~al.(2021)Thies, Wouters, Hennig, and Vanroose}]{Thies2021}
\bibinfo{author}{J.~Thies}, \bibinfo{author}{M.~Wouters},
  \bibinfo{author}{R.~S. Hennig}, \bibinfo{author}{W.~Vanroose},
\newblock \bibinfo{title}{{Towards Scalable Automatic Exploration of
  Bifurcation Diagrams for Large-Scale Applications}},
\newblock \bibinfo{journal}{Lecture Notes in Computational Science and
  Engineering} \bibinfo{volume}{139} (\bibinfo{year}{2021})
  \bibinfo{pages}{981--989}.
\bibitem[{Wouters and Vanroose(2019)}]{Wouters2019}
\bibinfo{author}{M.~Wouters}, \bibinfo{author}{W.~Vanroose},
\newblock \bibinfo{title}{{Automatic Exploration Techniques of Numerical
  Bifurcation Diagrams Illustrated by the Ginzburg--Landau Equation}},
\newblock \bibinfo{journal}{https://doi.org/10.1137/19M1248467}
  \bibinfo{volume}{18} (\bibinfo{year}{2019}) \bibinfo{pages}{2047--2098}.
\bibitem[{Lions et~al.(2001)Lions, Maday, and Turinici}]{Lions2001}
\bibinfo{author}{J.~L. Lions}, \bibinfo{author}{Y.~Maday},
  \bibinfo{author}{G.~Turinici},
\newblock \bibinfo{title}{{R{\'{e}}solution d'EDP par un sch{\'{e}}ma en temps
  << parar{\'{e}}el >>}},
\newblock \bibinfo{journal}{Comptes Rendus de l'Acad{\'{e}}mie des Sciences -
  Series I - Mathematics} \bibinfo{volume}{332} (\bibinfo{year}{2001})
  \bibinfo{pages}{661--668}.
\bibitem[{Falgout et~al.(2014)Falgout, Friedhoff, Kolev, MacLachlan, and
  Schroder}]{Falgout2014}
\bibinfo{author}{R.~D. Falgout}, \bibinfo{author}{S.~Friedhoff},
  \bibinfo{author}{T.~V. Kolev}, \bibinfo{author}{S.~P. MacLachlan},
  \bibinfo{author}{J.~B. Schroder},
\newblock \bibinfo{title}{{Parallel time integration with multigrid}},
\newblock \bibinfo{journal}{SIAM Journal on Scientific Computing}
  \bibinfo{volume}{36} (\bibinfo{year}{2014}) \bibinfo{pages}{C635--C661}.
\bibitem[{Cyr et~al.(2019)Cyr, G{\"u}nther, and Schroder}]{cyr2019multilevel}
\bibinfo{author}{E.~C. Cyr}, \bibinfo{author}{S.~G{\"u}nther},
  \bibinfo{author}{J.~B. Schroder},
\newblock \bibinfo{title}{Multilevel initialization for layer-parallel deep
  neural network training},
\newblock \bibinfo{journal}{arXiv preprint arXiv:1912.08974}
  (\bibinfo{year}{2019}).
\bibitem[{Hessenthaler et~al.(2021)Hessenthaler, Falgout, Schroder, de~Vecchi,
  Nordsletten, and R{\"o}hrle}]{hessenthaler2021time}
\bibinfo{author}{A.~Hessenthaler}, \bibinfo{author}{R.~D. Falgout},
  \bibinfo{author}{J.~B. Schroder}, \bibinfo{author}{A.~de~Vecchi},
  \bibinfo{author}{D.~Nordsletten}, \bibinfo{author}{O.~R{\"o}hrle},
\newblock \bibinfo{title}{Time-periodic steady-state solution of
  fluid-structure interaction and cardiac flow problems through
  multigrid-reduction-in-time},
\newblock \bibinfo{journal}{arXiv preprint arXiv:2105.00305}
  (\bibinfo{year}{2021}).
\bibitem[{Aruliah et~al.(2016)Aruliah, {Van Veen}, and Dubitski}]{Aruliah2016}
\bibinfo{author}{D.~A. Aruliah}, \bibinfo{author}{L.~{Van Veen}},
  \bibinfo{author}{A.~Dubitski},
\newblock \bibinfo{title}{{Algorithm 956: PAMPAC, a parallel adaptive method
  for pseudo-arclength continuation}},
\newblock \bibinfo{journal}{ACM Transactions on Mathematical Software}
  \bibinfo{volume}{42} (\bibinfo{year}{2016}).
\bibitem[{Rafsanjani et~al.(2015)Rafsanjani, Akbarzadeh, and
  Pasini}]{Rafsanjani2015}
\bibinfo{author}{A.~Rafsanjani}, \bibinfo{author}{A.~Akbarzadeh},
  \bibinfo{author}{D.~Pasini},
\newblock \bibinfo{title}{Snapping {{Mechanical Metamaterials}} under
  {{Tension}}},
\newblock \bibinfo{journal}{Advanced Materials} \bibinfo{volume}{27}
  (\bibinfo{year}{2015}) \bibinfo{pages}{5931--5935}.
\bibitem[{Ritto-Corr{\^{e}}a and Camotim(2008)}]{Ritto-Correa2008}
\bibinfo{author}{M.~Ritto-Corr{\^{e}}a}, \bibinfo{author}{D.~Camotim},
\newblock \bibinfo{title}{{On the arc-length and other quadratic control
  methods: Established, less known and new implementation procedures}},
\newblock \bibinfo{journal}{Computers and Structures} \bibinfo{volume}{86}
  (\bibinfo{year}{2008}) \bibinfo{pages}{1353--1368}.
\bibitem[{Ragon et~al.(2002)Ragon, G{\"{u}}rdal, and Watson}]{Ragon2002}
\bibinfo{author}{S.~A. Ragon}, \bibinfo{author}{Z.~G{\"{u}}rdal},
  \bibinfo{author}{L.~T. Watson},
\newblock \bibinfo{title}{{A comparison of three algorithms for tracing
  nonlinear equilibrium paths of structural systems}},
\newblock \bibinfo{journal}{International Journal of Solids and Structures}
  \bibinfo{volume}{39} (\bibinfo{year}{2002}) \bibinfo{pages}{689--698}.
\bibitem[{Schweizerhof and Wriggers(1986)}]{Schweizerhof1986}
\bibinfo{author}{K.~Schweizerhof}, \bibinfo{author}{P.~Wriggers},
\newblock \bibinfo{title}{{Consistent linearization for path following methods
  in nonlinear fe analysis}},
\newblock \bibinfo{journal}{Computer Methods in Applied Mechanics and
  Engineering} \bibinfo{volume}{59} (\bibinfo{year}{1986})
  \bibinfo{pages}{261--279}.
\bibitem[{Bellini and Chulya(1987)}]{Bellini1987}
\bibinfo{author}{P.~Bellini}, \bibinfo{author}{A.~Chulya},
\newblock \bibinfo{title}{{An improved automatic incremental algorithm for the
  efficient solution of nonlinear finite element equations}},
\newblock \bibinfo{journal}{Computers \& Structures} \bibinfo{volume}{26}
  (\bibinfo{year}{1987}) \bibinfo{pages}{99--110}.
\bibitem[{Carrera(1994)}]{Carrera1994}
\bibinfo{author}{E.~Carrera},
\newblock \bibinfo{title}{{A study on arc-length-type methods and their
  operation failures illustrated by a simple model}},
\newblock \bibinfo{journal}{Computers and Structures} \bibinfo{volume}{50}
  (\bibinfo{year}{1994}) \bibinfo{pages}{217--229}.
\bibitem[{Lam and Morley(1992)}]{Lam1992}
\bibinfo{author}{W.~F. Lam}, \bibinfo{author}{C.~T. Morley},
\newblock \bibinfo{title}{{Arc‐Length Method for Passing Limit Points in
  Structural Calculation}},
\newblock \bibinfo{journal}{Journal of Structural Engineering}
  \bibinfo{volume}{118} (\bibinfo{year}{1992}) \bibinfo{pages}{169--185}.
\bibitem[{Zhou and Murray(1995)}]{Zhou1995}
\bibinfo{author}{Z.~Zhou}, \bibinfo{author}{D.~Murray},
\newblock \bibinfo{title}{{An incremental solution technique for unstable
  equilibrium paths of shell structures}},
\newblock \bibinfo{journal}{Computers and Structures} \bibinfo{volume}{55}
  (\bibinfo{year}{1995}) \bibinfo{pages}{749--759}.
\bibitem[{Piegl and Tiller(1995)}]{Piegl1995}
\bibinfo{author}{L.~Piegl}, \bibinfo{author}{W.~Tiller}, \bibinfo{title}{{The
  NURBS Book}}, Monographs in Visual Communications,
  \bibinfo{publisher}{Springer Berlin Heidelberg}, \bibinfo{address}{Berlin,
  Heidelberg}, \bibinfo{year}{1995}. \DOIprefix\doi{10.1007/978-3-642-97385-7}.
\bibitem[{Kiendl et~al.(2009)Kiendl, Bletzinger, Linhard, and
  W{\"{u}}chner}]{Kiendl2009}
\bibinfo{author}{J.~Kiendl}, \bibinfo{author}{K.-U. Bletzinger},
  \bibinfo{author}{J.~Linhard}, \bibinfo{author}{R.~W{\"{u}}chner},
\newblock \bibinfo{title}{{Isogeometric shell analysis with Kirchhoff–Love
  elements}},
\newblock \bibinfo{journal}{Computer Methods in Applied Mechanics and
  Engineering} \bibinfo{volume}{198} (\bibinfo{year}{2009})
  \bibinfo{pages}{3902--3914}.
\bibitem[{Kiendl et~al.(2015)Kiendl, Hsu, Wu, and Reali}]{Kiendl2015}
\bibinfo{author}{J.~Kiendl}, \bibinfo{author}{M.-C. Hsu},
  \bibinfo{author}{M.~C. Wu}, \bibinfo{author}{A.~Reali},
\newblock \bibinfo{title}{{Isogeometric Kirchhoff–Love shell formulations for
  general hyperelastic materials}},
\newblock \bibinfo{journal}{Computer Methods in Applied Mechanics and
  Engineering} \bibinfo{volume}{291} (\bibinfo{year}{2015})
  \bibinfo{pages}{280--303}.
\bibitem[{Verhelst et~al.(2021)Verhelst, M{\"{o}}ller, {Den Besten},
  Mantzaflaris, and Kaminski}]{Verhelst2021}
\bibinfo{author}{H.~Verhelst}, \bibinfo{author}{M.~M{\"{o}}ller},
  \bibinfo{author}{J.~{Den Besten}}, \bibinfo{author}{A.~Mantzaflaris},
  \bibinfo{author}{M.~Kaminski},
\newblock \bibinfo{title}{{Stretch-Based Hyperelastic Material Formulations for
  Isogeometric Kirchhoff–Love Shells with Application to Wrinkling}},
\newblock \bibinfo{journal}{Computer-Aided Design} \bibinfo{volume}{139}
  (\bibinfo{year}{2021}) \bibinfo{pages}{103075}.
\bibitem[{{D}elft {H}igh {P}erformance {C}omputing~{C}entre
  ({DHPC})(2022)}]{DHPC2022}
\bibinfo{author}{{D}elft {H}igh {P}erformance {C}omputing~{C}entre ({DHPC})},
  \bibinfo{title}{{D}elft{B}lue {S}upercomputer ({P}hase 1)},
  \bibinfo{howpublished}{\url{https://www.tudelft.nl/dhpc/ark:/44463/DelftBluePhase1}},
  \bibinfo{year}{2022}.
\bibitem[{J{\"u}ttler et~al.(2014)J{\"u}ttler, Langer, Mantzaflaris, Moore, and
  Zulehner}]{Juttler2014}
\bibinfo{author}{B.~J{\"u}ttler}, \bibinfo{author}{U.~Langer},
  \bibinfo{author}{A.~Mantzaflaris}, \bibinfo{author}{S.~E. Moore},
  \bibinfo{author}{W.~Zulehner},
\newblock \bibinfo{title}{Geometry + {{Simulation Modules}}: {{Implementing
  Isogeometric Analysis}}},
\newblock \bibinfo{journal}{PAMM} \bibinfo{volume}{14} (\bibinfo{year}{2014})
  \bibinfo{pages}{961--962}.
\bibitem[{Leonetti et~al.(2019)Leonetti, Magisano, Madeo, Garcea, Kiendl, and
  Reali}]{Leonetti2019}
\bibinfo{author}{L.~Leonetti}, \bibinfo{author}{D.~Magisano},
  \bibinfo{author}{A.~Madeo}, \bibinfo{author}{G.~Garcea},
  \bibinfo{author}{J.~Kiendl}, \bibinfo{author}{A.~Reali},
\newblock \bibinfo{title}{{A simplified Kirchhoff–Love large deformation
  model for elastic shells and its effective isogeometric formulation}},
\newblock \bibinfo{journal}{Computer Methods in Applied Mechanics and
  Engineering} \bibinfo{volume}{354} (\bibinfo{year}{2019})
  \bibinfo{pages}{369--396}.
\bibitem[{Herrema et~al.(2019)Herrema, Johnson, Proserpio, Wu, Kiendl, and
  Hsu}]{Herrema2019}
\bibinfo{author}{A.~J. Herrema}, \bibinfo{author}{E.~L. Johnson},
  \bibinfo{author}{D.~Proserpio}, \bibinfo{author}{M.~C. Wu},
  \bibinfo{author}{J.~Kiendl}, \bibinfo{author}{M.-C. Hsu},
\newblock \bibinfo{title}{{Penalty coupling of non-matching isogeometric
  Kirchhoff–Love shell patches with application to composite wind turbine
  blades}},
\newblock \bibinfo{journal}{Computer Methods in Applied Mechanics and
  Engineering} \bibinfo{volume}{346} (\bibinfo{year}{2019})
  \bibinfo{pages}{810--840}.
\bibitem[{Başar and Itskov(1998)}]{Basar1998}
\bibinfo{author}{Y.~Başar}, \bibinfo{author}{M.~Itskov},
\newblock \bibinfo{title}{{Finite element formulation of the Ogden material
  model with application to ruber-like shells}},
\newblock \bibinfo{journal}{International Journal for Numerical Methods in
  Engineering}  (\bibinfo{year}{1998}).
\bibitem[{Pagani and Carrera(2018)}]{Pagani2018}
\bibinfo{author}{A.~Pagani}, \bibinfo{author}{E.~Carrera},
\newblock \bibinfo{title}{{Unified formulation of geometrically nonlinear
  refined beam theories}},
\newblock \bibinfo{journal}{Mechanics of Advanced Materials and Structures}
  \bibinfo{volume}{25} (\bibinfo{year}{2018}) \bibinfo{pages}{15--31}.
\bibitem[{Verhelst et~al.(2020)Verhelst, Moller, {Den Besten}, Vermolen, and
  Kaminski}]{Verhelst2019}
\bibinfo{author}{H.~M. Verhelst}, \bibinfo{author}{M.~Moller},
  \bibinfo{author}{J.~{Den Besten}}, \bibinfo{author}{F.~Vermolen},
  \bibinfo{author}{M.~Kaminski},
\newblock \bibinfo{title}{{Equilibrium Path Analysis Including Bifurcations
  with an Arc-Length Method Avoiding A Priori Perturbations}},
\newblock \bibinfo{journal}{Proceedings of ENUMATH2019 Conference}
  (\bibinfo{year}{2020}).
\bibitem[{Borden et~al.(2014)Borden, Hughes, Landis, and
  Verhoosel}]{Borden2014}
\bibinfo{author}{M.~J. Borden}, \bibinfo{author}{T.~J. Hughes},
  \bibinfo{author}{C.~M. Landis}, \bibinfo{author}{C.~V. Verhoosel},
\newblock \bibinfo{title}{{A higher-order phase-field model for brittle
  fracture: Formulation and analysis within the isogeometric analysis
  framework}},
\newblock \bibinfo{journal}{Computer Methods in Applied Mechanics and
  Engineering} \bibinfo{volume}{273} (\bibinfo{year}{2014})
  \bibinfo{pages}{100--118}.
\bibitem[{Magisano et~al.(2017{\natexlab{a}})Magisano, Leonetti, and
  Garcea}]{Magisano2017}
\bibinfo{author}{D.~Magisano}, \bibinfo{author}{L.~Leonetti},
  \bibinfo{author}{G.~Garcea},
\newblock \bibinfo{title}{{How to improve efficiency and robustness of the
  Newton method in geometrically non-linear structural problem discretized via
  displacement-based finite elements}},
\newblock \bibinfo{journal}{Comput. Methods Appl. Mech. Eng.}
  \bibinfo{volume}{313} (\bibinfo{year}{2017}{\natexlab{a}})
  \bibinfo{pages}{986--1005}.
\bibitem[{Magisano et~al.(2017{\natexlab{b}})Magisano, Leonetti, and
  Garcea}]{Magisano2017b}
\bibinfo{author}{D.~Magisano}, \bibinfo{author}{L.~Leonetti},
  \bibinfo{author}{G.~Garcea},
\newblock \bibinfo{title}{{Advantages of the mixed format in geometrically
  nonlinear analysis of beams and shells using solid finite elements}},
\newblock \bibinfo{journal}{Int. J. Numer. Methods Eng.} \bibinfo{volume}{109}
  (\bibinfo{year}{2017}{\natexlab{b}}) \bibinfo{pages}{1237--1262}.
\bibitem[{Leonetti and Kiendl(2023)}]{Leonetti2023}
\bibinfo{author}{L.~Leonetti}, \bibinfo{author}{J.~Kiendl},
\newblock \bibinfo{title}{{A mixed integration point (MIP) formulation for
  hyperelastic Kirchhoff–Love shells for nonlinear static and dynamic
  analysis}},
\newblock \bibinfo{journal}{Comput. Methods Appl. Mech. Eng.}
  \bibinfo{volume}{416} (\bibinfo{year}{2023}) \bibinfo{pages}{116325}.
\bibitem[{Leonetti et~al.(2018)Leonetti, Magisano, Liguori, and
  Garcea}]{Leonetti2018a}
\bibinfo{author}{L.~Leonetti}, \bibinfo{author}{D.~Magisano},
  \bibinfo{author}{F.~Liguori}, \bibinfo{author}{G.~Garcea},
\newblock \bibinfo{title}{{An isogeometric formulation of the Koiter's theory
  for buckling and initial post-buckling analysis of composite shells}},
\newblock \bibinfo{journal}{Comput. Methods Appl. Mech. Eng.}
  \bibinfo{volume}{337} (\bibinfo{year}{2018}) \bibinfo{pages}{387--410}.

\end{thebibliography}
